\documentclass[11pt]{article}
\usepackage{amssymb}
\usepackage{float}
\usepackage{graphicx,color}
\usepackage{amsmath}
\usepackage{amsbsy}
\usepackage{epsfig}
\usepackage{caption}
\usepackage{subcaption}

\newtheorem{lemma}{Lemma}[section]

\newtheorem{theorem}{Theorem}[section]

\def\proclaim#1{\par \bigskip\noindent {\bf #1}\bgroup\it\ }
\def\endproclaim{\egroup\par\bigskip}

\newbox\TempBox \newbox\TempBoxA

\def\text#1{\mbox{\rm #1}}
\def\overset#1#2{\stackrel{#1}{#2} }

\def\underwiggle 1{
\ifmmode\setbox\TempBox=\hbox{$ 1$}\else\setbox\TempBox=\hbox{ 1}\fi
\setbox\TempBoxA=\hbox to \wd\TempBox{\hss\char'176\hss}
\rlap{\copy\TempBox}\smash{\lower9pt\hbox{\copy\TempBoxA}} }

\newcommand{\be}{\begin{eqnarray}}
\newcommand{\ee}{\end{eqnarray}}
\newcommand{\by}{\begin{eqnarray*}}
\newcommand{\ey}{\end{eqnarray*}}
\newcommand{\bn}{\begin{enumerate}}
\newcommand{\en}{\end{enumerate}}
\newcommand{\bi}{\begin{itemize}}
\newcommand{\ei}{\end{itemize}}
\newcommand{\bds}{\begin{description}}
\newcommand{\eds}{\end{description}}
\newcommand{\bcen}{\begin{center}}
\newcommand{\ecen}{\end{center}}

\input{tcilatex}
\begin{document}

\title{\textbf{Asymptotic inferences for an AR(1) model with a change point:
stationary and nearly non-stationary cases\thanks{%
Tianxiao Pang and Danna Zhang's research was supported by the National
Natural Science Foundation of China (No. 11001236, No. 11071213 and No.
J1210038) and the Fundamental Research Funds for the Central Universities.%
\newline
\hspace*{5mm}E-mail addresses: txpang@zju.edu.cn, zhangdanna0507@gmail.com, chong2064@cuhk.edu.hk.}}\\
}
\date{}
\author{}
\maketitle


\begin{center}
\vskip -1cm {{\small Tianxiao Pang$^{1,}$\footnote[2]{Corresponding author.}, Danna Zhang$^{1,2}$ and Terence Tai-Leung Chong$^{3}$}}
{\small\centerline{$^{1}$Department of Mathematics, Yuquan Campus, Zhejiang University, Hangzhou 310027, P.R.China}}
{\small\centerline{$^{2}$Department of Statistics, University of Chicago, 5734 S. University Ave., Chicago, Illinois 60637, USA}}
{\small\centerline{$^{3}$Department of Economics and Institute of Global Economics and Finance, }}
{\small\centerline{The Chinese University of Hong Kong, Hong Kong}}
\end{center}

\bigskip \textbf{Abstract.} This paper examines the asymptotic inference for
AR(1) models with a possible structural break in the AR parameter $\beta $\
near the unity at an unknown time $k_{0}$. Consider the model $y_{t}=\beta
_{1}y_{t-1}I\{t\leq k_{0}\}+\beta _{2}y_{t-1}I\{t>k_{0}\}+\varepsilon
_{t},~t=1,2,\cdots ,T,$ where $I\{\cdot \}$ denotes the indicator function.
We examine two cases: Case (I) $|\beta _{1}|<1,\beta _{2}=\beta _{2T}=1-c/T$%
; and case (II)~$\beta _{1}=\beta _{1T}=1-c/T,|\beta _{2}|<1$, where $c$\ is
a fixed constant, and $\{\varepsilon _{t},t\geq 1\}$\ is a sequence of
i.i.d. random variables which are in the domain of attraction of the normal
law with zero means and possibly infinite variances. We derive the limiting
distributions of the least squares estimators of $\beta _{1}$ and $\beta
_{2} $, and that of the break-point estimator for shrinking break for the
aforementioned cases. Monte Carlo simulations are conducted to demonstrate
the finite sample properties of the estimators. Our theoretical results are supported by
Monte Carlo simulations.\newline



\section{Introduction}

\setcounter{equation}{0}

The change-point problem has received considerable attention in the literature over the past three decades (Mankiw and Miron, 1986; Mankiw, Miron and Weil, 1987; Hansen, 1992; Chong, 2001). This paper extends the work of Chong (2001), who studies an AR(1) model with a structural break in the AR parameter $\beta $ at an unknown time $k_{0}$. We consider the following model:
\begin{equation}
y_{t}=\beta _{1}y_{t-1}I\{t\leq k_{0}\}+\beta
_{2}y_{t-1}I\{t>k_{0}\}+\varepsilon _{t},~t=1,2,\cdots ,T,  \label{model1}
\end{equation}%
where $I\{\cdot \}$ denotes the indicator function and $\{\varepsilon
_{t},t\geq 1\}$ is a sequence of i.i.d. random variables. Under some
regularity conditions that $E\varepsilon _{t}^{4}<\infty $\ and $%
Ey_{0}^{2}<\infty $, Chong (2001) proves the consistency and derives the
limiting distributions of the least squares estimators of $\beta _{1},\beta
_{2}$\ and $\tau _{0}$\ for three cases: (1) $|\beta _{1}|<1$ and $|\beta
_{2}|<1$; (2) $|\beta _{1}|<1$ and $\beta _{2}=1$; (3) $\beta _{1}=1$ and $%
|\beta _{2}|<1$.

In the present paper, we focus on Model (\ref{model1}) where one of
the pre-shift and post-shift AR parameters is less than one in absolute
value while the other is local to unity. This case is omitted in Chong
(2001). Specifically, we focus on the following two cases: (I) $|\beta
_{1}|<1$, $\beta _{2}=\beta _{2T}=1-c/T$; (II) $\beta _{1}=\beta
_{1T}=1-c/T,|\beta _{2}|<1$, where $c$ is a fixed constant. The case of
local to unity in AR(1) model was first independently studied by Chan and
Wei (1987) and Phillips (1987). Their studies bridge the gap between
stationary AR(1) model and unit root model. Moreover, since heavy-tailed
distributions, such as Student's $t$ distribution with degrees of freedom 2
and Pareto distribution with index 2, are commonly found in insurance,
econometrics and other literature, it is more appropriate to impose
weaker moment conditions on the $\varepsilon _{t}$'s and $y_{0}$ than those
in Chong (2001). The primary contribution of this paper is to derive the
consistency and limiting distributions of the least squares estimators of $%
\beta _{1}$, $\beta _{2}$ and the estimator of $\tau _{0}$ under a more
general setting.

Throughout the rest of the present paper, we shall focus on the random
variables which are in the domain of attraction of the normal law (DAN),
which is an important sub-class of heavy-tailed random variables.  A
sequence of i.i.d. random variables $\{X_{i},i\geq 1\}$ belongs to DAN if there exist two constant sequences $\{A_{n},n\geq 1\}$ and $\{B_{n},n\geq 1\}$ such that $Z_{n}:=B_{n}^{-1}(X_{1}+\cdots +X_{n})-A_{n}$ converges to a standard normal
random variable in distribution (Feller, 1971), where $B_{n}$ takes the form $\sqrt{n}h(n)$ and $h(n)$ is a slowly
varying function at infinity. We make the following assumptions:

\begin{itemize}
\item C1:~$\{\varepsilon _{t},t\geq 1\}$ is a sequence of i.i.d. random
variables which are in the domain of attraction of the normal law with zero
means and possibly infinite variances.

\item C2:~$y_0$ is an arbitrary random variable such that $y_0=o_p(\sqrt{T})$%
, where $T$ is the sample size.

\item C3:~$\tau _{0}\in \lbrack \underline{\tau },\overline{\tau }]\subset
(0,1)$.
\end{itemize}

\medskip

\noindent \textbf{Remark 1}.~ Assumption C2 is a weak initial condition. It
not only allows $y_{0}$ to be a finite random variable, but also allows it
to be a random variable of order smaller than $\sqrt{T}$ in probability.
\newline

For any given $\tau $, the ordinary least squares estimators of parameters $%
\beta _{1}$ and $\beta _{2}$ are given by
\begin{equation*}
\hat{\beta}_{1}(\tau )=\frac{\sum_{t=1}^{[\tau T]}y_{t}y_{t-1}}{%
\sum_{t=1}^{[\tau T]}y_{t-1}^{2}},~~\hat{\beta}_{2}(\tau )=\frac{%
\sum_{t=[\tau T]+1}^{T}y_{t}y_{t-1}}{\sum_{t=[\tau T]+1}^{T}y_{t-1}^{2}},
\end{equation*}%
The symbol $[a]$ denotes the integer part of $a$ and the change-point
estimator satisfies
\begin{equation*}
\hat{\tau}_{T}=\mathop{\arg\min}_{\tau \in (0,1)}RSS_{T}(\tau ),
\end{equation*}%
where
\begin{equation*}
RSS_{T}(\tau )=\sum_{t=1}^{[\tau T]}\Big(y_{t}-\hat{\beta}_{1}(\tau )y_{t-1}%
\Big)^{2}+\sum_{t=[\tau T]+1}^{T}\Big(y_{t}-\hat{\beta}_{2}(\tau )y_{t-1}%
\Big)^{2}.
\end{equation*}

We introduce some notations before presenting our main results. Let $%
W_{1}(\cdot )$\ and $W_{2}(\cdot )$\ be two independent Brownian motions
defined on the non-negative half real $R_{+}$; $W(\cdot )$ and $\overline{W}%
(\cdot )$ be two independent Brownian motions defined on $[0,1]$ and $R_{+}$
respectively; "$\Rightarrow $" signifies the weak convergence of the
associated probability measures; "$\overset{p}{\rightarrow }$" represents
convergence in probability; "$\overset{d}{=}$" denotes identical in
distribution. Let $C$ be a finite constant. The
limits in this paper are all taken as $T\rightarrow \infty $ unless
specified otherwise.

Under assumptions C1-C3, we have
{\theorem\label{thm}
In Model (\ref{model1}), if $|\beta _{1}|<1$, $\beta _{2}=\beta _{2T}=1-c/T$, where $c$ is a
fixed constant, and the assumptions C1-C3 are satisfied, then the estimators
$\hat{\tau}_{T},\hat{\beta}_{1}(\hat{\tau}_{T})$ and $\hat{\beta}_{2}(\hat{%
\tau}_{T})$ are all consistent, and
\begin{equation}
\left\{
\begin{array}{ll}
|\hat{\tau}_{T}-\tau _{0}|=O_{p}(1/T), &  \\
\sqrt{T}(\hat{\beta}_{1}(\hat{\tau}_{T})-\beta _{1})\Rightarrow N(0,(1-\beta
_{1}^{2})/\tau _{0}), &  \\
\displaystyle T(\hat{\beta}_{2}(\hat{\tau}_{T})-\beta _{2})\Rightarrow \frac{%
\frac{1}{2}F^{2}(W,c,\tau _{0},1)+c\int_{\tau
_{0}}^{1}e^{2c(1-t)}F^{2}(W,c,\tau _{0},t)dt-\frac{1}{2}(1-\tau _{0})}{%
\int_{\tau _{0}}^{1}e^{2c(1-t)}F^{2}(W,c,\tau _{0},t)dt}, &
\end{array}%
\right.   \label{result1}
\end{equation}%
where
\begin{equation*}
F(W,c,\tau _{0},t)=e^{-c(1-t)}(W(t)-W(\tau _{0}))-c\int_{\tau
_{0}}^{t}e^{-c(1-s)}(W(s)-W(\tau _{0}))ds.
\end{equation*}%
If we also let $\beta _{1T}$ be a sequence of $\beta _{1}$ such that $|\beta
_{2T}-\beta _{1T}|\rightarrow 0$ and $T(\beta _{2T}-\beta _{1T})\rightarrow
\infty $, then the limiting distribution of $\hat{\tau}_{T}$ is given by
\begin{equation*}
(\beta _{2T}-\beta _{1T})T(\hat{\tau}_{T}-\tau _{0})\Rightarrow %
\mathop{\arg\max}_{\nu \in R}\bigg\{\frac{C^{\ast }(\nu )}{B_{a}(\frac{1}{2})%
}-\frac{|\nu |}{2}\bigg\},
\end{equation*}%
where $B_{a}(\frac{1}{2})$ is generated by $\int_{0}^{\infty }\exp {(-s)}%
dW_{1}(s)$ and $C^{\ast }(\nu )$ is defined to be $C^{\ast }(\nu
)=W_{1}(-\nu )$ for $\nu \leq 0$ and
\begin{eqnarray*}
C^{\ast }(\nu ) &=&-I(W_{2},c,\tau _{0},\nu )-\int_{0}^{\nu }\frac{%
I(W_{2},c,\tau _{0},t)}{B_{a}(\frac{1}{2})}dI(W_{2},c,\tau _{0},t) \\
&&-\int_{0}^{\nu }\Big(\frac{I(W_{2},c,\tau _{0},t)}{2B_{a}(\frac{1}{2})}+1%
\Big)I(W_{2},c,\tau _{0},t)dt
\end{eqnarray*}%
for $\nu >0$ with
\begin{equation*}
I(W_{2},c,\tau _{0},t)=W_{2}(\tau _{0}+t)-W_{2}(\tau _{0})-c\int_{\tau
_{0}}^{\tau _{0}+t}e^{-c(\tau _{0}+t-s)}(W_{2}(s)-W_{2}(\tau _{0}))ds.
\end{equation*}%
}

{\theorem\label{thm2}  In Model (\ref{model1}), if $\beta_{1}=\beta_{1T}=1-c/T$, where $c$ is a fixed constant and $|\beta_{2}|<1$%
, and the assumptions C1-C3 are satisfied, then the estimators $\hat{\tau}_{T},\hat{\beta}_{1}(\hat{\tau}_{T})$ and $\hat{\beta}_{2}(\hat{\tau}_{T})$
are all consistent and
\begin{equation}
\left\{
\begin{array}{ll}
\mathbb{P}(\hat{k}\neq k_{0})\rightarrow 0, &  \\
\displaystyle T(\hat{\beta}_{1}(\hat{\tau}_{T})-\beta _{1})\Rightarrow \frac{%
\frac{1}{2}e^{2c(1-\tau _{0})}G^{2}(W,c,\tau _{0})+c\int_{0}^{\tau
_{0}}e^{2c(1-t)}G^{2}(W,c,t)dt-\frac{\tau _{0}}{2}}{\int_{0}^{\tau
_{0}}e^{2c(1-t)}G^{2}(W,c,t)dt}, &  \\
\displaystyle\sqrt{T}(\hat{\beta}_{2}(\hat{\tau}_{T})-\beta _{2})\Rightarrow
\frac{\sqrt{1-\beta _{2}^{2}}\cdot \overline{W}(B(c,\tau _{0}))}{1-\tau
_{0}+e^{2c(1-\tau _{0})}G^{2}(W,c,\tau _{0})}, &
\end{array}%
\right.   \label{result2}
\end{equation}%
where
\begin{equation*}
G(W,c,t)=e^{-c(1-t)}W(t)-c\int_{0}^{t}e^{-c(1-s)}W(s)ds
\end{equation*}%
and
\begin{equation*}
B(c,\tau _{0})=(1-e^{-2c\tau _{0}})/(2c)+1-\tau _{0}.
\end{equation*}%
Suppose we also let $\beta _{2T}$ be a sequence of $\beta _{2}$ such that $%
\sqrt{T}(\beta _{2T}-\beta _{1T})\rightarrow 0$ and $T^{3/4}(\beta
_{1T}-\beta _{2T})\rightarrow \infty $, then the limiting distribution of $%
\hat{\tau}_{T}$ is given by
\begin{equation*}
(\beta _{2T}-\beta _{1T})^{2}T^{2}(\hat{\tau}_{T}-\tau _{0})\Rightarrow %
\mathop{\arg\max}_{\nu \in R}\bigg\{\frac{B^{\ast }(\nu )}{e^{c(1-\tau
_{0})}G(W_{1},c,\tau _{0})}-\frac{|\nu |}{2}\bigg\},
\end{equation*}%
where $B^{\ast }(\nu )$ is a two-sided Brownian motion on $R$ defined to be $B^{\ast }(\nu )=W_{1}(-\nu )$ for $\nu \leq 0$ and $B^{\ast }(\nu)=W_{2}(\nu )$ for $\nu >0$. }\bigskip \bigskip

\noindent \textbf{Remark 2}. In Theorem \ref{thm}, letting $c=0$, it is
clear that
\begin{eqnarray*}
&&\frac{\frac{1}{2}F^{2}(W,c,\tau _{0},1)+c\int_{\tau
_{0}}^{1}e^{2c(1-t)}F^{2}(W,c,\tau _{0},t)dt-\frac{1}{2}(1-\tau _{0})}{%
\int_{\tau _{0}}^{1}e^{2c(1-t)}F^{2}(W,c,\tau _{0},t)dt}\Big|_{c=0} \\
&=&\frac{\frac{1}{2}(W(1)-W(\tau _{0}))^{2}-\frac{1}{2}(1-\tau _{0})}{%
\int_{\tau _{0}}^{1}(W(t)-W(\tau _{0}))^{2}dt} \\
&\overset{d}{=} &\frac{W^{2}(1)-1}{2(1-\tau _{0})\int_{0}^{1}W^{2}(t)dt}
\end{eqnarray*}%
and for $\nu >0$
\begin{eqnarray*}
&&\Big\{-I(W_{2},c,\tau _{0},\nu )-\int_{0}^{\nu }\frac{I(W_{2},c,\tau
_{0},t)}{B_{a}(\frac{1}{2})}dI(W_{2},c,\tau _{0},t) \\
&&~~~-\int_{0}^{\nu }\Big(\frac{I(W_{2},c,\tau _{0},t)}{2B_{a}(\frac{1}{2})}%
+1\Big)I(W_{2},c,\tau _{0},t)dt\Big\}\Big|_{c=0} \\
&=&-(W_{2}(\tau _{0}+\nu )-W_{2}(\tau _{0}))-\int_{0}^{\nu }\frac{W_{2}(\tau
_{0}+t)-W_{2}(\tau _{0})}{B_{a}(\frac{1}{2})}d(W_{2}(\tau _{0}+t)-W_{2}(\tau
_{0})) \\
&&~~-\int_{0}^{\nu }\Big(\frac{W_{2}(\tau _{0}+t)-W_{2}(\tau _{0})}{2B_{a}(%
\frac{1}{2})}+1\Big)(W_{2}(\tau _{0}+t)-W_{2}(\tau _{0}))dt \\
&\overset{d}{=} &-W_{2}(\nu )-\int_{0}^{\nu }\frac{W_{2}(t)}{B_{a}(\frac{1}{2%
})}dW_{2}(t)-\int_{0}^{\nu }\Big(\frac{W_{2}(t)}{2B_{a}(\frac{1}{2})}+1\Big)%
W_{2}(t)dt.
\end{eqnarray*}%
The above two expressions coincide with the third term of (15) and $C^{\ast
}(\nu )$ with $\nu >0$ in Chong (2001), respectively. Hence, our Theorem \ref%
{thm} is reduced to Theorem 3 in Chong (2001) by taking $c=0$.

Similarly, letting $c=0$ in Theorem \ref{thm2}, we have
\begin{equation*}
\frac{\frac{1}{2}e^{2c(1-\tau _{0})}G^{2}(W,c,\tau _{0})+c\int_{0}^{\tau
_{0}}e^{2c(1-t)}G^{2}(W,c,t)dt-\frac{\tau _{0}}{2}}{\int_{0}^{\tau
_{0}}e^{2c(1-t)}G^{2}(W,c,t)dt}\Big|_{c=0}=\frac{W^{2}(\tau _{0})-\tau _{0}}{%
2\int_{0}^{\tau _{0}}W^{2}(t)dt},
\end{equation*}%
\begin{equation*}
\frac{\sqrt{1-\beta _{2}^{2}}\cdot \overline{W}(B(c,\tau _{0}))}{1-\tau
_{0}+e^{2c(1-\tau _{0})}G^{2}(W,c,\tau _{0})}\Big|_{c=0}=\frac{\sqrt{1-\beta
_{2}^{2}}\overline{W}(1)}{1-\tau _{0}+W^{2}(\tau _{0})}
\end{equation*}%
and
\begin{equation*}
e^{c(1-\tau _{0})}G(W_{1},c,\tau _{0})\big|_{c=0}=W_{1}(\tau _{0}),
\end{equation*}%
indicating that Theorem \ref{thm2} is reduced to Theorem 4 in Chong
(2001) when $c=0$.

Note that the assumptions on the $\varepsilon _{t}$'s and $y_{0}$
are weaker than those in Chong (2001).\newline

\noindent \textbf{Remark 3}. The limiting distributions of $\hat{\beta}_{1}(\hat{\tau}_{T})$ and $\hat{\beta}_{2}(\hat{\tau}_{T})$
in Theorem \ref{thm2} could be simplified
if assumption C2 is more specific. For example, if the initial
value $y_{0}$ is defined as $y_{0}=y_{T,0}=\sum_{j=0}^{\infty }\rho
_{T}^{j}\varepsilon _{-j}$ with $\rho _{T}$ satisfying $T(1-\rho
_{T})=h_{T}\rightarrow 0$ and $\{\varepsilon _{-j},j\geq 0\}$ being a
sequence of i.i.d. random variables sharing the same distribution with $\varepsilon_1$, then similar
arguments of Lemma 3 in Andrews and Guggenberger (2008) will lead to $\sqrt{%
2h_{T}}y_{0}/\sqrt{Tl(\eta _{T})}\Rightarrow N(0,1)$, where the definitions of $\eta_T$ and the function $l(\cdot)$ can be found at the beginning of Section \ref{sec3}. Since $y_{0}$ dominates the asymptotic distribution of $\hat{%
\beta}_{1}(\hat{\tau}_{T})$, we have
\begin{eqnarray*}
&&\frac{2h_{T}}{Tl(\eta _{T})}\sum_{t=1}^{[\tau _{0}T]}y_{t-1}\varepsilon
_{t}\overset{p}{\rightarrow }0, \\
&&\frac{2h_{T}}{T^{2}l(\eta _{T})}\sum_{t=1}^{[\tau
_{0}T]}y_{t-1}^{2}\Rightarrow \frac{1-e^{-2c\tau _{0}}}{2c}W^{2}(1).
\end{eqnarray*}%
Consequently, we have
\begin{equation*}
T(\hat{\beta}_{1}(\hat{\tau}_{T})-\beta _{1})\overset{p}{\rightarrow }0.
\end{equation*}%
Similarly, from the proof of Lemma \ref{lem4}, it can be shown that
\begin{equation*}
\sqrt{\frac{T}{2h_{T}}}(\hat{\beta}_{1}(\hat{\tau}_{T})-\beta
_{1})\Rightarrow (1-\beta _{2}^{2})e^{c\tau _{0}}\pi (\beta _{2})/W(1)
\end{equation*}%
if the stationary distribution (denoted by $\pi (\beta _{2})$) of the AR(1) process: $y_{t}=\beta _{2}y_{t-1}+\varepsilon _{t}/\sqrt{l(\eta _{T})}$ with $%
y_{0}=0$ for $t=1,\cdots ,T-[\tau _{0}T]$, exists. Note that $\pi(\beta _{2})$ and $W(1)$ are independent.\newline

\noindent \textbf{Remark 4}. Chong (2001) proves that $|\hat{\tau}_{T}-\tau
_{0}|=O_{p}(1/T)$ in the case of $|\beta _{1}|<1$ and $\beta _{2}=1$, while $%
\mathbb{P}(\hat{k}\neq k_{0})\rightarrow 0$ in the case of $\beta _{1}=1$
and $|\beta _{2}|<1$. This result also holds in the present paper. Note that
the result about the estimator of $k_{0}$ in Theorem \ref{thm2} is stronger
than that in Theorem \ref{thm}. This is because the signal
from the regressor $y_{t-1}$ when the serial correlation coefficient is $%
1-c/T$ is stronger than that from the regressor $y_{t-1}$ when the serial
correlation coefficient is a fixed constant smaller than one in absolute
value (as implied by the faster convergence rate of $\hat{\beta}_{2}(\hat{%
\tau}_{T})$ in Theorem \ref{thm} and the faster convergence rate of $\hat{%
\beta}_{1}(\hat{\tau}_{T})$ in Theorem \ref{thm2}), meanwhile, the signal
from the regressor $y_{t-1}$ under the situation of $(\beta _{1},\beta
_{2})=(1-c/T,c_{0})$ is stronger than that under the situation of $(\beta
_{1},\beta _{2})=(c_{0},1-c/T)$, where $c_{0}$ is fixed and $|c_{0}|<1$.
\newline

\noindent \textbf{Remark 5}. The statistical inference on the least squares
estimators of $\beta _{1}$, $\beta _{2}$\ and $\hat{\tau}_{T}$\ for the
following cases: (I) $\beta _{1}=\beta _{1T}=1-c/T$, $\beta _{2}=1$; (II) $%
\beta _{1}=1,\beta _{2}=\beta _{2T}=1-c/T$\ are much more complicated and
would be left for future research.\newline

The rest of the paper is organized as follows: Section 2 presents the
simulation results for the finite sample properties of the estimators in
Theorems \ref{thm} and \ref{thm2}. Section 3 states some useful lemmas and
provides the proof for Theorem \ref{thm}. Section 4 provides the proof for
Theorem \ref{thm2}. \bigskip

\section{Simulations}

We perform the following experiments to see how well our asymptotic results
match the finite-sample properties of the estimators. In all
experiments, the sample size is set at $T=200$ and the number of
replications is set at $N=20,000$; $\{y_{t}\}_{t=1}^{T}$ is generated from
Model (\ref{model1}); $y_{0}$ has the following probability density
function
\begin{equation*}
f(x)=%
\begin{cases}
0 & \mathrm{if}~x\leq -2, \\
\frac{3}{2(x+3)^{5/2}} & \mathrm{if}~x>-2.%
\end{cases}%
\end{equation*}%
Note that $E|y_{0}|<\infty $ and $E|y_{0}|^{3/2+\delta }=\infty $ for any $%
\delta \geq 0$. Note also that assumption C2 holds. The true change point is
set at $\tau _{0}=0.3$ and $0.5$. For the constant $c$ and the distribution
of $\varepsilon _{t}$'s, we consider the following numerical setup:
\begin{equation}
\left\{
\begin{array}{l}
c=1, \\
\varepsilon _{t}\in \left\{ t(3),\ \,t(2)\right\} ,%
\end{array}%
\right.   \label{numerical-setup}
\end{equation}%
where $t(3)$ and $t(2)$ denote the student-$t$ random variables with degrees of
freedom 3 and 2 respectively. It is easy to verify that $t(3)$ and $t(2)$
are both in the domain of attraction of the normal law, and that $t(3)$\ has
finite variance but infinite fourth moment, while $t(2)$ has infinite
variance. When the AR parameter which depends on the sample size $T$ and the
constant $c$ is determined, the values of the fixed AR parameter are
set to $0.5$, $0.75$ and $0.8$. Then, both the second and
third results in (\ref{result1}) and (\ref{result2}) are conducted under (\ref{numerical-setup}).


First, we conduct experiments to verify the second and the third results in (%
\ref{result1}) under (\ref{numerical-setup}), which predict that the
finite-sample distribution of $\hat{\beta}_{1}(\hat{\tau}_{T})$ is
approximately normal, whereas $\hat{\beta}_{2}(\hat{\tau}_{T})$ appears to
have a Dickey-Fuller distribution. Figure 1-Figure 12 agree with our results.

Second, we conduct experiments to verify the second and the third results in
(\ref{result2}) under (\ref{numerical-setup}), which predict that
the finite-sample distribution of $\hat{\beta}_{1}(\hat{\tau}_{T})$
is approximately the Dickey-Fuller type, whereas $\hat{\beta}_{2}(\hat{\tau}%
_{T})$ will have a normal distribution. Figure 13-Figure 24 also
agree with our results.

In all of the above experiments, we have the following
observations: (I) the performance for the simulations on $\hat{\beta}_{1}(\hat{\tau}_{T})$ is better for the case where $\tau_{0}=0.5$
compared to the case where $\tau_{0}=0.3$, since more data are used to generate the first subsample. Analogously, the performance
for the simulations on $\hat{\beta}_{2}(\hat{\tau}_{T})$ is better for the case where $\tau _{0}=0.3$ compared to the case where $\tau_{0}=0.5$.~ (II) the performance is better for the case where $\{\varepsilon _{t}\}_{t=1}^{T}\sim t(3)$\ compared to the case where $\{\varepsilon _{t}\}_{t=1}^{T}\sim t(2)$ under (\ref{numerical-setup}). This is not surprising since $t(2)$ is a heavy-tailed distribution. We have also conducted the experiments for larger sample size when $\{\varepsilon_{t}\}_{t=1}^{T}\sim t(2)$, the performance does not improve much. Given the results, one should be cautious when conducting statistical inference under heavy-tailed innovations such as $t(2)$. The experiments when $c=-1$\ are also studied. However, since the results are very similar to the case where $c=1$, we do not report those simulations here to conserve space.

In the following figures, we let $c=1$. The solid line shows the finite sample distribution while the dashed line shows the asymptotic distribution.\\
\bigskip

\begin{figure}[H]
\label{fig1}
\begin{subfigure}[b]{0.5\textwidth}
\vspace{0.7cm}\centering
\includegraphics[trim=5cm 5cm 2cm 4cm, width=0.55\textwidth]{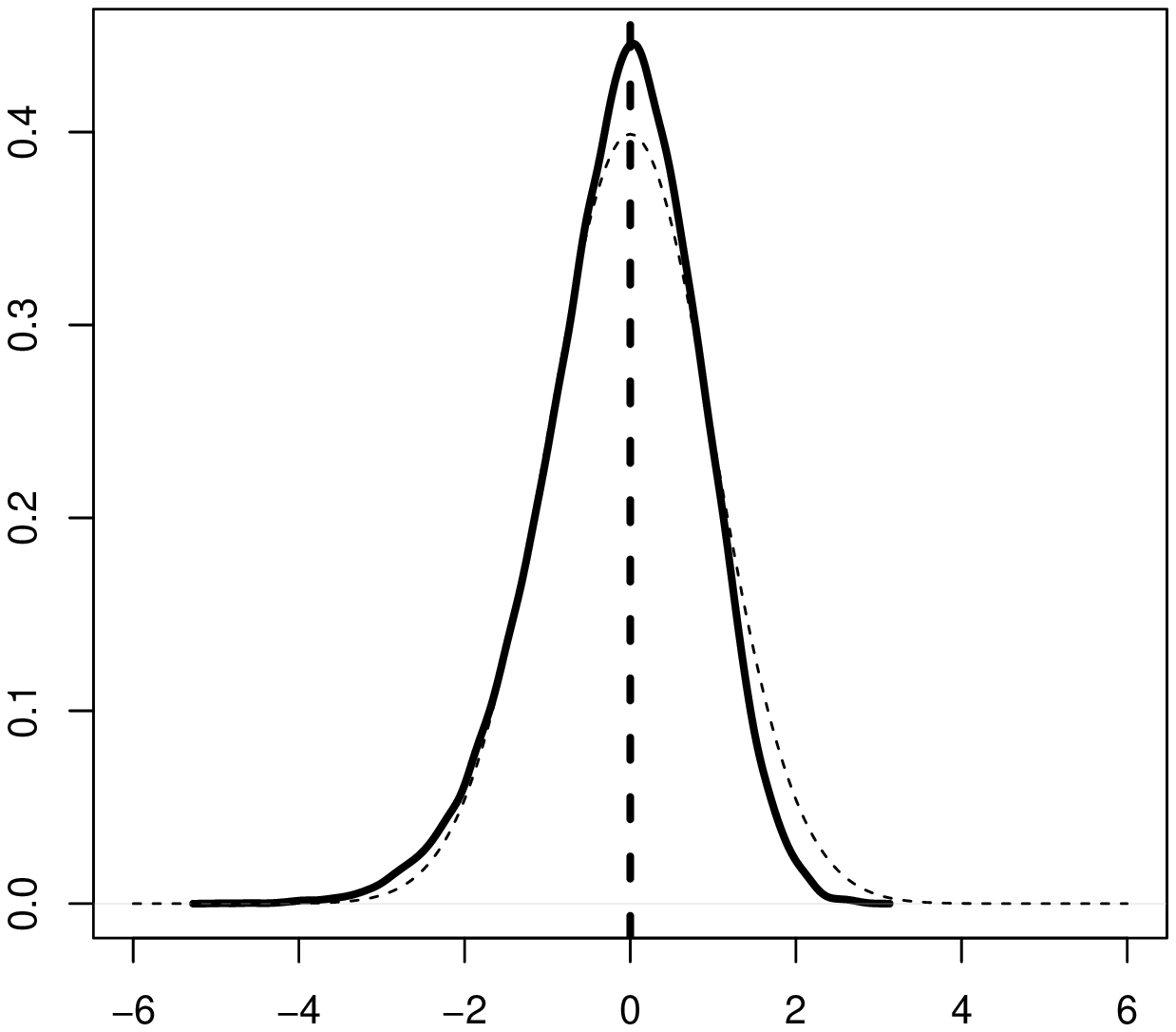}
\end{subfigure}\hspace{0.7cm}
\begin{subfigure}[b]{0.5\textwidth}
\vspace{0.7cm}\centering
\includegraphics[trim=5cm 5cm 2cm 4cm, width=0.55\textwidth]{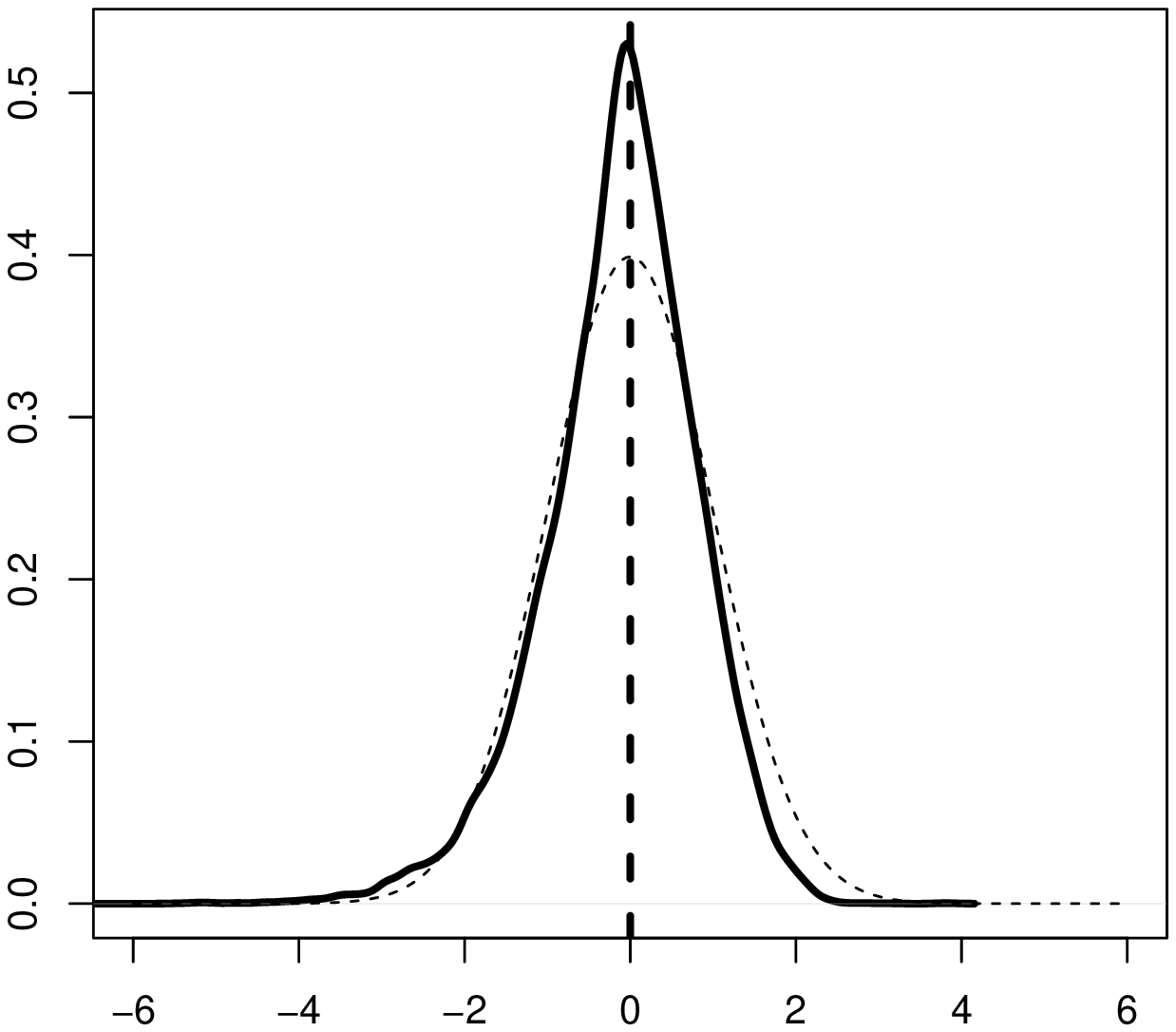}
\end{subfigure}\vspace{1.5cm}
\caption{{\small Distribution of $\sqrt{T}(\hat{\beta}_1(\hat{\tau}_T)-\beta_1)$ when $\tau_0=0.3$, $\beta_1=0.5$, $\beta_2=1-1/T$. Left: $\{\varepsilon_t\}_{t=1}^T\sim t(3)$; Right: $\{\varepsilon_t\}_{t=1}^T\sim t(2)$.}}
\end{figure}

\begin{figure}[H]
\label{fig1}
\begin{subfigure}[b]{0.5\textwidth}
\vspace{0.7cm}\centering
\includegraphics[trim=5cm 5cm 2cm 4cm, width=0.55\textwidth]{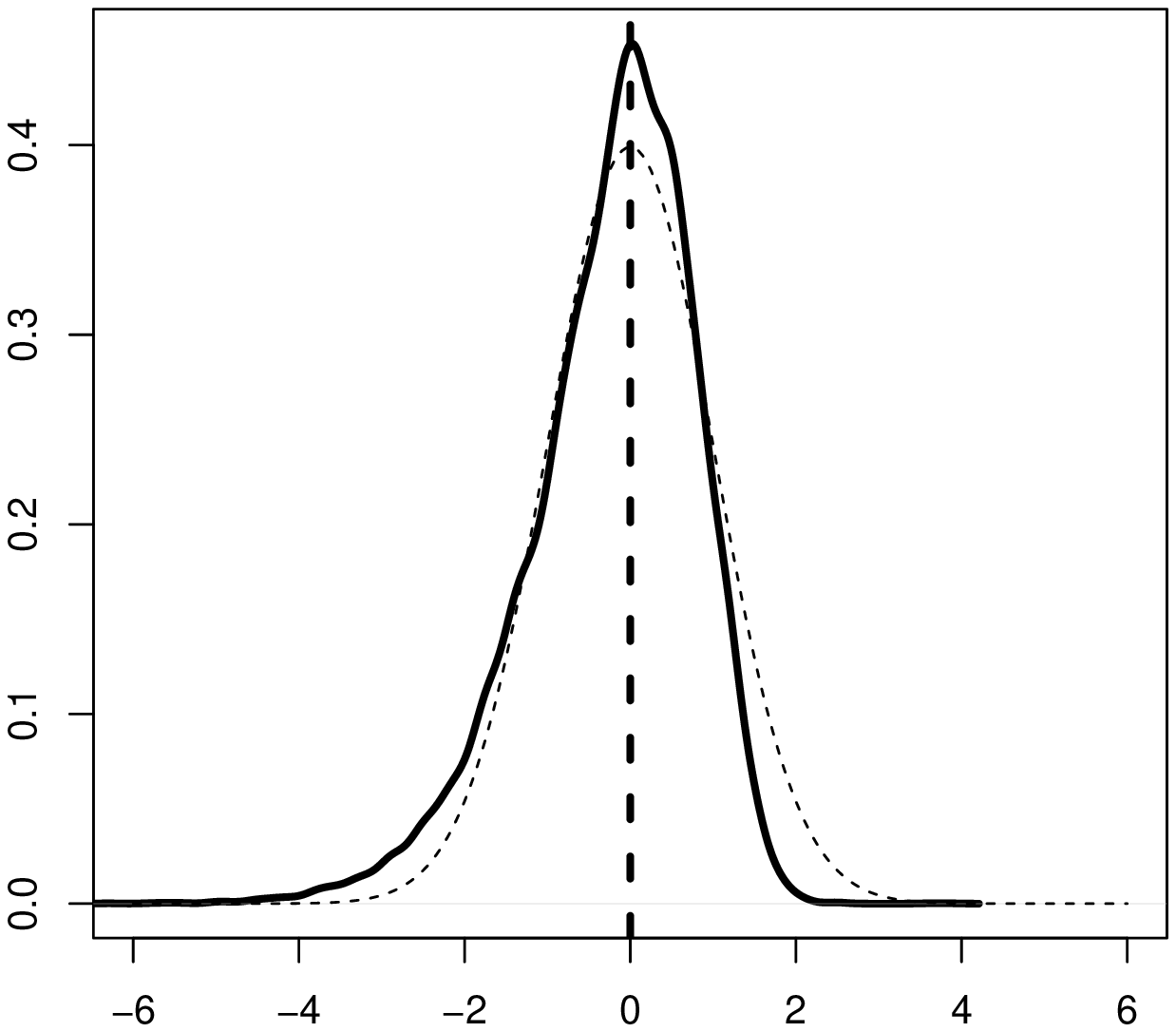}
\end{subfigure}\hspace{0.7cm}
\begin{subfigure}[b]{0.5\textwidth}
\vspace{0.7cm}\centering
\includegraphics[trim=5cm 5cm 2cm 4cm, width=0.55\textwidth]{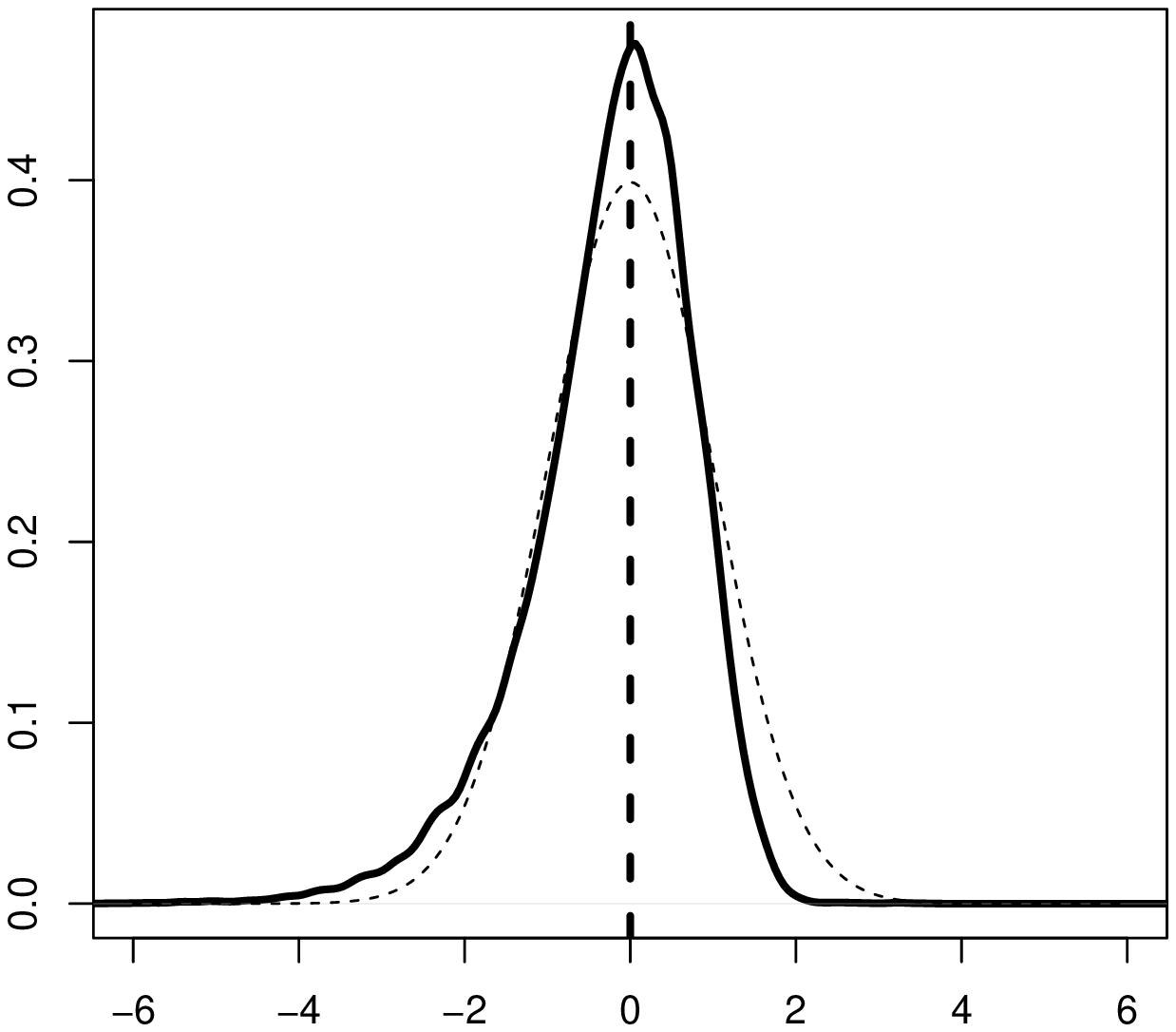}
\end{subfigure}\vspace{1.5cm}
\caption{{\small Distribution of $\sqrt{T}(\hat{\beta}_1(\hat{\tau}_T)-\beta_1)$ when $\tau_0=0.3$, $\beta_1=0.75$, $\beta_2=1-1/T$. Left: $\{\varepsilon_t\}_{t=1}^T\sim t(3)$; Right: $\{\varepsilon_t\}_{t=1}^T\sim t(2)$.}}
\end{figure}

\begin{figure}[H]
\label{fig1}
\begin{subfigure}[b]{0.5\textwidth}
\vspace{0.7cm}\centering
\includegraphics[trim=5cm 5cm 2cm 4cm, width=0.55\textwidth]{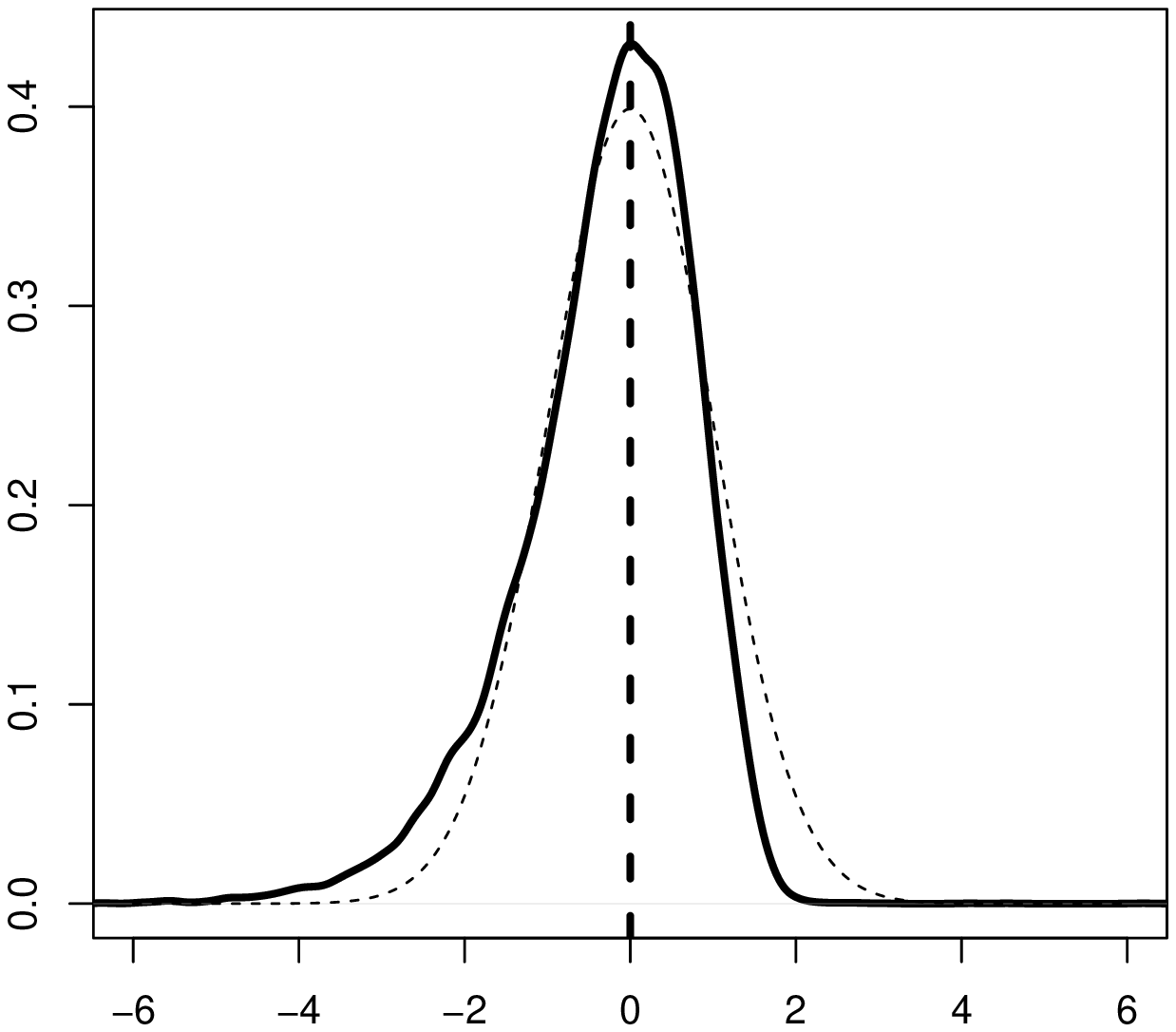}
\end{subfigure}\hspace{0.7cm}
\begin{subfigure}[b]{0.5\textwidth}
\vspace{0.7cm}\centering
\includegraphics[trim=5cm 5cm 2cm 4cm, width=0.55\textwidth]{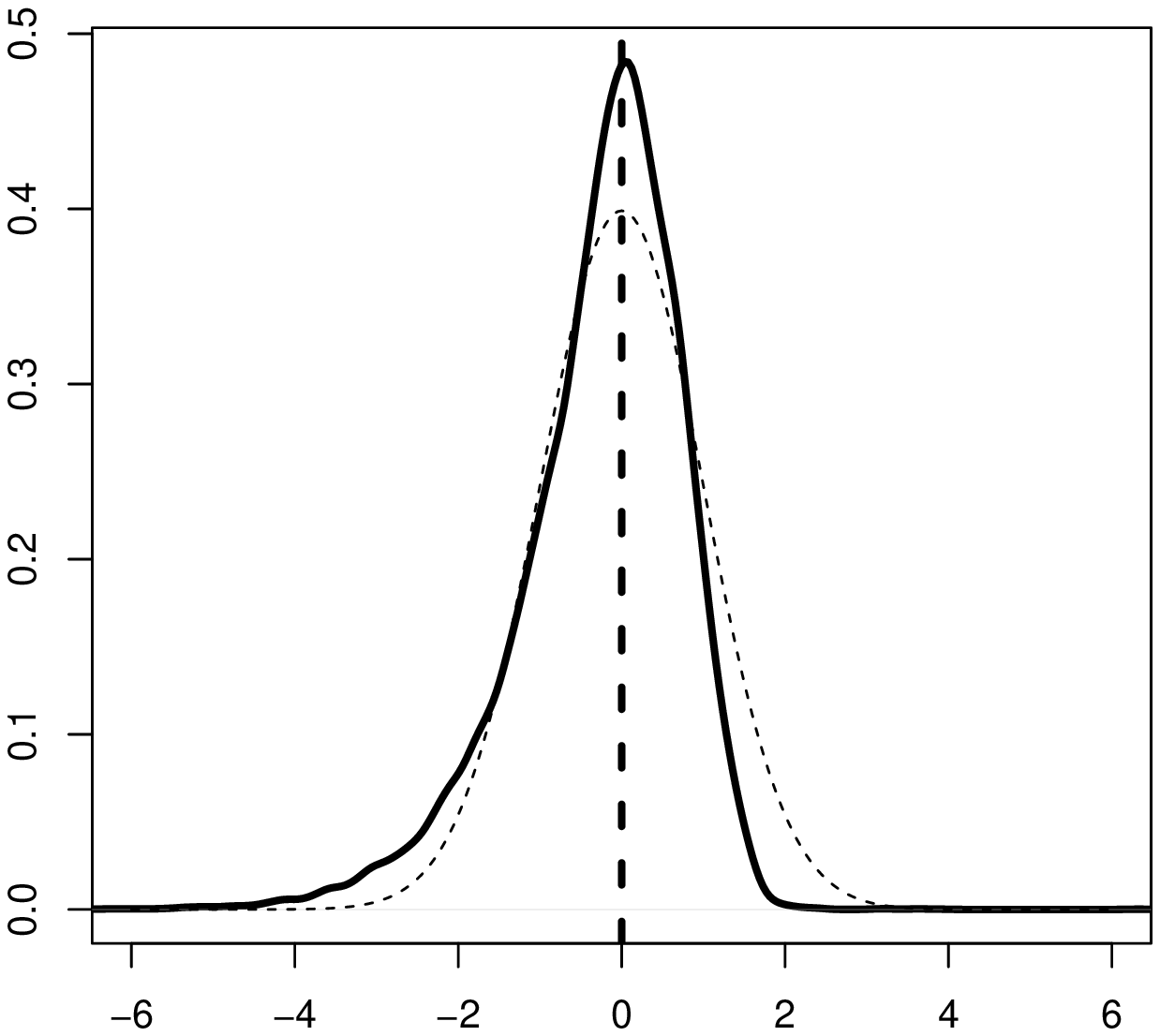}
\end{subfigure}\vspace{1.5cm}
\caption{{\small Distribution of $\sqrt{T}(\hat{\beta}_1(\hat{\tau}_T)-\beta_1)$ when $\tau_0=0.3$, $\beta_1=0.8$, $\beta_2=1-1/T$. Left: $\{\varepsilon_t\}_{t=1}^T\sim t(3)$; Right: $\{\varepsilon_t\}_{t=1}^T\sim t(2)$.}}
\end{figure}

\begin{figure}[H]
\label{fig1}
\begin{subfigure}[b]{0.5\textwidth}
\vspace{0.7cm}\centering
\includegraphics[trim=5cm 5cm 2cm 4cm, width=0.55\textwidth]{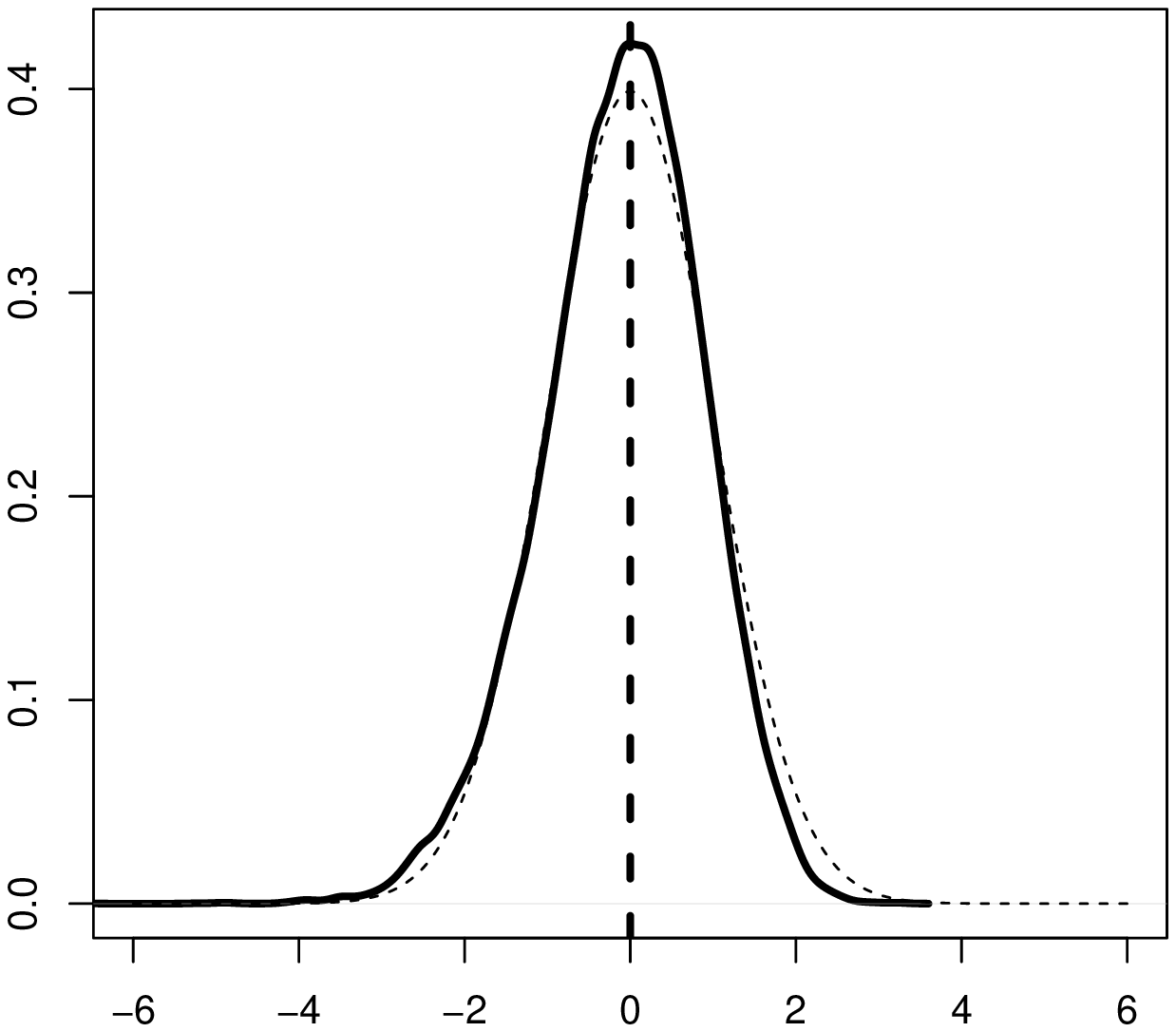}
\end{subfigure}\hspace{0.7cm}
\begin{subfigure}[b]{0.5\textwidth}
\vspace{0.7cm}\centering
\includegraphics[trim=5cm 5cm 2cm 4cm, width=0.55\textwidth]{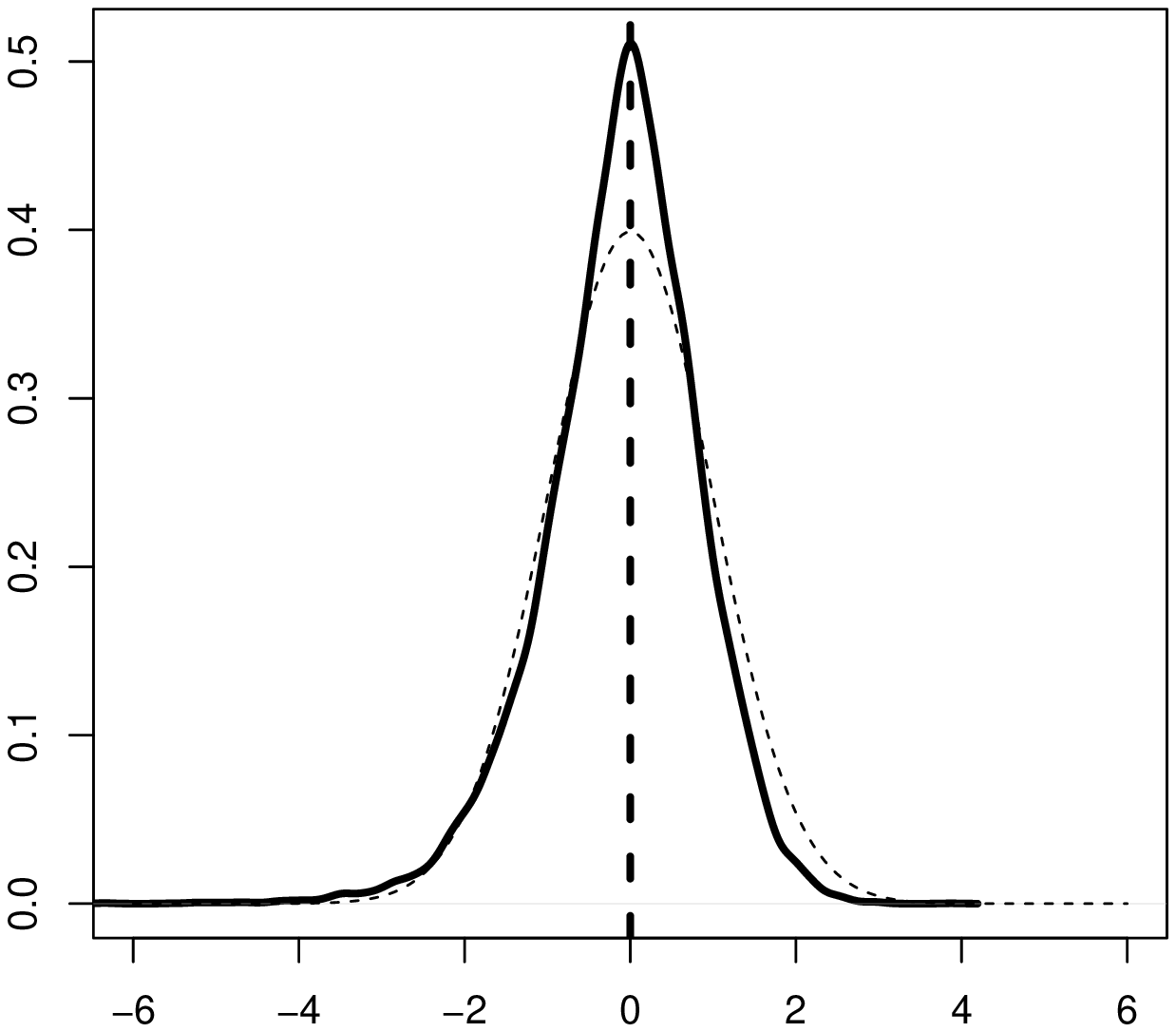}
\end{subfigure}\vspace{1.5cm}
\caption{{\small Distribution of $\sqrt{T}(\hat{\beta}_1(\hat{\tau}_T)-\beta_1)$ when $\tau_0=0.5$, $\beta_1=0.5$, $\beta_2=1-1/T$. Left: $\{\varepsilon_t\}_{t=1}^T\sim t(3)$; Right: $\{\varepsilon_t\}_{t=1}^T\sim t(2)$.}}
\end{figure}

\begin{figure}[H]
\label{fig1}
\begin{subfigure}[b]{0.5\textwidth}
\vspace{0.7cm}\centering
\includegraphics[trim=5cm 5cm 2cm 4cm, width=0.55\textwidth]{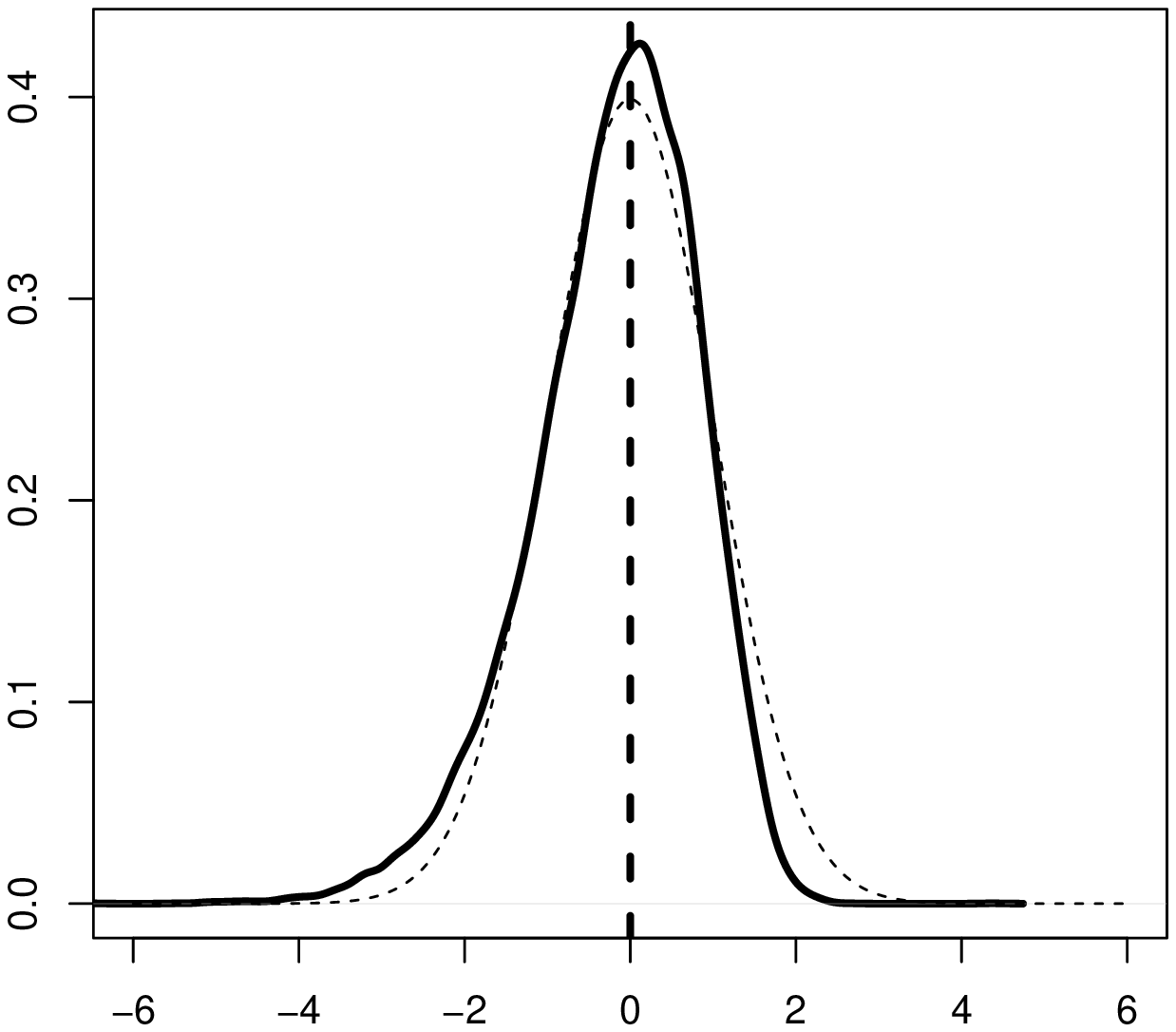}
\end{subfigure}\hspace{0.7cm}
\begin{subfigure}[b]{0.5\textwidth}
\vspace{0.7cm}\centering
\includegraphics[trim=5cm 5cm 2cm 4cm, width=0.55\textwidth]{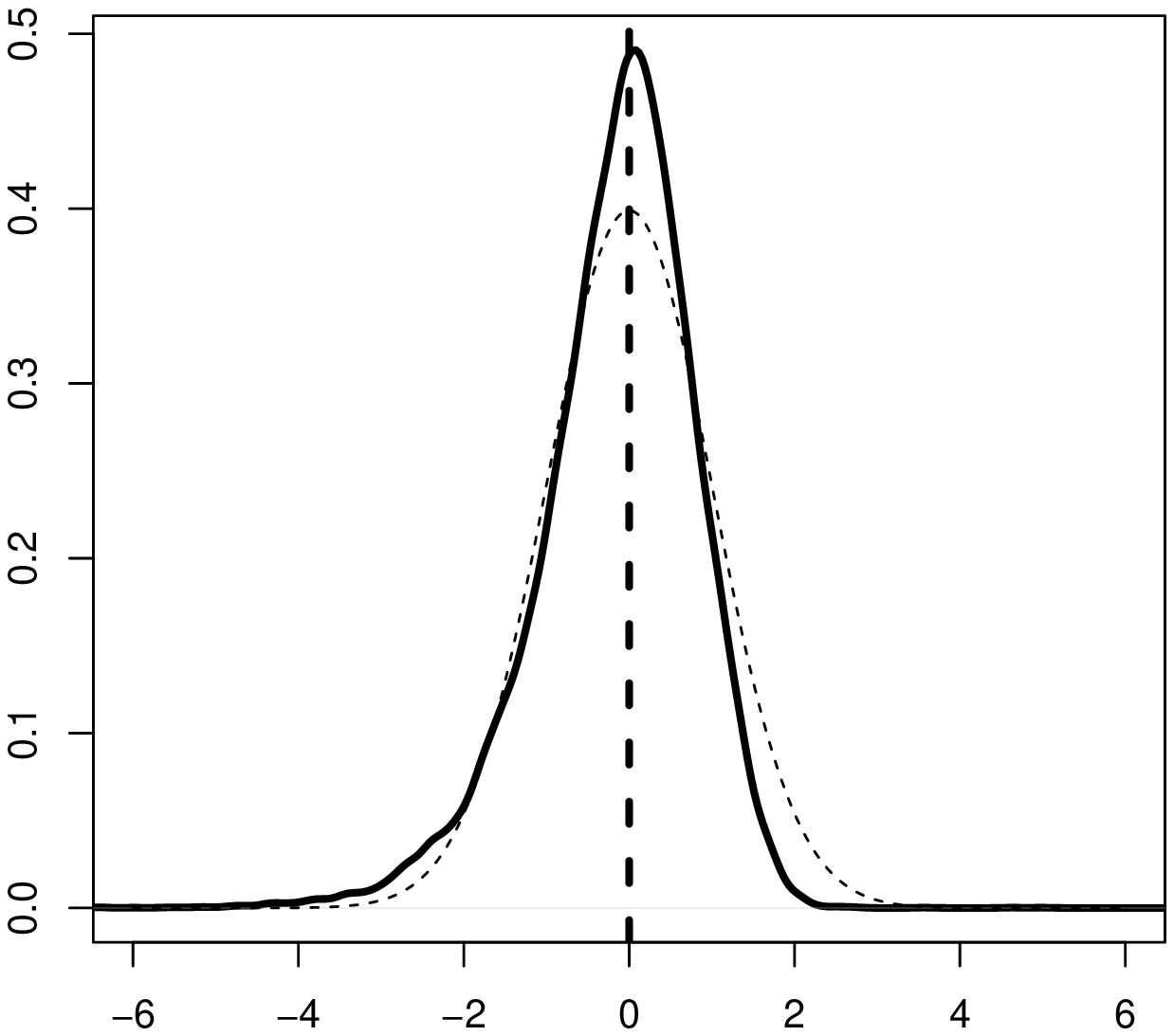}
\end{subfigure}\vspace{1.5cm}
\caption{{\small Distribution of $\sqrt{T}(\hat{\beta}_1(\hat{\tau}_T)-\beta_1)$ when $\tau_0=0.5$, $\beta_1=0.75$, $\beta_2=1-1/T$. Left: $\{\varepsilon_t\}_{t=1}^T\sim t(3)$; Right: $\{\varepsilon_t\}_{t=1}^T\sim t(2)$.}}
\end{figure}

\begin{figure}[H]
\label{fig1}
\begin{subfigure}[b]{0.5\textwidth}
\vspace{0.7cm}\centering
\includegraphics[trim=5cm 5cm 2cm 4cm, width=0.55\textwidth]{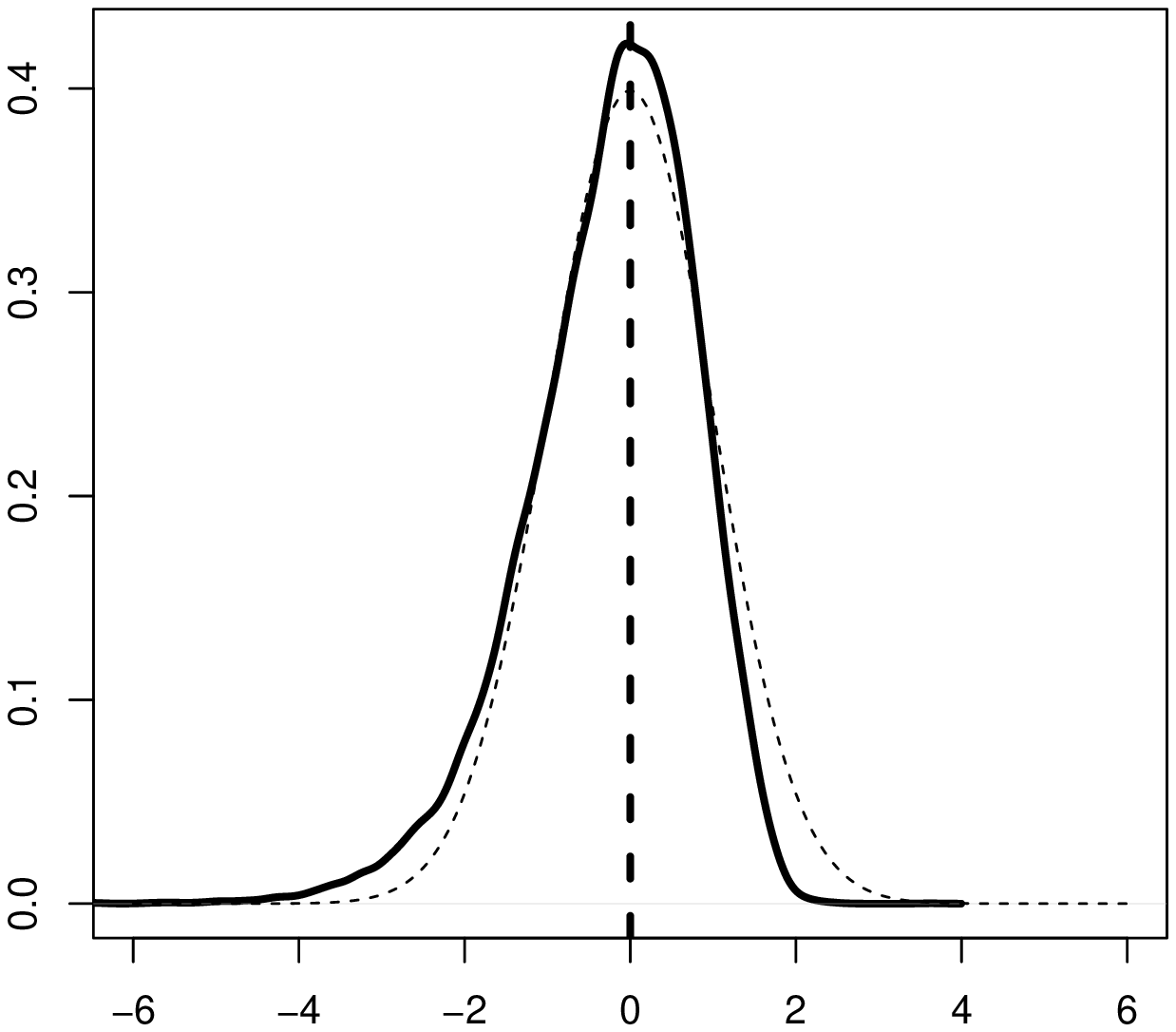}
\end{subfigure}\hspace{0.7cm}
\begin{subfigure}[b]{0.5\textwidth}
\vspace{0.7cm}\centering
\includegraphics[trim=5cm 5cm 2cm 4cm, width=0.55\textwidth]{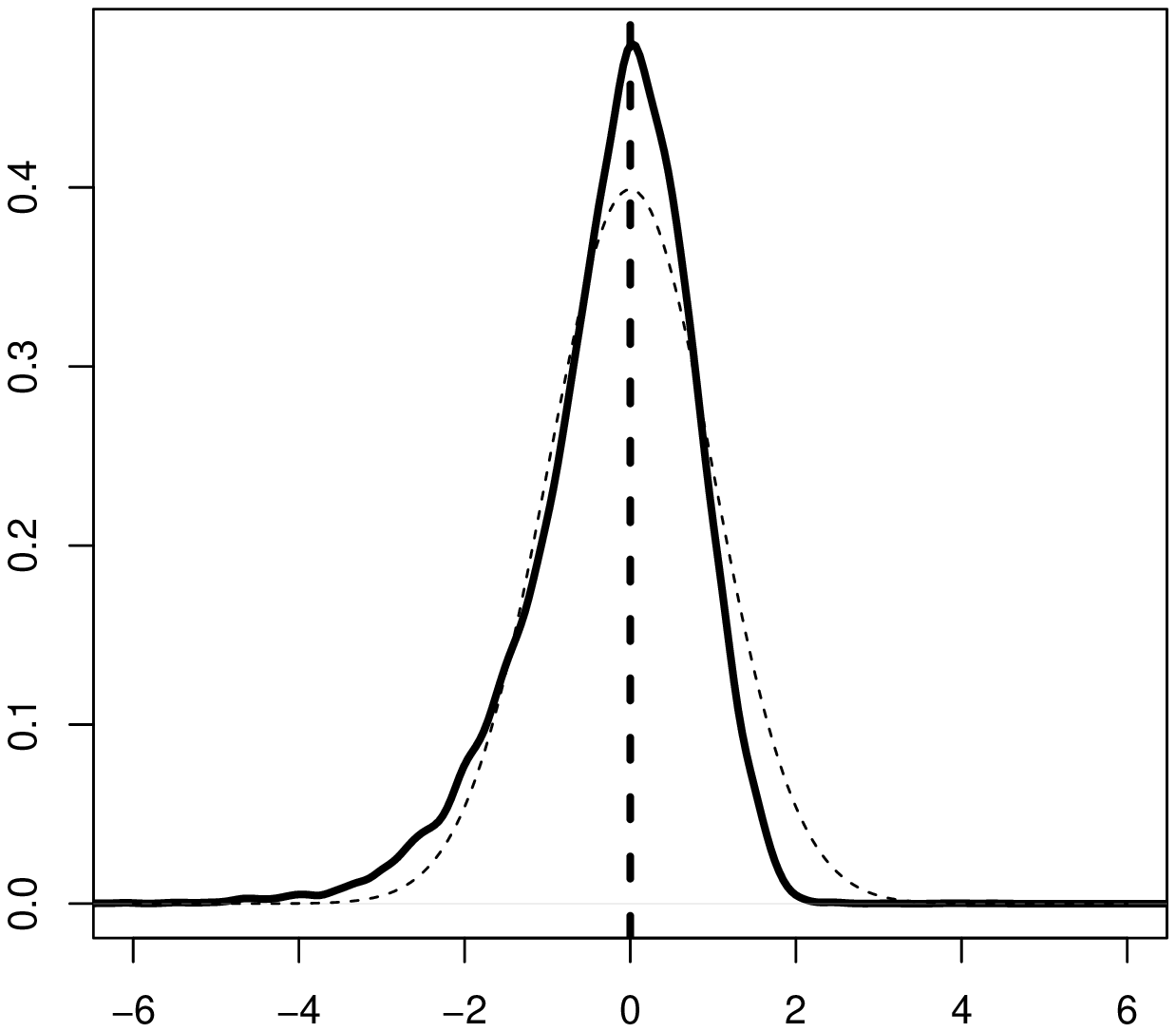}
\end{subfigure}\vspace{1.5cm}
\caption{{\small Distribution of $\sqrt{T}(\hat{\beta}_1(\hat{\tau}_T)-\beta_1)$ when $\tau_0=0.5$, $\beta_1=0.8$, $\beta_2=1-1/T$. Left: $\{\varepsilon_t\}_{t=1}^T\sim t(3)$; Right: $\{\varepsilon_t\}_{t=1}^T\sim t(2)$.}}
\end{figure}

\begin{figure}[H]
\label{fig7}
\begin{subfigure}[b]{0.5\textwidth}
\vspace{0.7cm}\centering
\includegraphics[trim=5cm 5cm 2cm 4cm, width=0.55\textwidth]{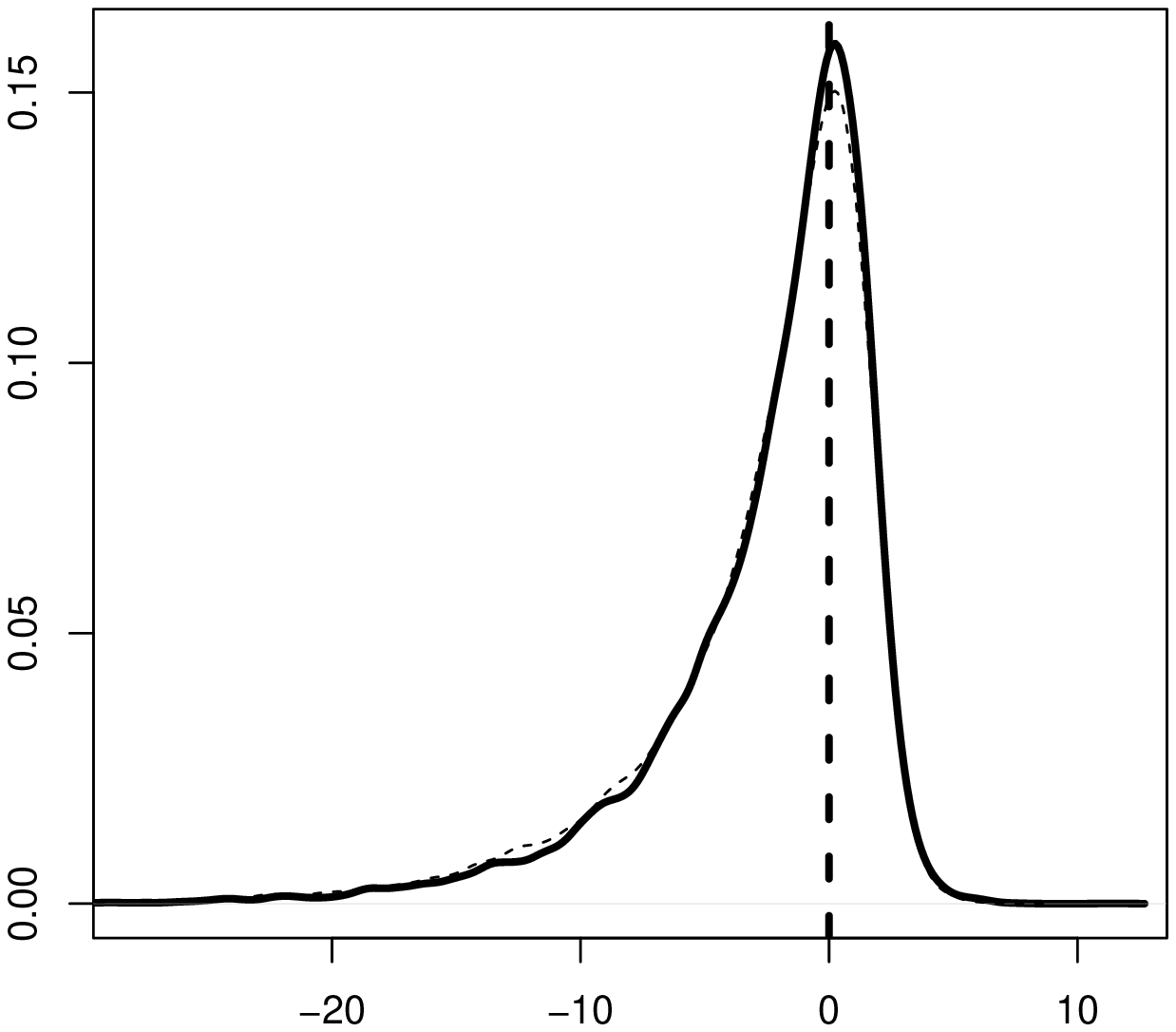}
\end{subfigure}\hspace{0.7cm}
\begin{subfigure}[b]{0.5\textwidth}
\vspace{0.7cm}\centering
\includegraphics[trim=5cm 5cm 2cm 4cm, width=0.55\textwidth]{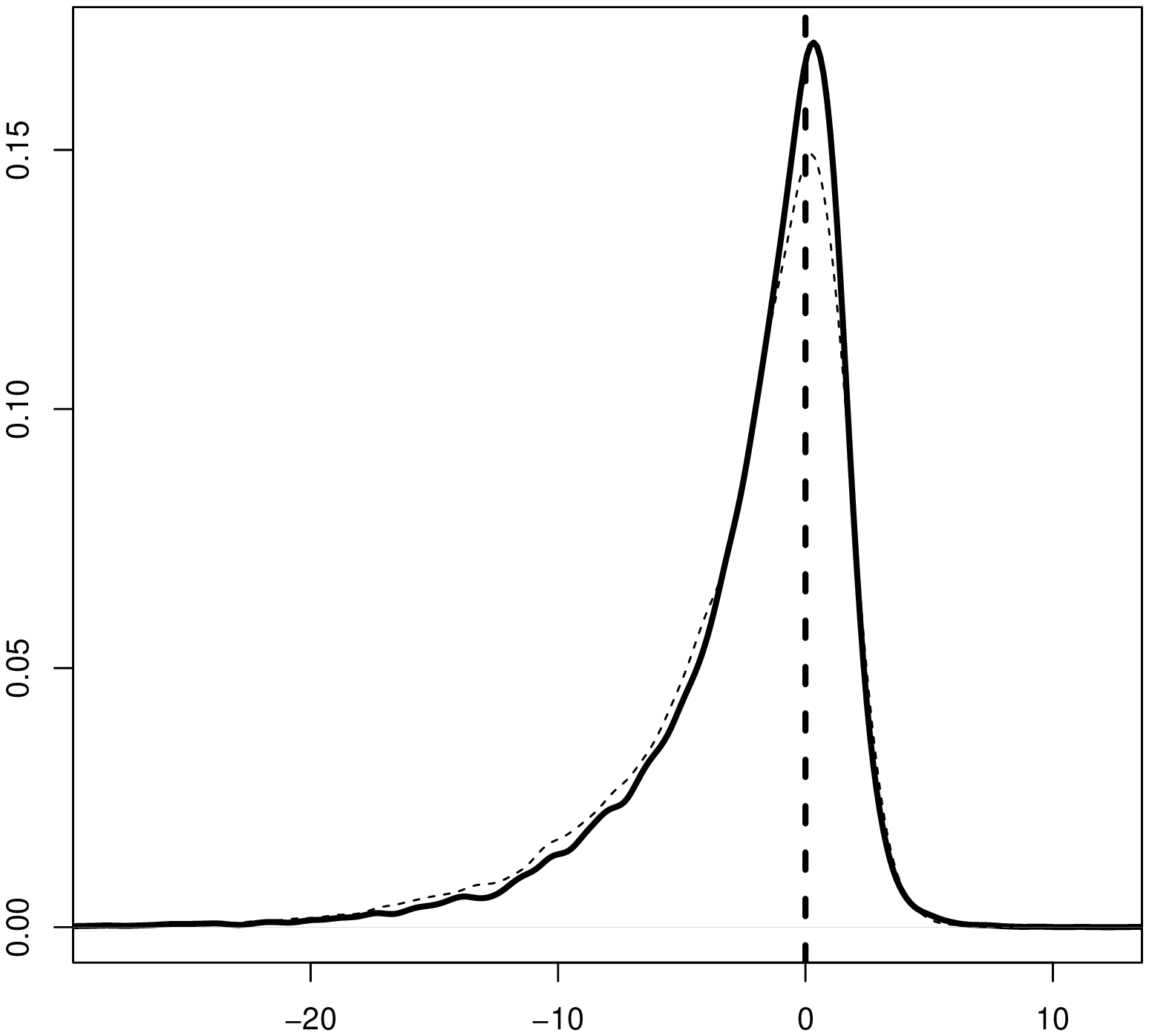}
\end{subfigure}\vspace{1.5cm}
\caption{{\small Distribution of $T(\hat{\beta}_2(\hat{\tau}_T)-\beta_2)$ when $\tau_0=0.3$, $\beta_1=0.5$, $\beta_2=1-1/T$. Left: $\{\varepsilon_t\}_{t=1}^T\sim t(3)$; Right: $\{\varepsilon_t\}_{t=1}^T\sim t(2)$.}}
\end{figure}

\begin{figure}[H]
\label{fig7}
\begin{subfigure}[b]{0.5\textwidth}
\vspace{0.7cm}\centering
\includegraphics[trim=5cm 5cm 2cm 4cm, width=0.55\textwidth]{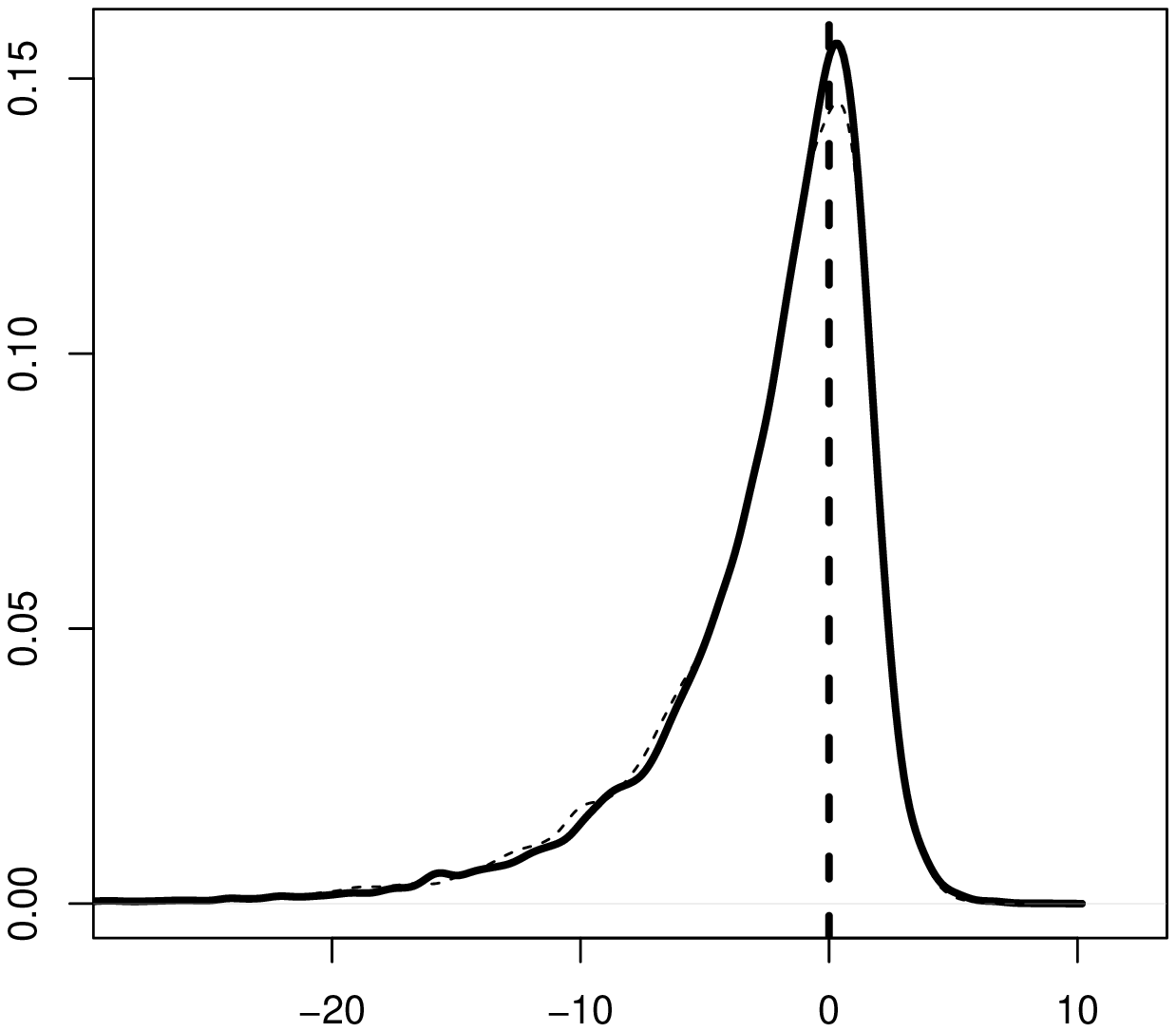}
\end{subfigure}\hspace{0.7cm}
\begin{subfigure}[b]{0.5\textwidth}
\vspace{0.7cm}\centering
\includegraphics[trim=5cm 5cm 2cm 4cm, width=0.55\textwidth]{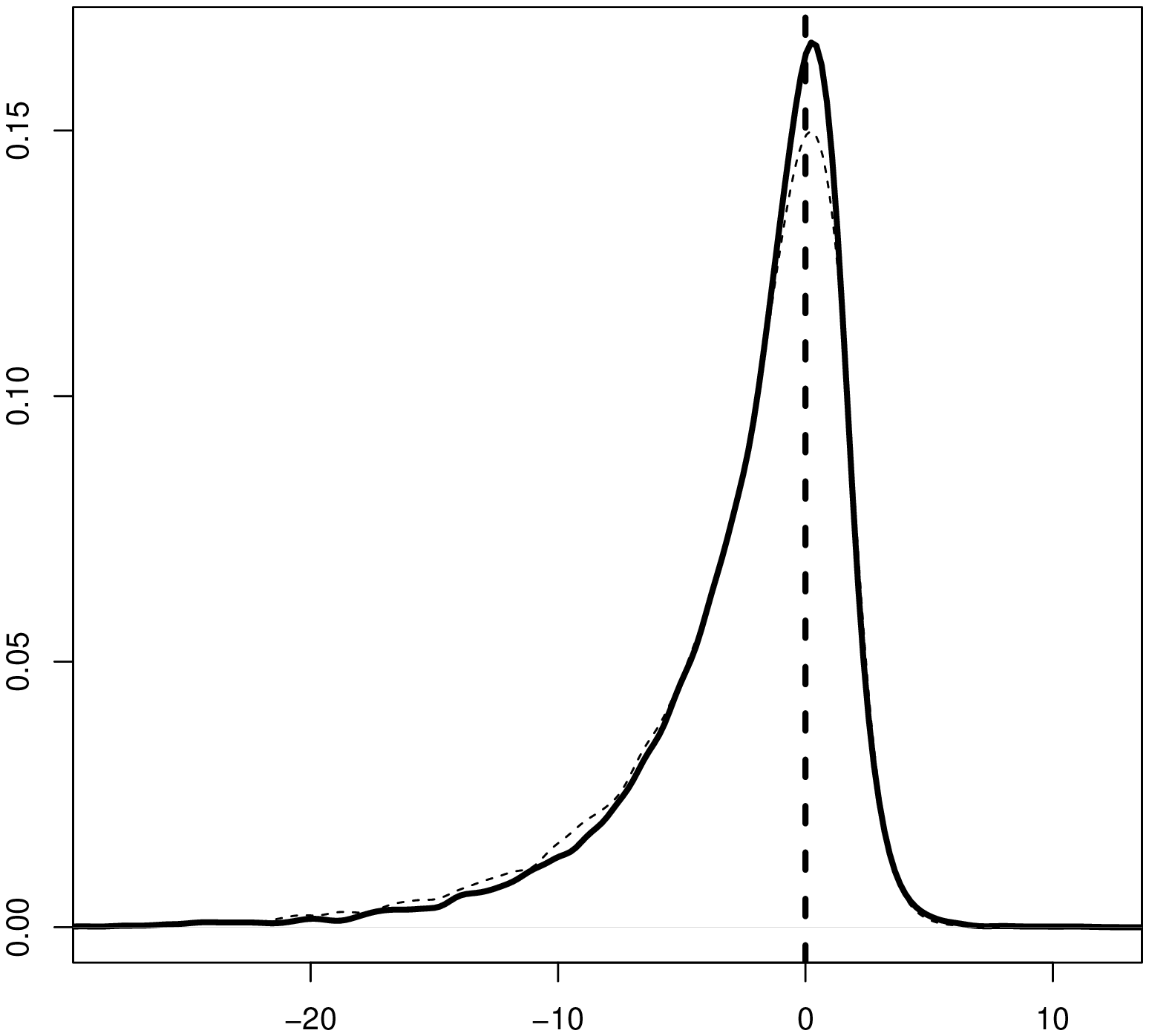}
\end{subfigure}\vspace{1.5cm}
\caption{{\small Distribution of $T(\hat{\beta}_2(\hat{\tau}_T)-\beta_2)$ when $\tau_0=0.3$, $\beta_1=0.75$, $\beta_2=1-1/T$. Left: $\{\varepsilon_t\}_{t=1}^T\sim t(3)$; Right: $\{\varepsilon_t\}_{t=1}^T\sim t(2)$.}}
\end{figure}

\begin{figure}[H]
\label{fig7}
\begin{subfigure}[b]{0.5\textwidth}
\vspace{0.7cm}\centering
\includegraphics[trim=5cm 5cm 2cm 4cm, width=0.55\textwidth]{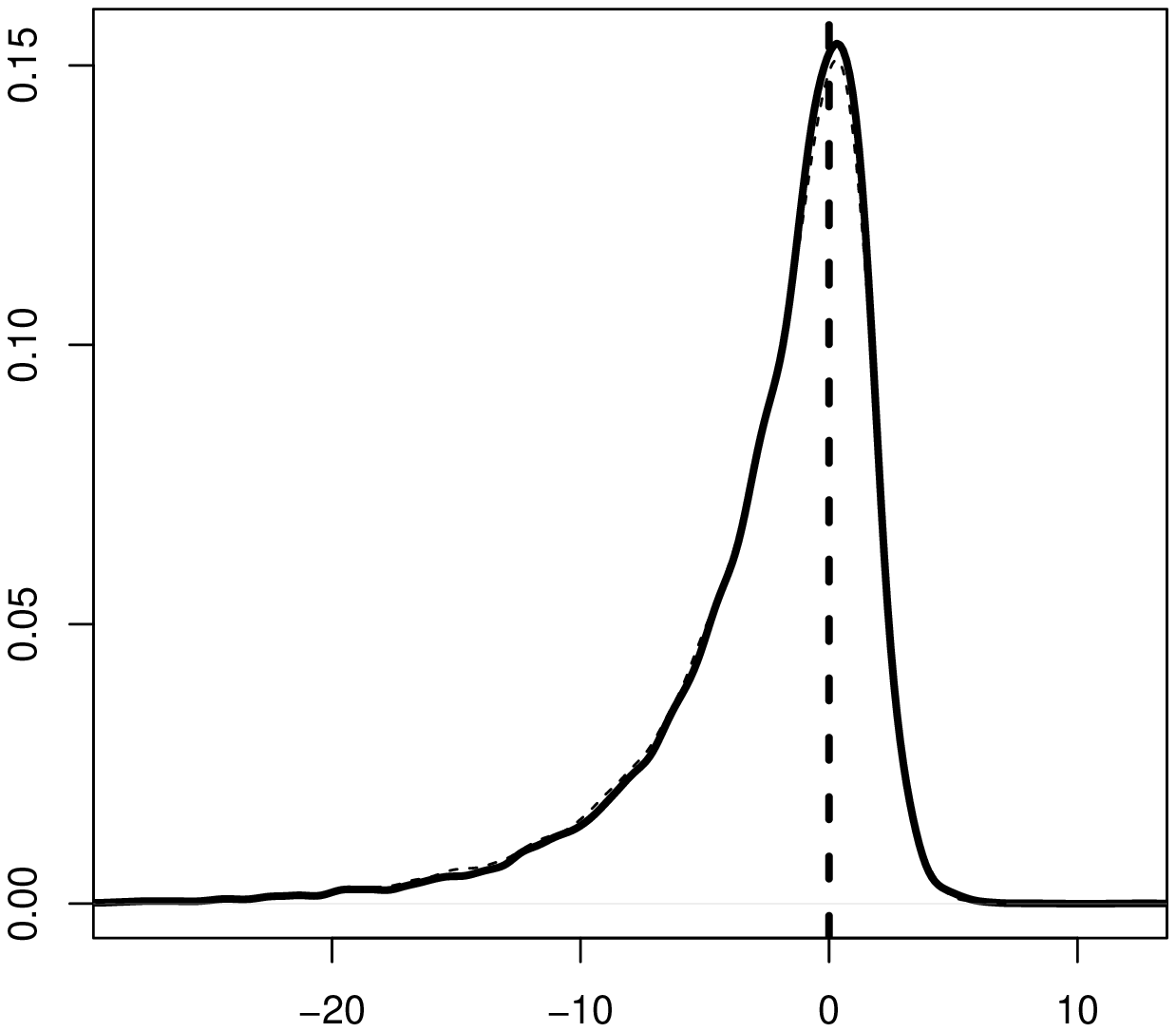}
\end{subfigure}\hspace{0.7cm}
\begin{subfigure}[b]{0.5\textwidth}
\vspace{0.7cm}\centering
\includegraphics[trim=5cm 5cm 2cm 4cm, width=0.55\textwidth]{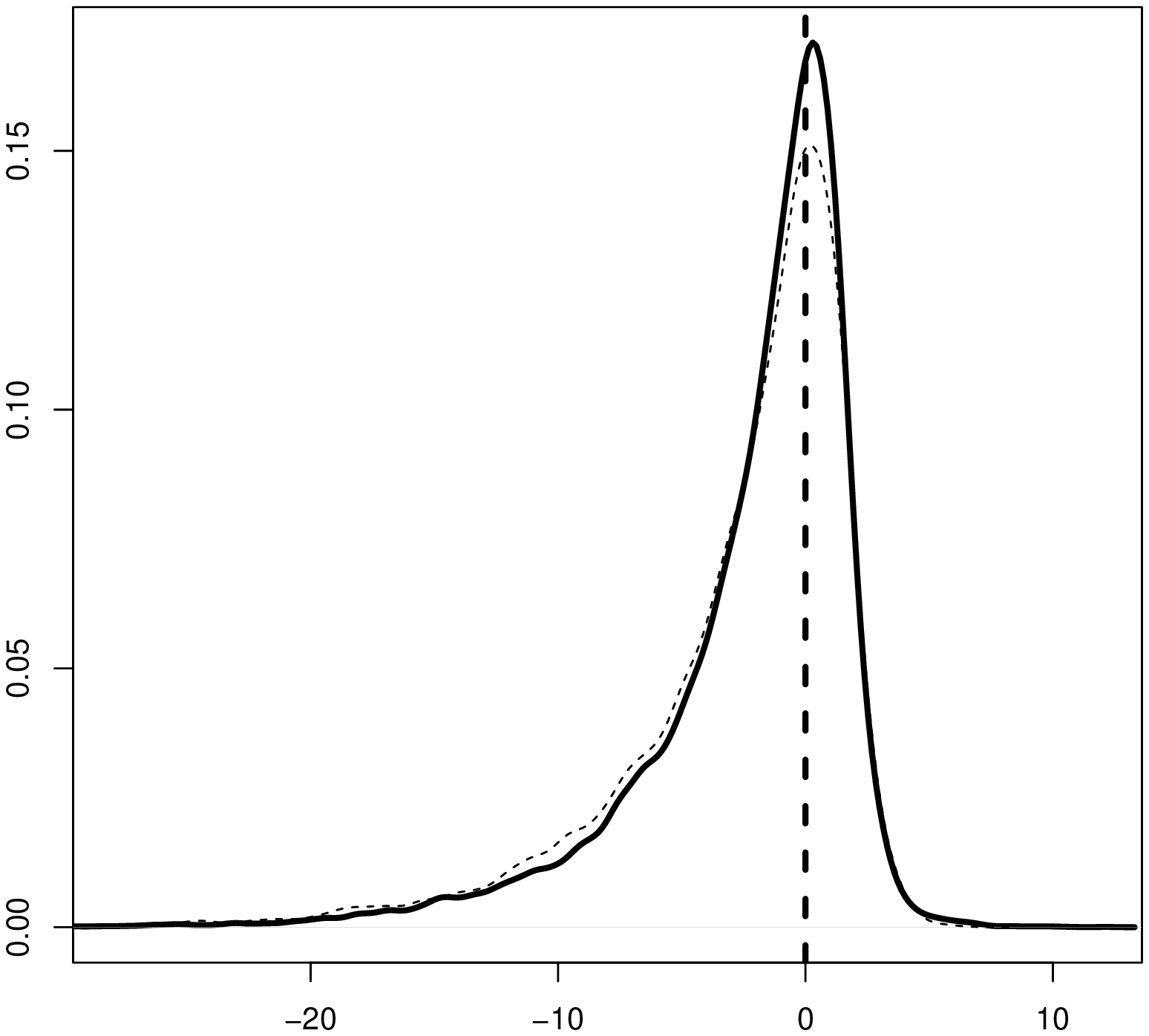}
\end{subfigure}\vspace{1.5cm}
\caption{{\small Distribution of $T(\hat{\beta}_2(\hat{\tau}_T)-\beta_2)$ when $\tau_0=0.3$, $\beta_1=0.8$, $\beta_2=1-1/T$. Left: $\{\varepsilon_t\}_{t=1}^T\sim t(3)$; Right: $\{\varepsilon_t\}_{t=1}^T\sim t(2)$.}}
\end{figure}

\begin{figure}[H]
\label{fig7}
\begin{subfigure}[b]{0.5\textwidth}
\vspace{0.7cm}\centering
\includegraphics[trim=5cm 5cm 2cm 4cm, width=0.55\textwidth]{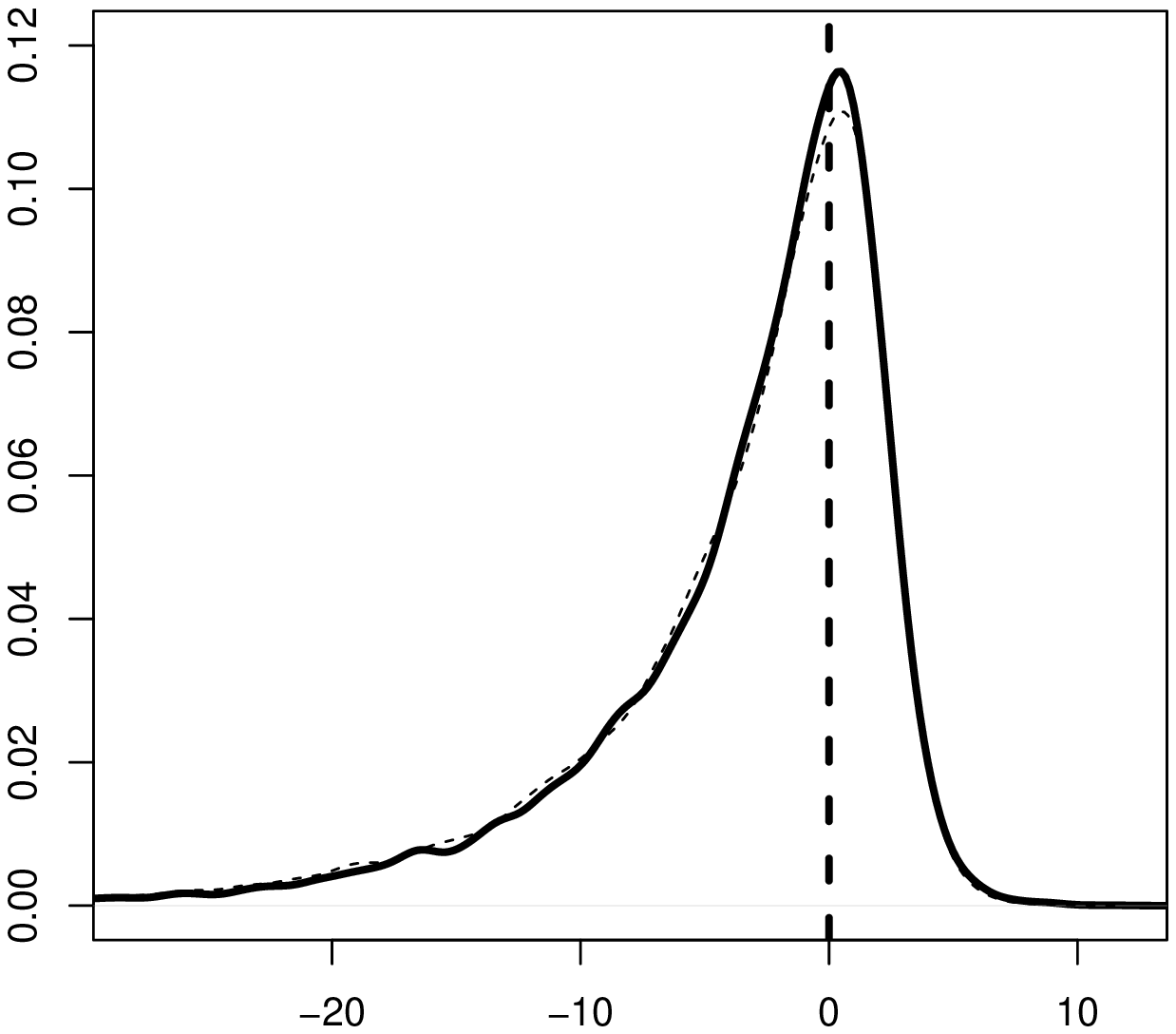}
\end{subfigure}\hspace{0.7cm}
\begin{subfigure}[b]{0.5\textwidth}
\vspace{0.7cm}\centering
\includegraphics[trim=5cm 5cm 2cm 4cm, width=0.55\textwidth]{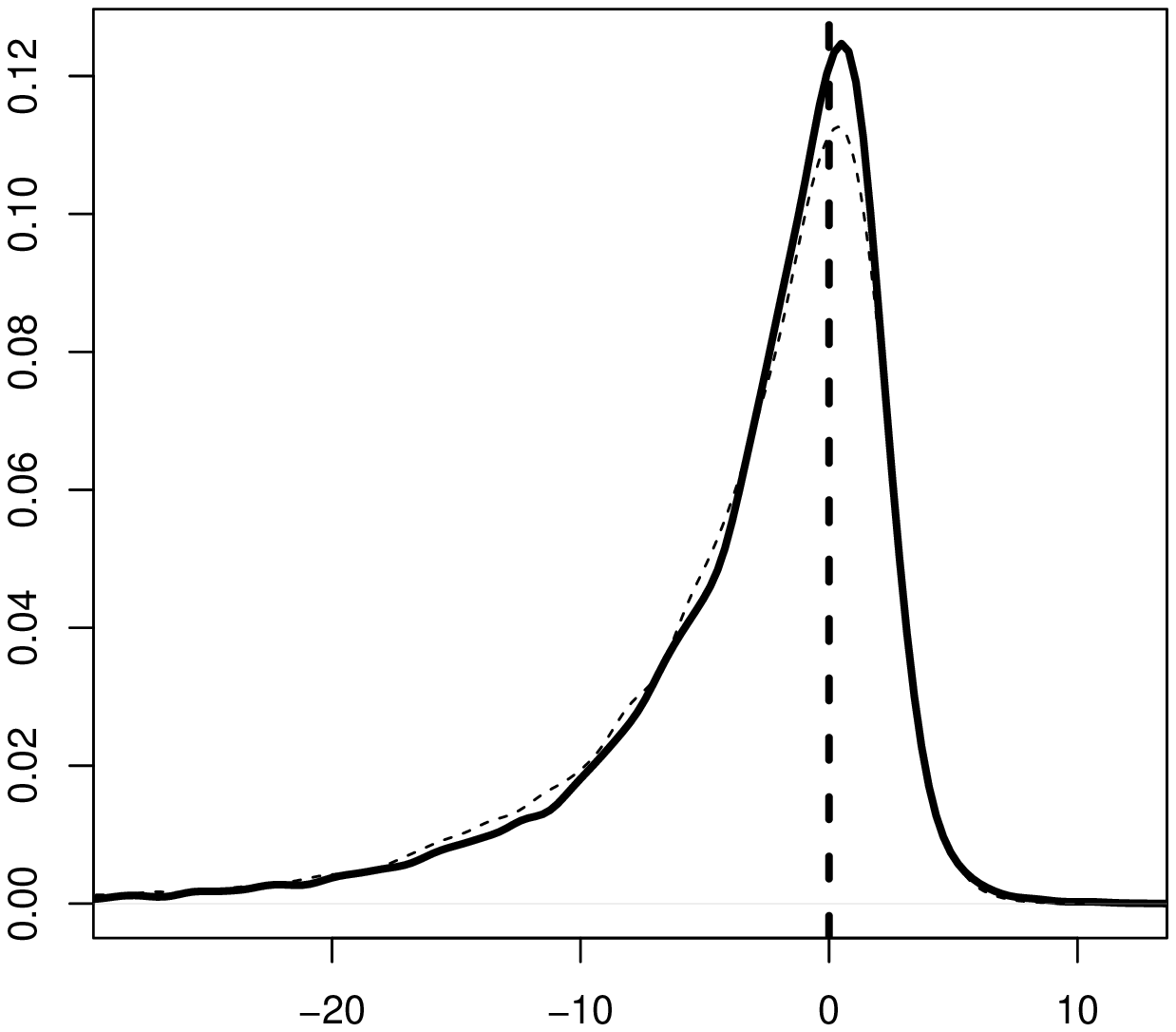}
\end{subfigure}\vspace{1.5cm}
\caption{{\small Distribution of $T(\hat{\beta}_2(\hat{\tau}_T)-\beta_2)$ when $\tau_0=0.5$, $\beta_1=0.5$, $\beta_2=1-1/T$. Left: $\{\varepsilon_t\}_{t=1}^T\sim t(3)$; Right: $\{\varepsilon_t\}_{t=1}^T\sim t(2)$.}}
\end{figure}

\begin{figure}[H]
\label{fig7}
\begin{subfigure}[b]{0.5\textwidth}
\vspace{0.7cm}\centering
\includegraphics[trim=5cm 5cm 2cm 4cm, width=0.55\textwidth]{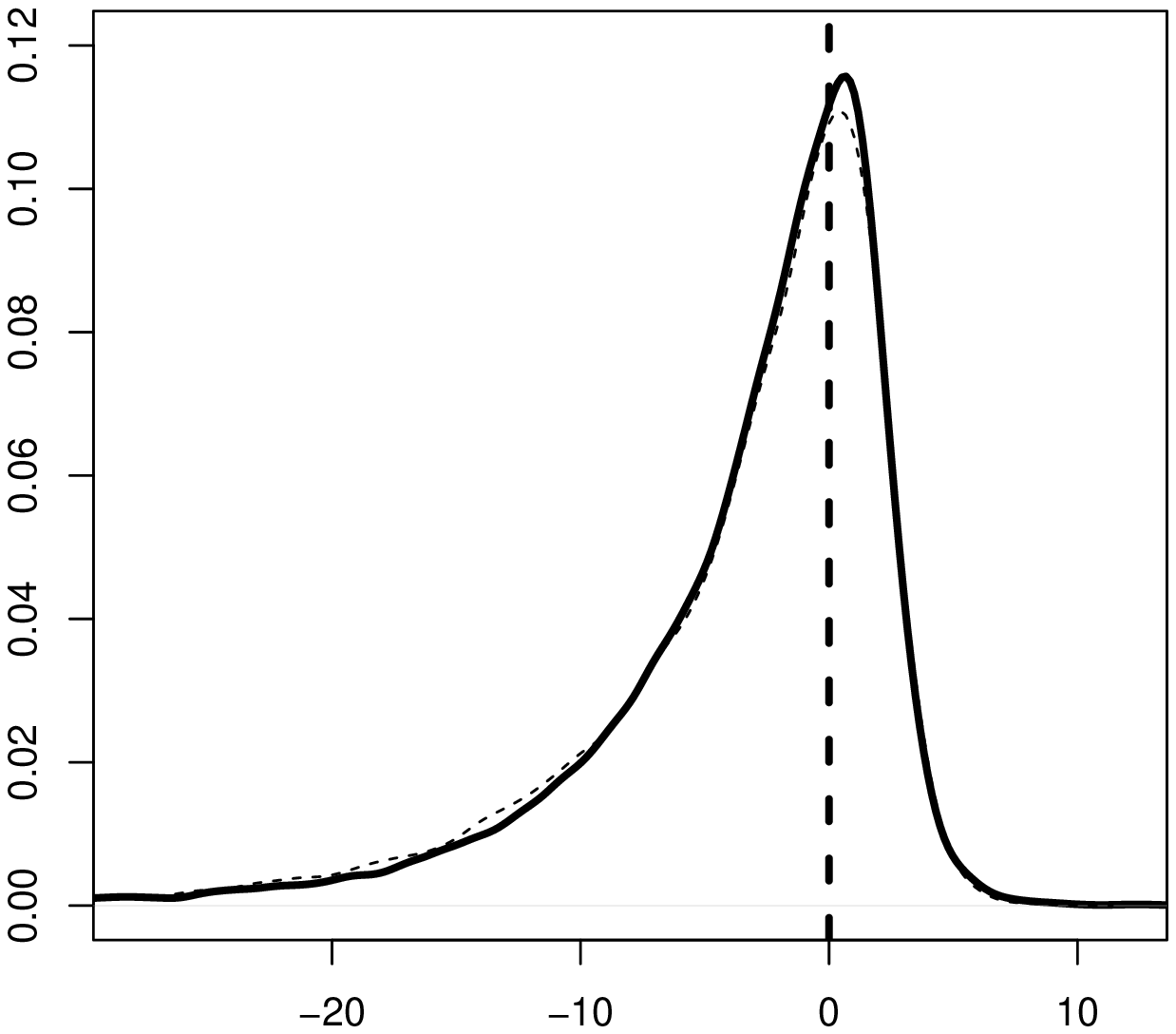}
\end{subfigure}\hspace{0.7cm}
\begin{subfigure}[b]{0.5\textwidth}
\vspace{0.7cm}\centering
\includegraphics[trim=5cm 5cm 2cm 4cm, width=0.55\textwidth]{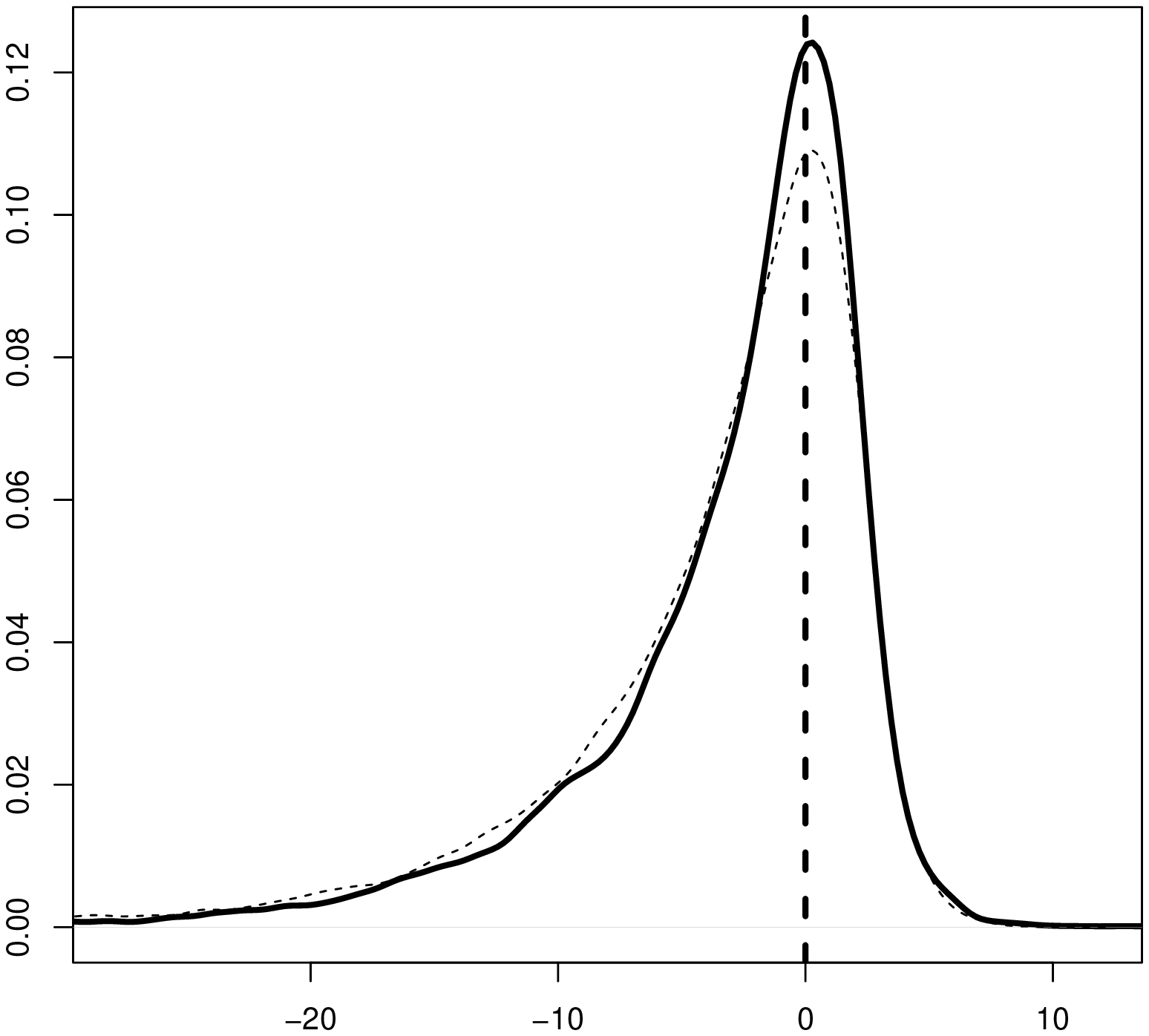}
\end{subfigure}\vspace{1.5cm}
\caption{{\small Distribution of $T(\hat{\beta}_2(\hat{\tau}_T)-\beta_2)$ when $\tau_0=0.5$, $\beta_1=0.75$, $\beta_2=1-1/T$. Left: $\{\varepsilon_t\}_{t=1}^T\sim t(3)$; Right: $\{\varepsilon_t\}_{t=1}^T\sim t(2)$.}}
\end{figure}

\begin{figure}[H]
\label{fig7}
\begin{subfigure}[b]{0.5\textwidth}
\vspace{0.7cm}\centering
\includegraphics[trim=5cm 5cm 2cm 4cm, width=0.55\textwidth]{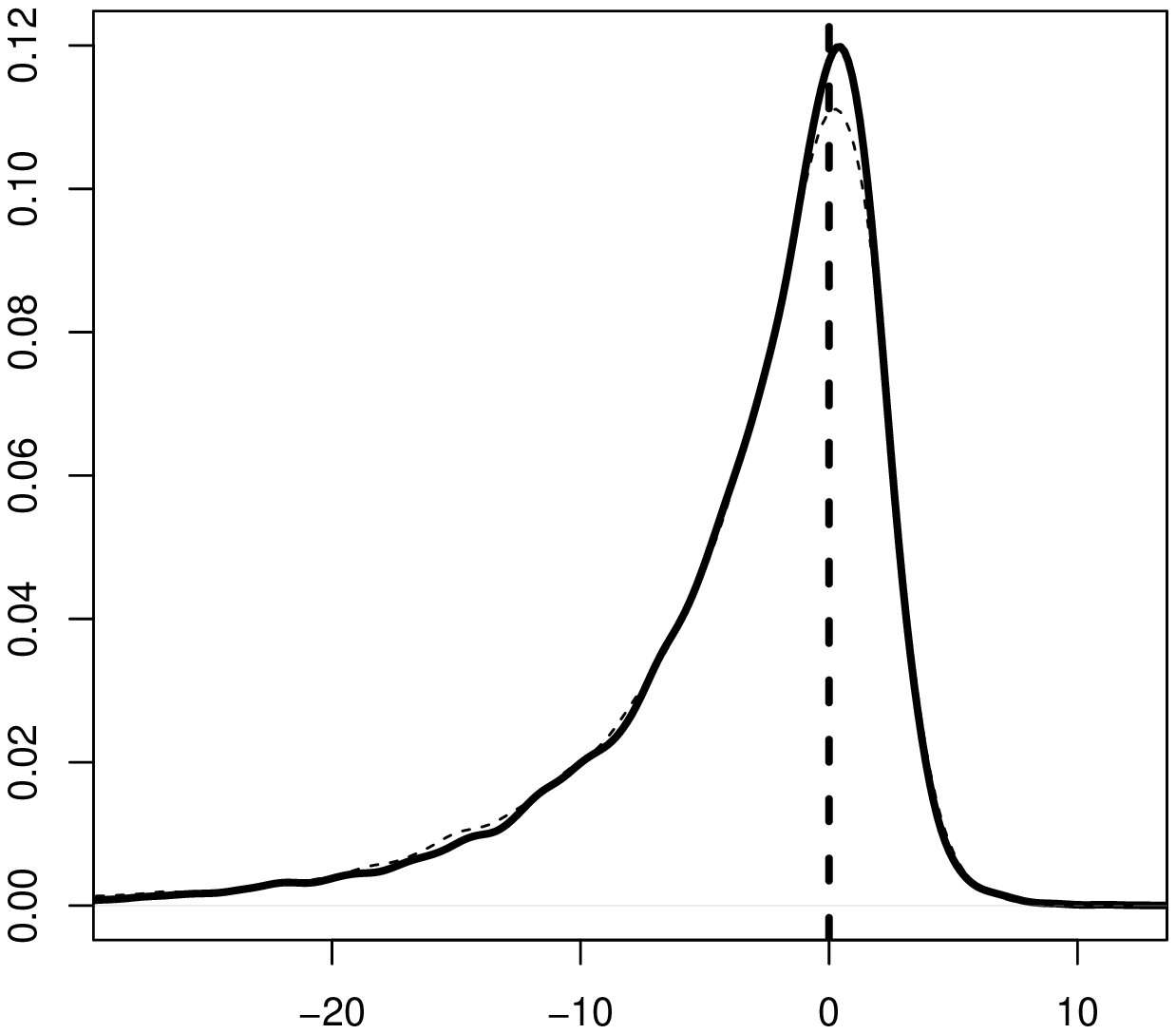}
\end{subfigure}\hspace{0.7cm}
\begin{subfigure}[b]{0.5\textwidth}
\vspace{0.7cm}\centering
\includegraphics[trim=5cm 5cm 2cm 4cm, width=0.55\textwidth]{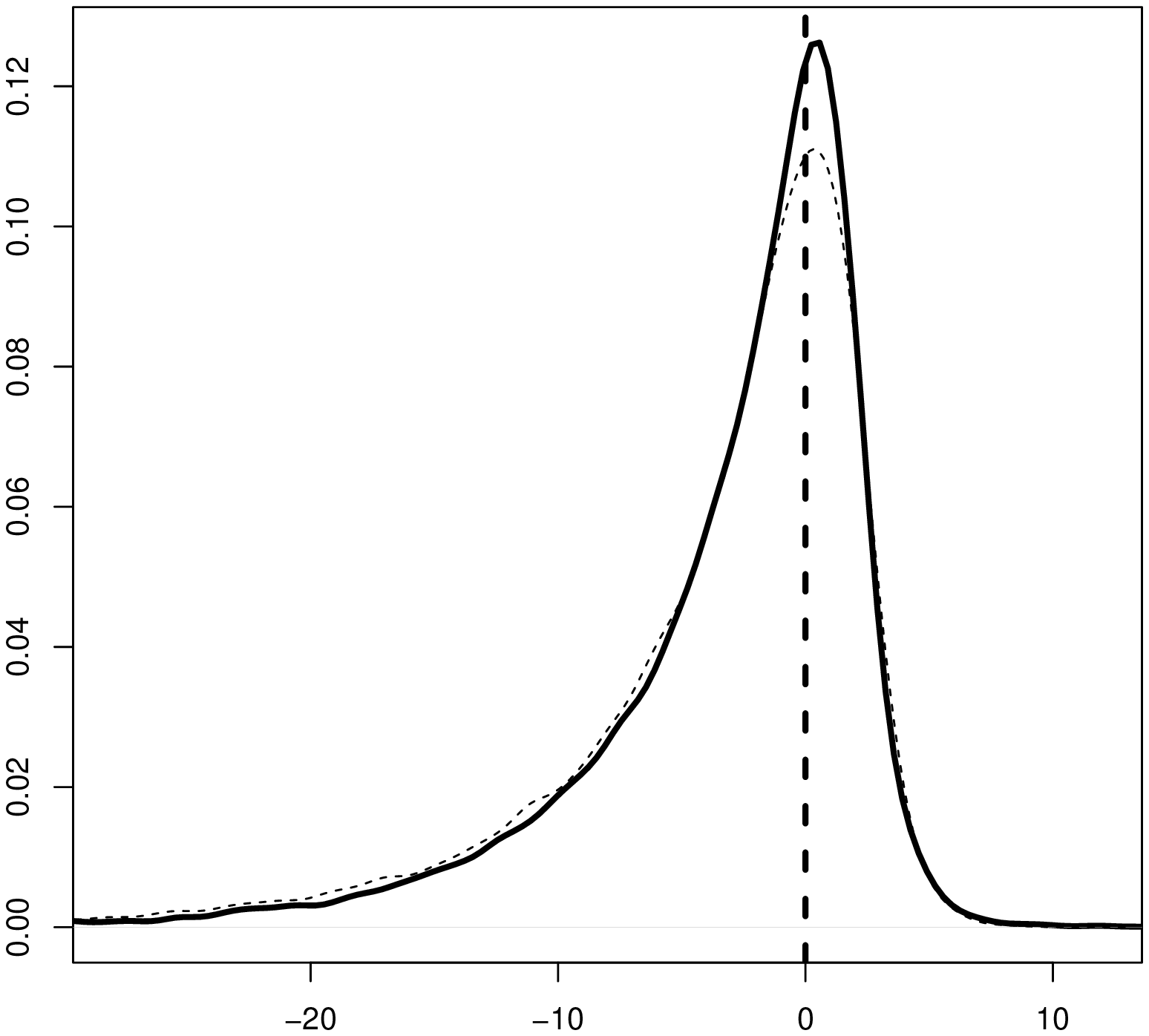}
\end{subfigure}\vspace{1.5cm}
\caption{{\small Distribution of $T(\hat{\beta}_2(\hat{\tau}_T)-\beta_2)$ when $\tau_0=0.5$, $\beta_1=0.8$, $\beta_2=1-1/T$. Left: $\{\varepsilon_t\}_{t=1}^T\sim t(3)$; Right: $\{\varepsilon_t\}_{t=1}^T\sim t(2)$.}}
\end{figure}

\begin{figure}[H]
\label{fig13}
\begin{subfigure}[b]{0.5\textwidth}
\vspace{0.7cm}\centering
\includegraphics[trim=5cm 5cm 2cm 4cm, width=0.55\textwidth]{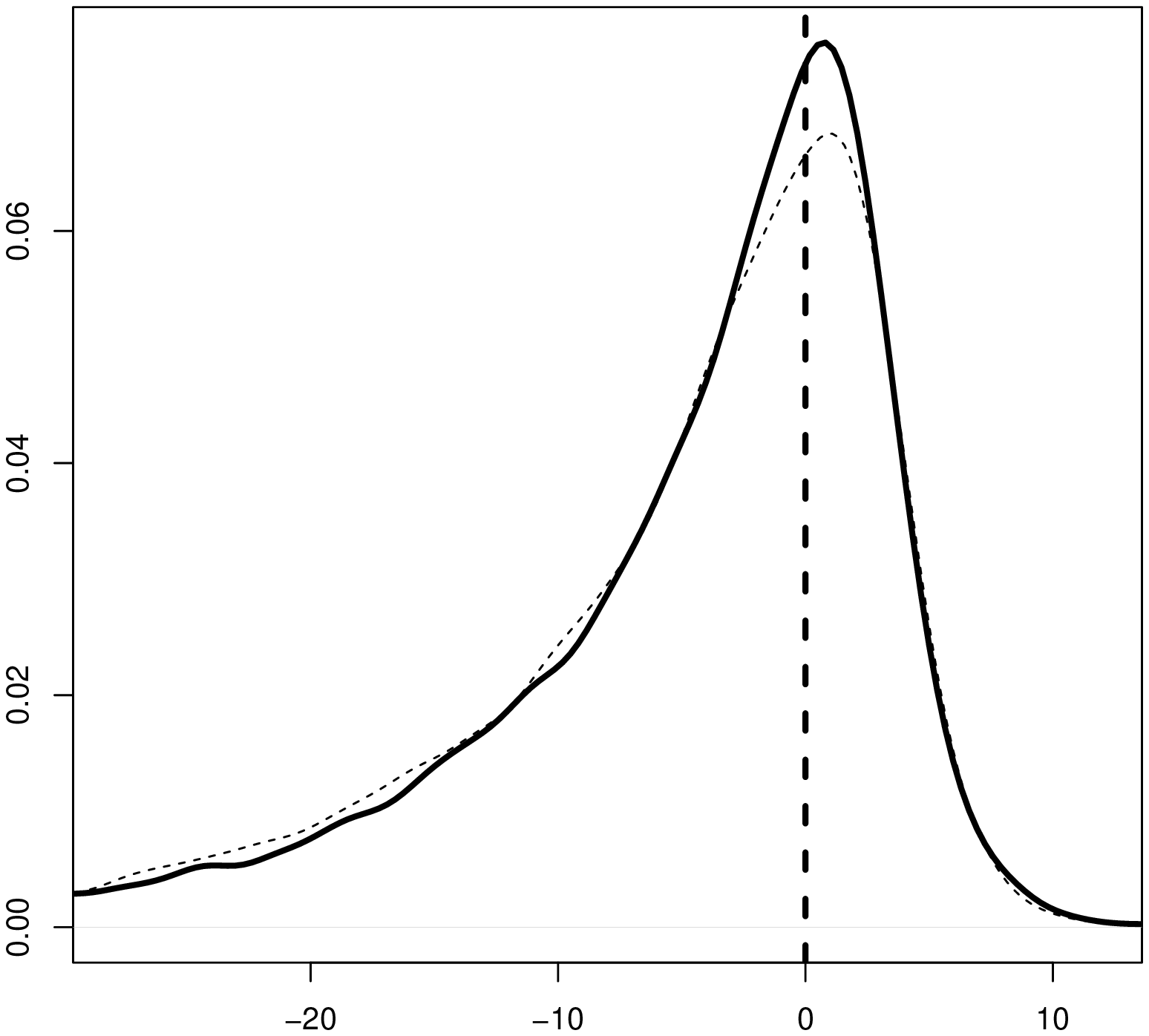}
\end{subfigure}\hspace{0.7cm}
\begin{subfigure}[b]{0.5\textwidth}
\vspace{0.7cm}\centering
\includegraphics[trim=5cm 5cm 2cm 4cm, width=0.55\textwidth]{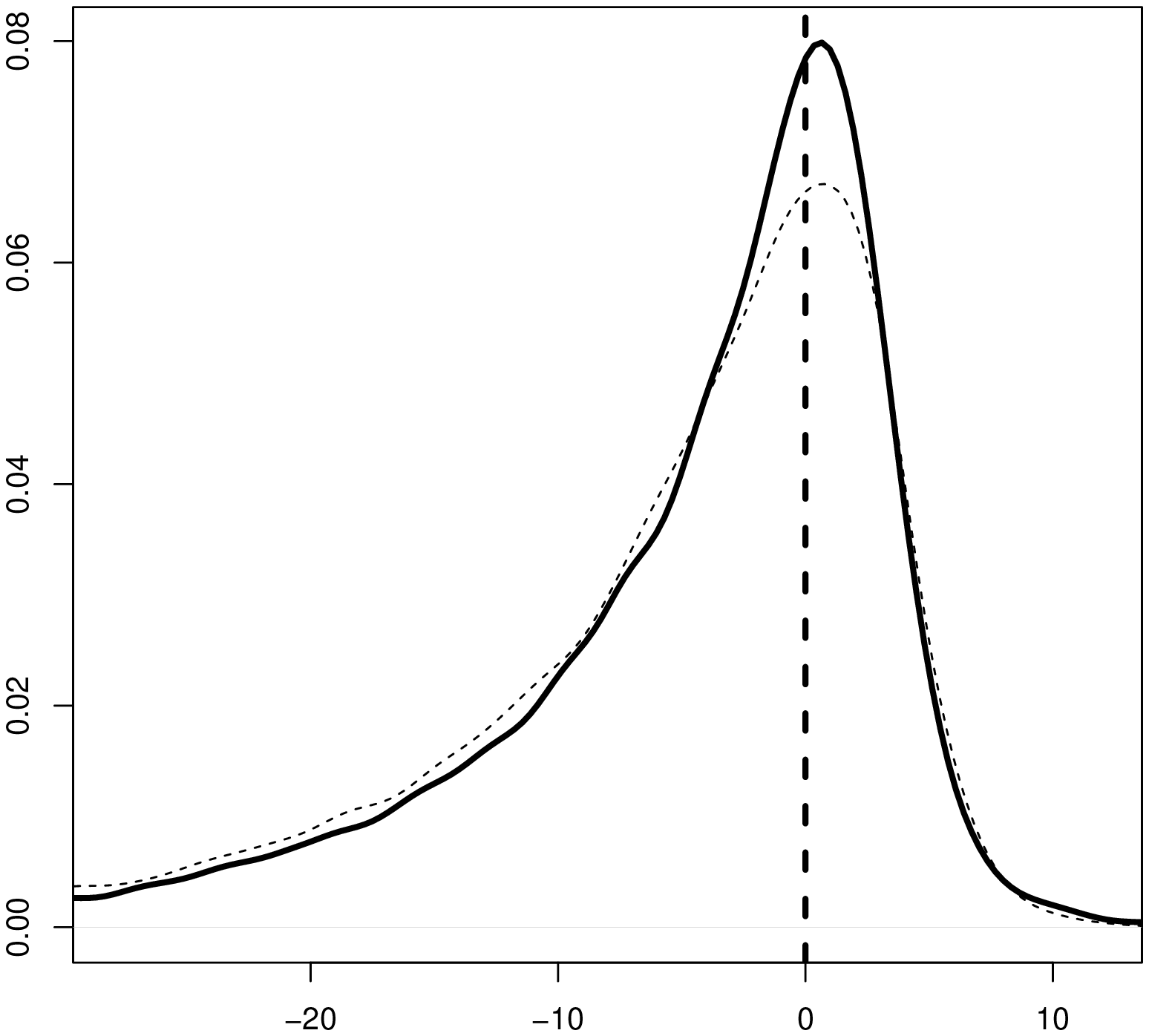}
\end{subfigure}\vspace{1.5cm}
\caption{{\small Distribution of $T(\hat{\beta}_1(\hat{\tau}_T)-\beta_1)$ when $\tau_0=0.3$, $\beta_2=0.5$, $\beta_1=1-1/T$. Left:
$\{\varepsilon_t\}_{t=1}^T\sim t(3)$; Right: $\{\varepsilon_t\}_{t=1}^T\sim t(2)$.}}
\end{figure}

\begin{figure}[H]
\label{fig13}
\begin{subfigure}[b]{0.5\textwidth}
\vspace{0.7cm}\centering
\includegraphics[trim=5cm 5cm 2cm 4cm, width=0.55\textwidth]{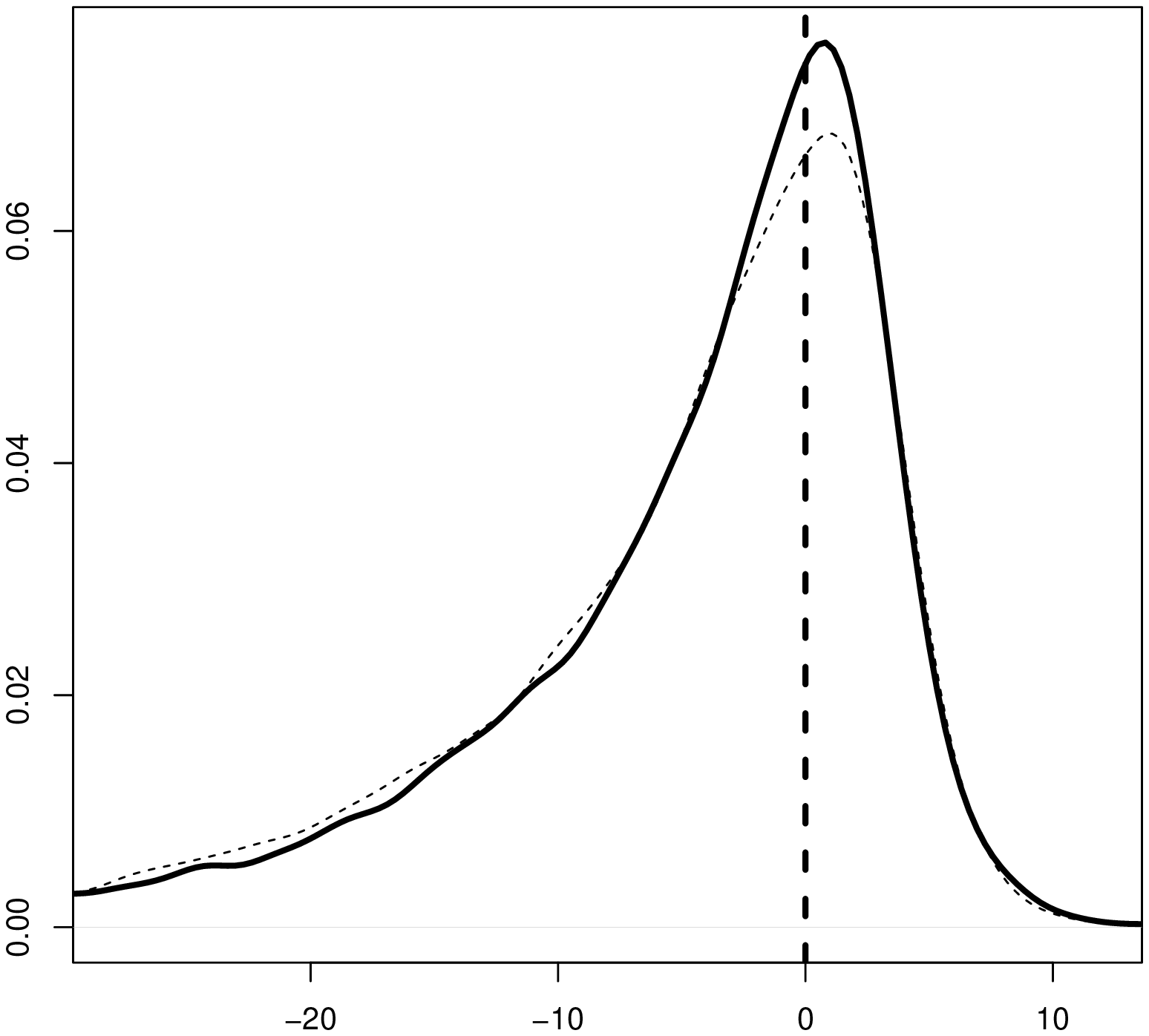}
\end{subfigure}\hspace{0.7cm}
\begin{subfigure}[b]{0.5\textwidth}
\vspace{0.7cm}\centering
\includegraphics[trim=5cm 5cm 2cm 4cm, width=0.55\textwidth]{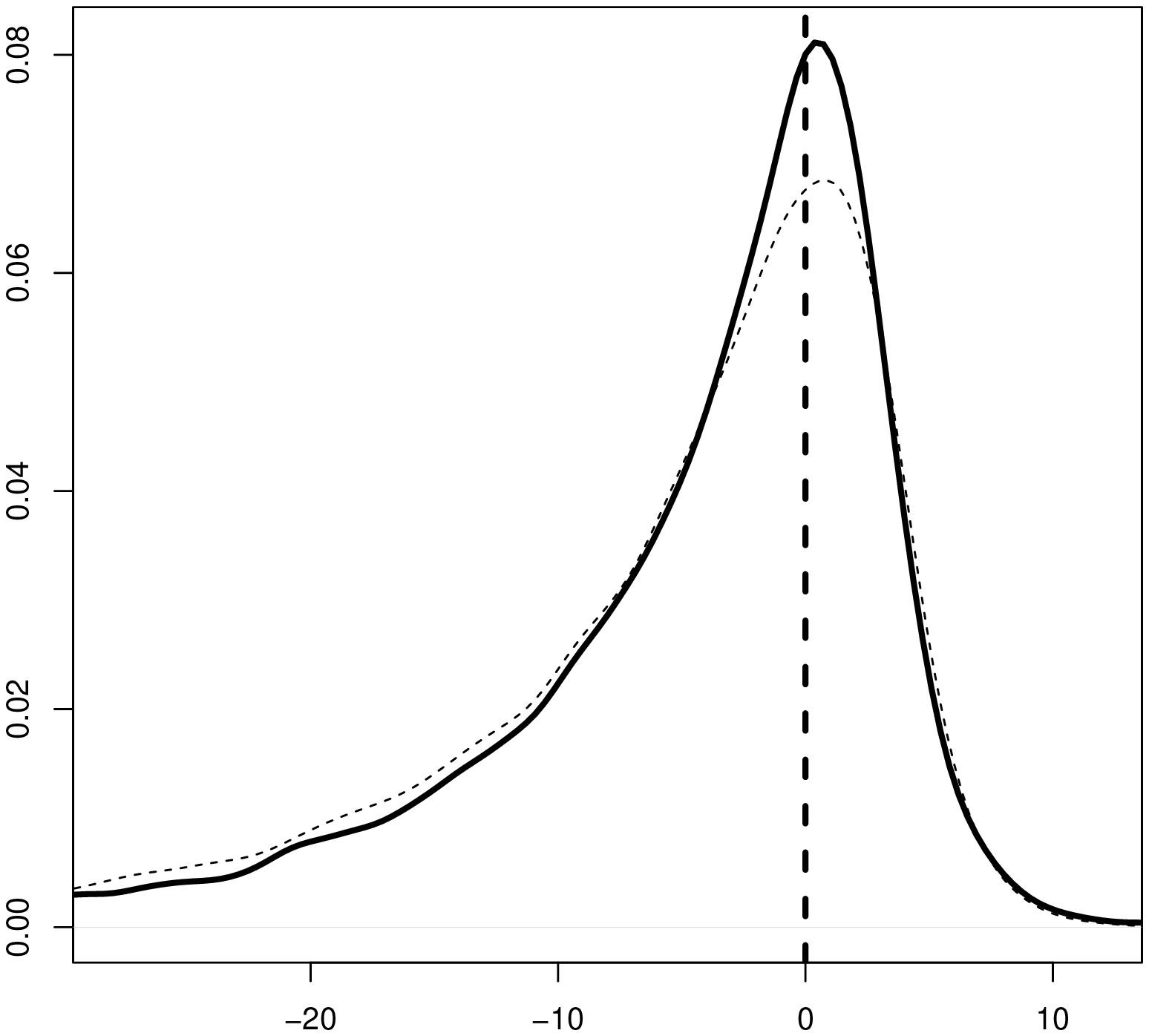}
\end{subfigure}\vspace{1.5cm}
\caption{{\small Distribution of $T(\hat{\beta}_1(\hat{\tau}_T)-\beta_1)$ when $\tau_0=0.3$, $\beta_2=0.75$, $\beta_1=1-1/T$. Left:
$\{\varepsilon_t\}_{t=1}^T\sim t(3)$; Right: $\{\varepsilon_t\}_{t=1}^T\sim t(2)$.}}
\end{figure}

\begin{figure}[H]
\label{fig13}
\begin{subfigure}[b]{0.5\textwidth}
\vspace{0.7cm}\centering
\includegraphics[trim=5cm 5cm 2cm 4cm, width=0.55\textwidth]{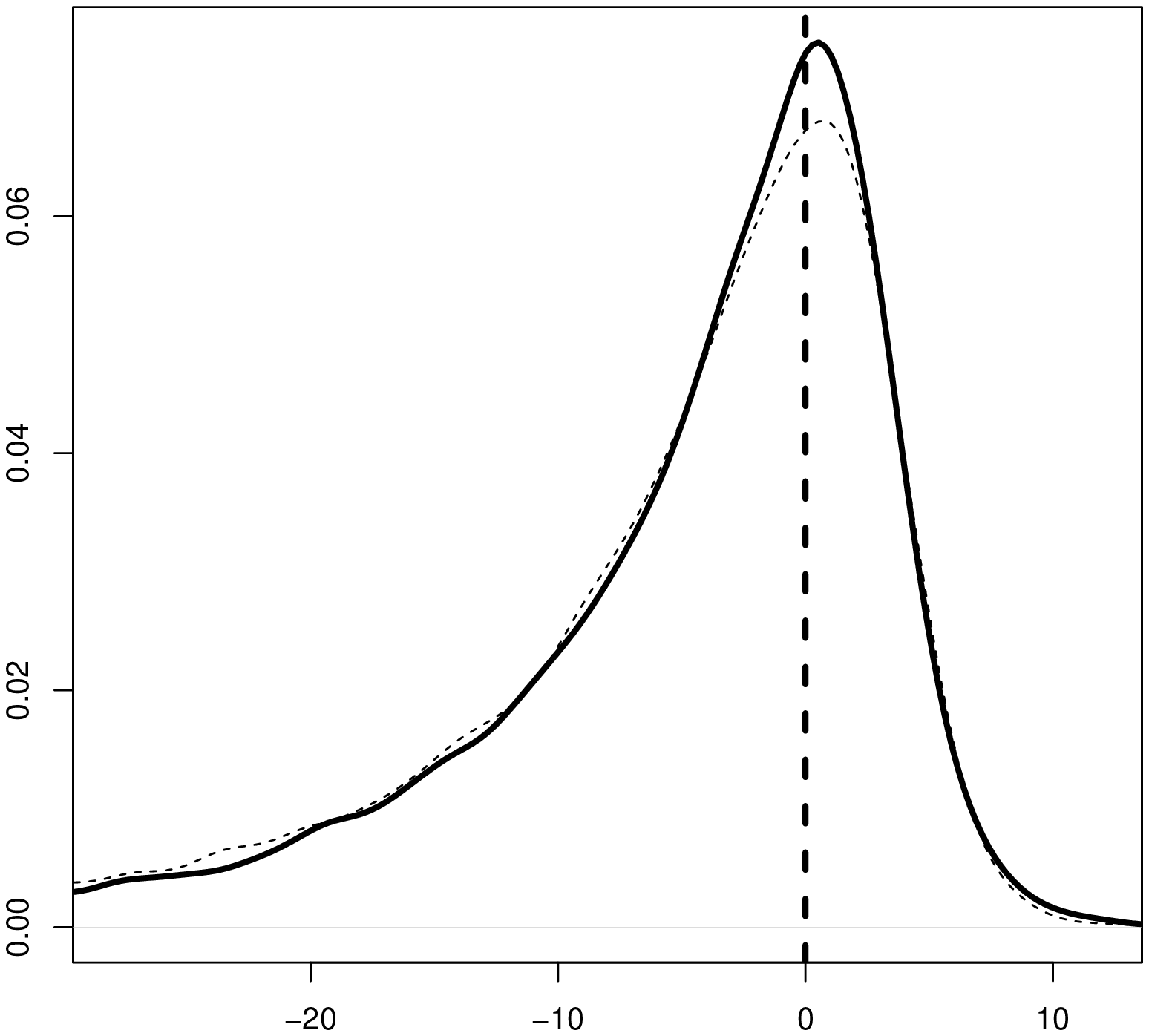}
\end{subfigure}\hspace{0.7cm}
\begin{subfigure}[b]{0.5\textwidth}
\vspace{0.7cm}\centering
\includegraphics[trim=5cm 5cm 2cm 4cm, width=0.55\textwidth]{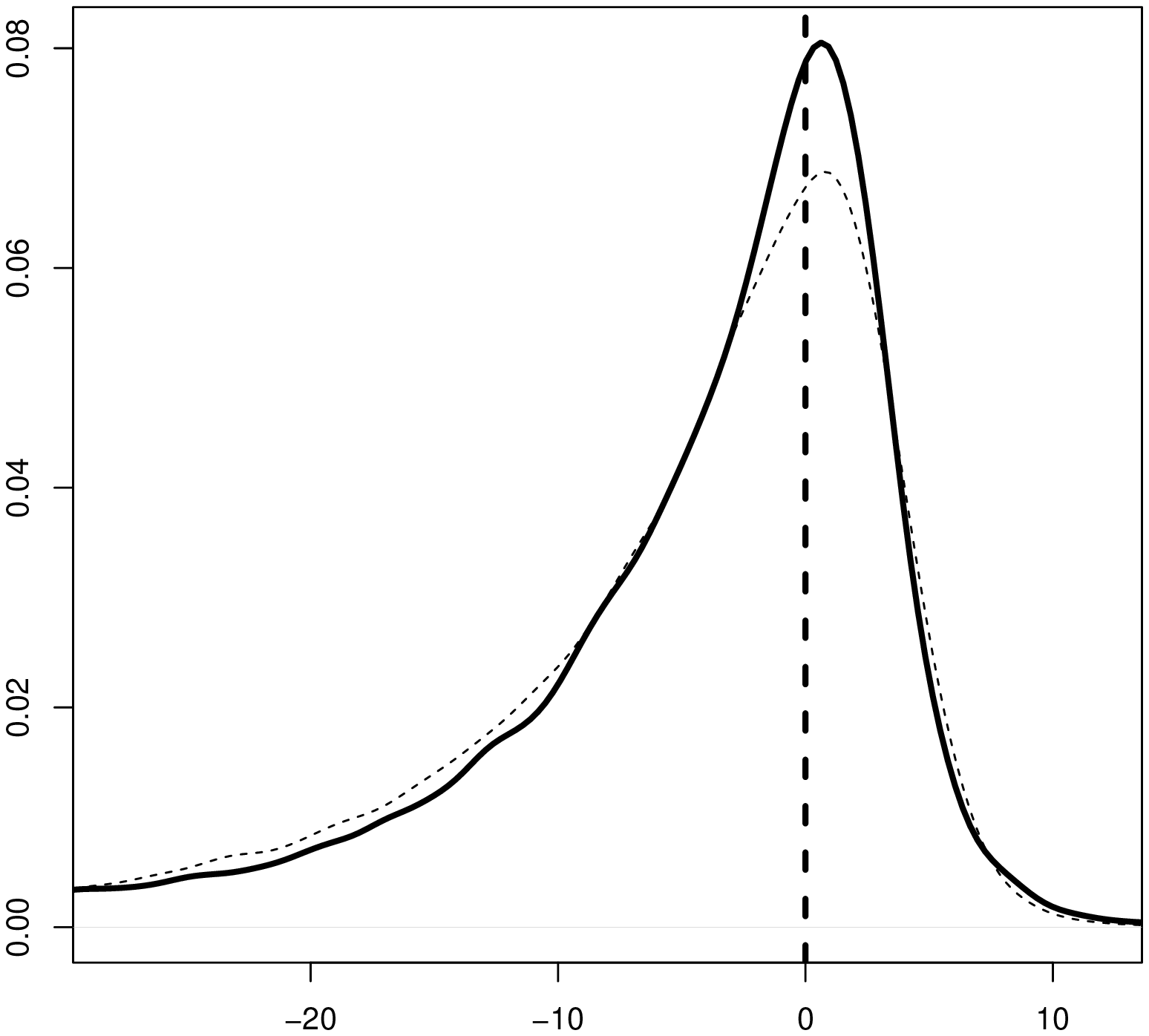}
\end{subfigure}\vspace{1.5cm}
\caption{{\small Distribution of $T(\hat{\beta}_1(\hat{\tau}_T)-\beta_1)$ when $\tau_0=0.3$, $\beta_2=0.8$, $\beta_1=1-1/T$. Left:
$\{\varepsilon_t\}_{t=1}^T\sim t(3)$; Right: $\{\varepsilon_t\}_{t=1}^T\sim t(2)$.}}
\end{figure}

\begin{figure}[H]
\label{fig13}
\begin{subfigure}[b]{0.5\textwidth}
\vspace{0.7cm}\centering
\includegraphics[trim=5cm 5cm 2cm 4cm, width=0.55\textwidth]{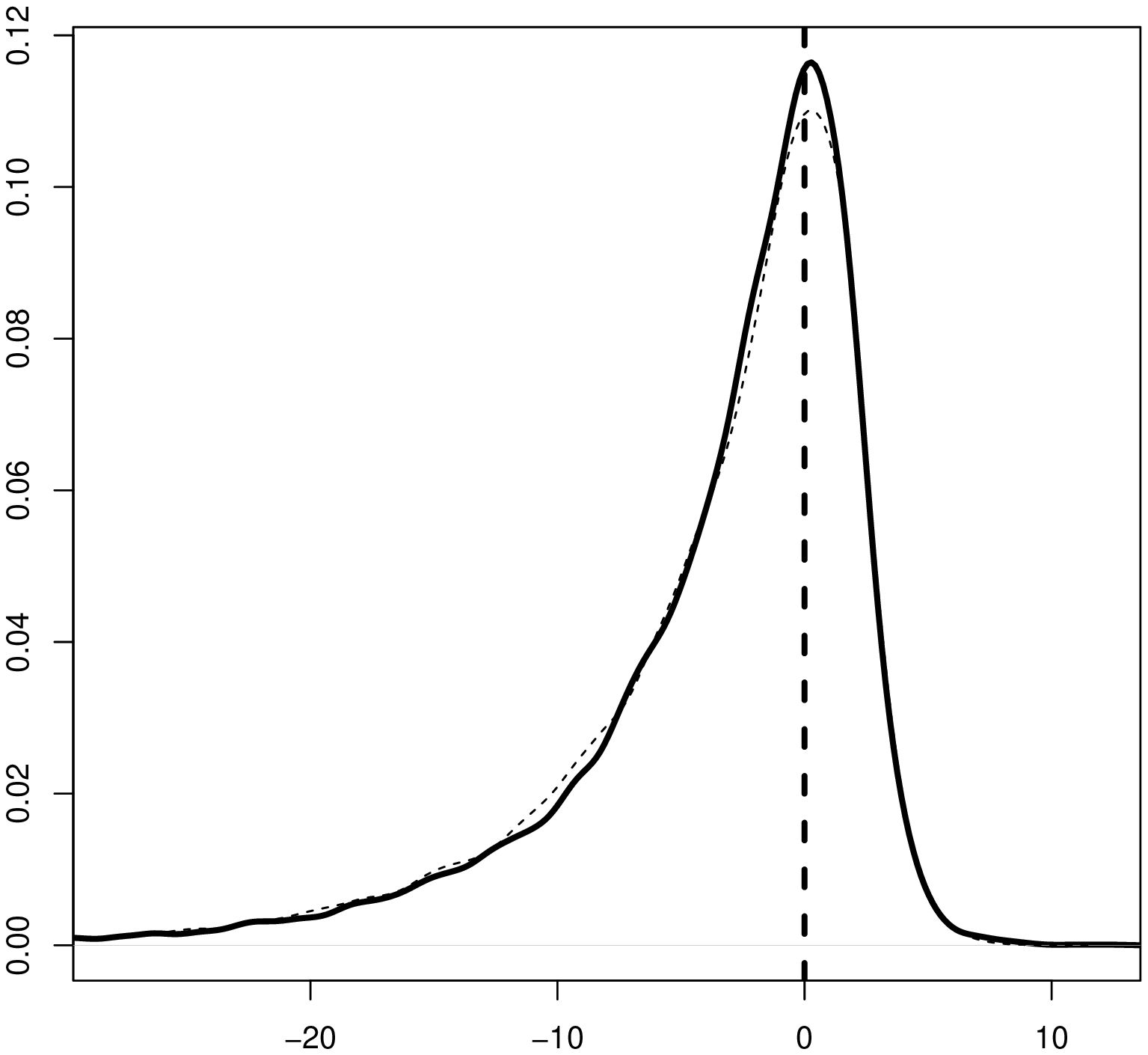}
\end{subfigure}\hspace{0.7cm}
\begin{subfigure}[b]{0.5\textwidth}
\vspace{0.7cm}\centering
\includegraphics[trim=5cm 5cm 2cm 4cm, width=0.55\textwidth]{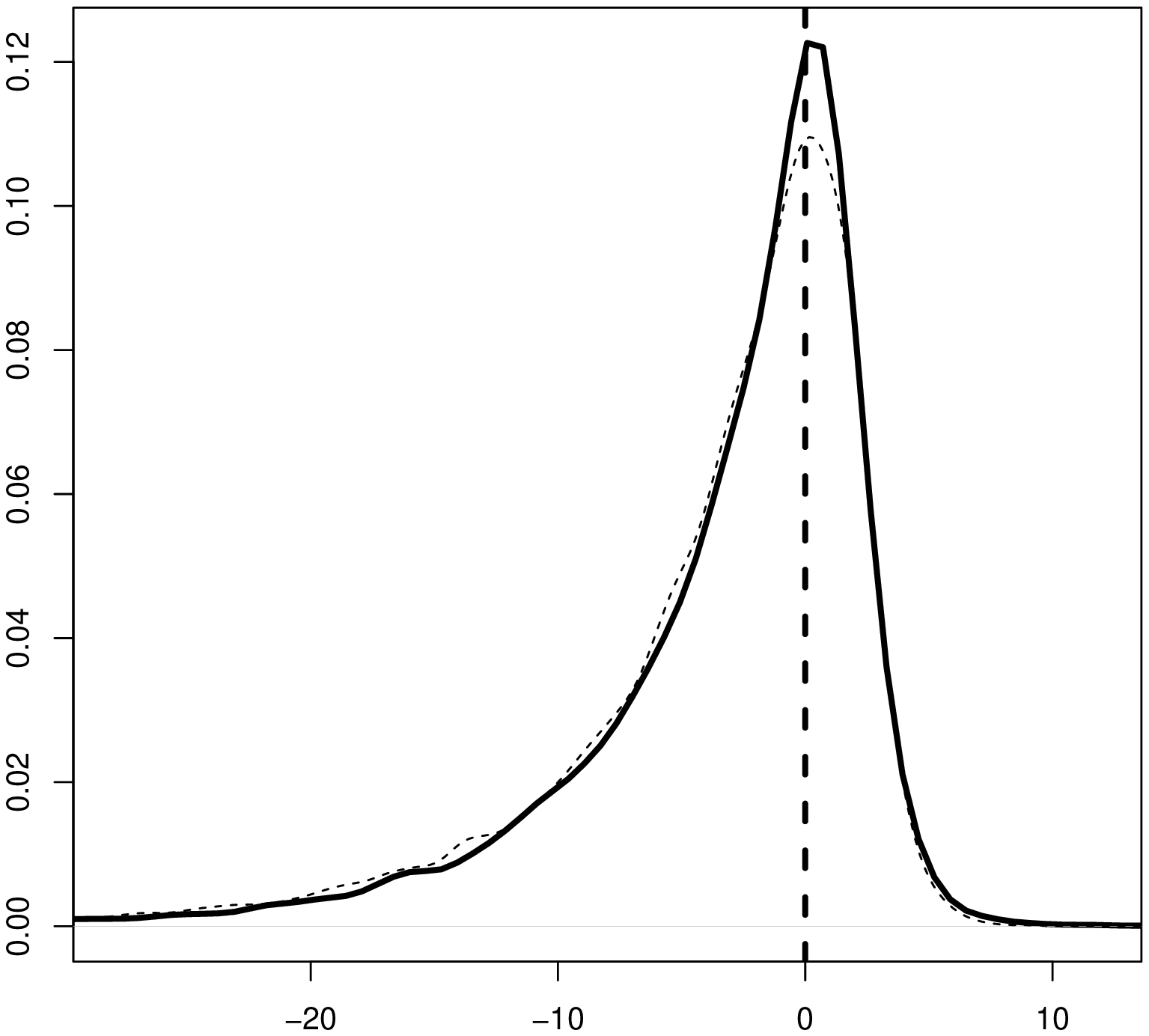}
\end{subfigure}\vspace{1.5cm}
\caption{{\small Distribution of $T(\hat{\beta}_1(\hat{\tau}_T)-\beta_1)$ when $\tau_0=0.5$, $\beta_2=0.5$, $\beta_1=1-1/T$. Left:
$\{\varepsilon_t\}_{t=1}^T\sim t(3)$; Right: $\{\varepsilon_t\}_{t=1}^T\sim t(2)$.}}
\end{figure}

\begin{figure}[H]
\label{fig13}
\begin{subfigure}[b]{0.5\textwidth}
\vspace{0.7cm}\centering
\includegraphics[trim=5cm 5cm 2cm 4cm, width=0.55\textwidth]{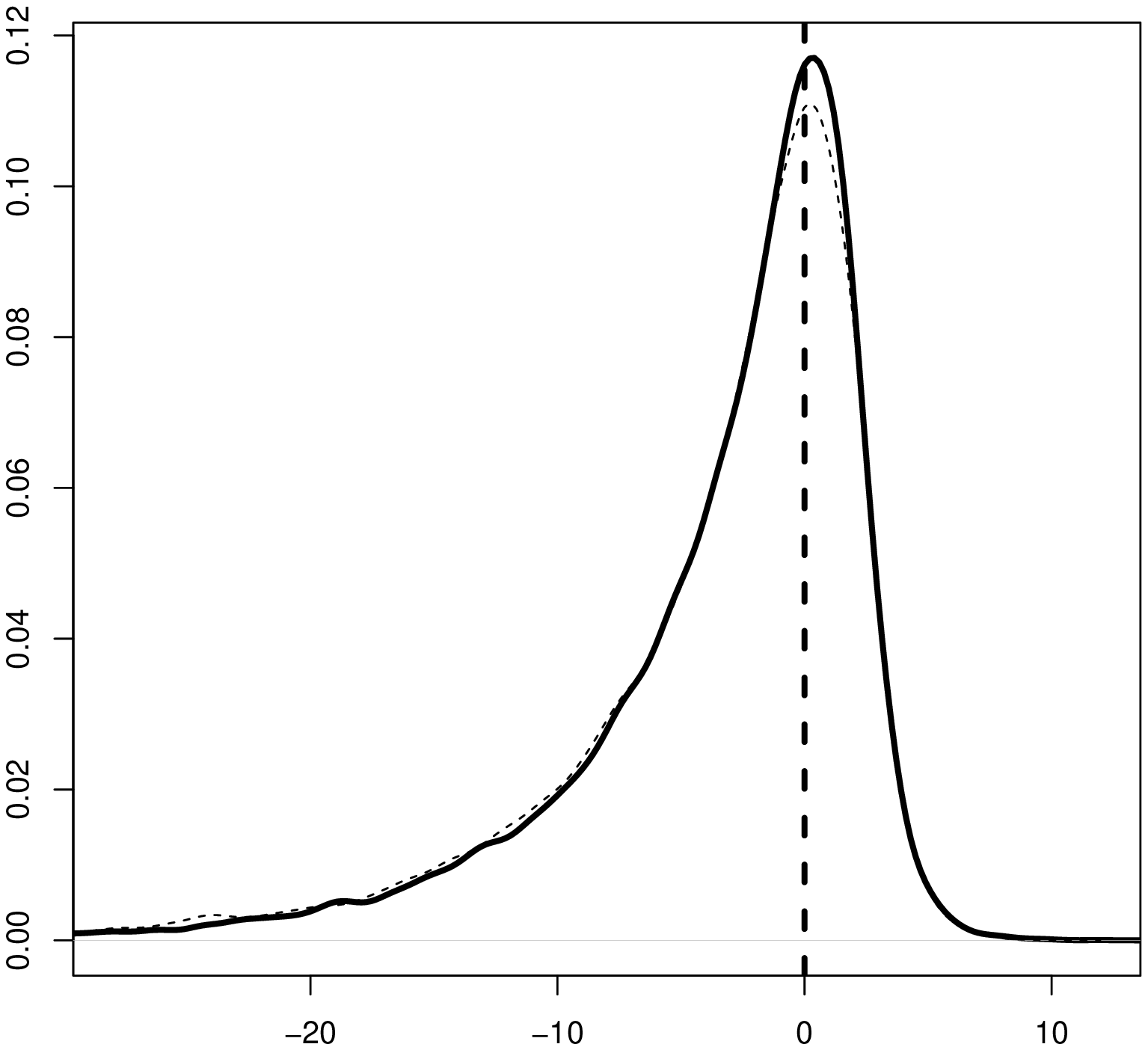}
\end{subfigure}\hspace{0.7cm}
\begin{subfigure}[b]{0.5\textwidth}
\vspace{0.7cm}\centering
\includegraphics[trim=5cm 5cm 2cm 4cm, width=0.55\textwidth]{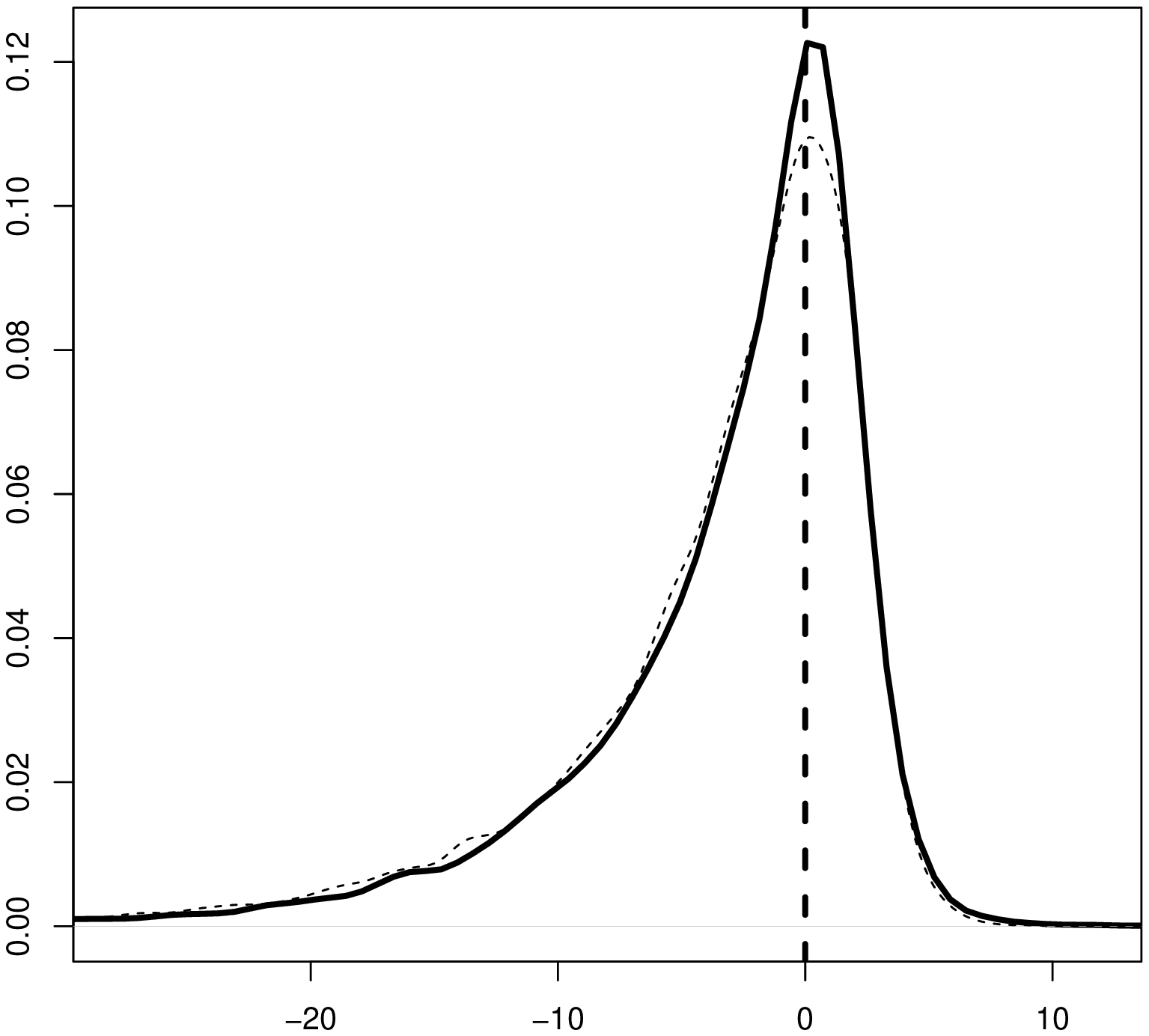}
\end{subfigure}\vspace{1.5cm}
\caption{{\small Distribution of $T(\hat{\beta}_1(\hat{\tau}_T)-\beta_1)$ when $\tau_0=0.5$, $\beta_2=0.75$, $\beta_1=1-1/T$. Left:
$\{\varepsilon_t\}_{t=1}^T\sim t(3)$; Right: $\{\varepsilon_t\}_{t=1}^T\sim t(2)$.}}
\end{figure}

\begin{figure}[H]
\label{fig13}
\begin{subfigure}[b]{0.5\textwidth}
\vspace{0.7cm}\centering
\includegraphics[trim=5cm 5cm 2cm 4cm, width=0.55\textwidth]{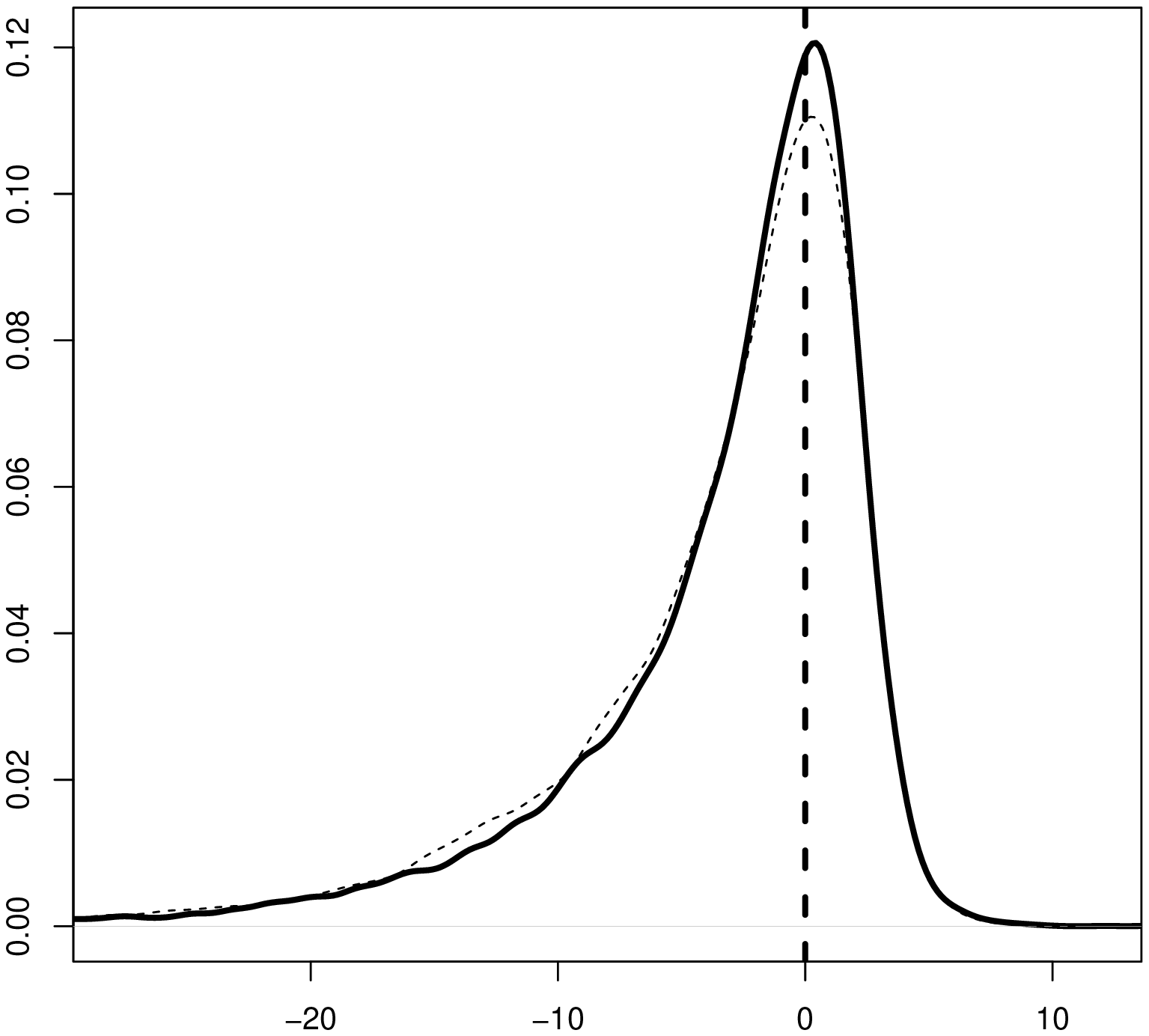}
\end{subfigure}\hspace{0.7cm}
\begin{subfigure}[b]{0.5\textwidth}
\vspace{0.7cm}\centering
\includegraphics[trim=5cm 5cm 2cm 4cm, width=0.55\textwidth]{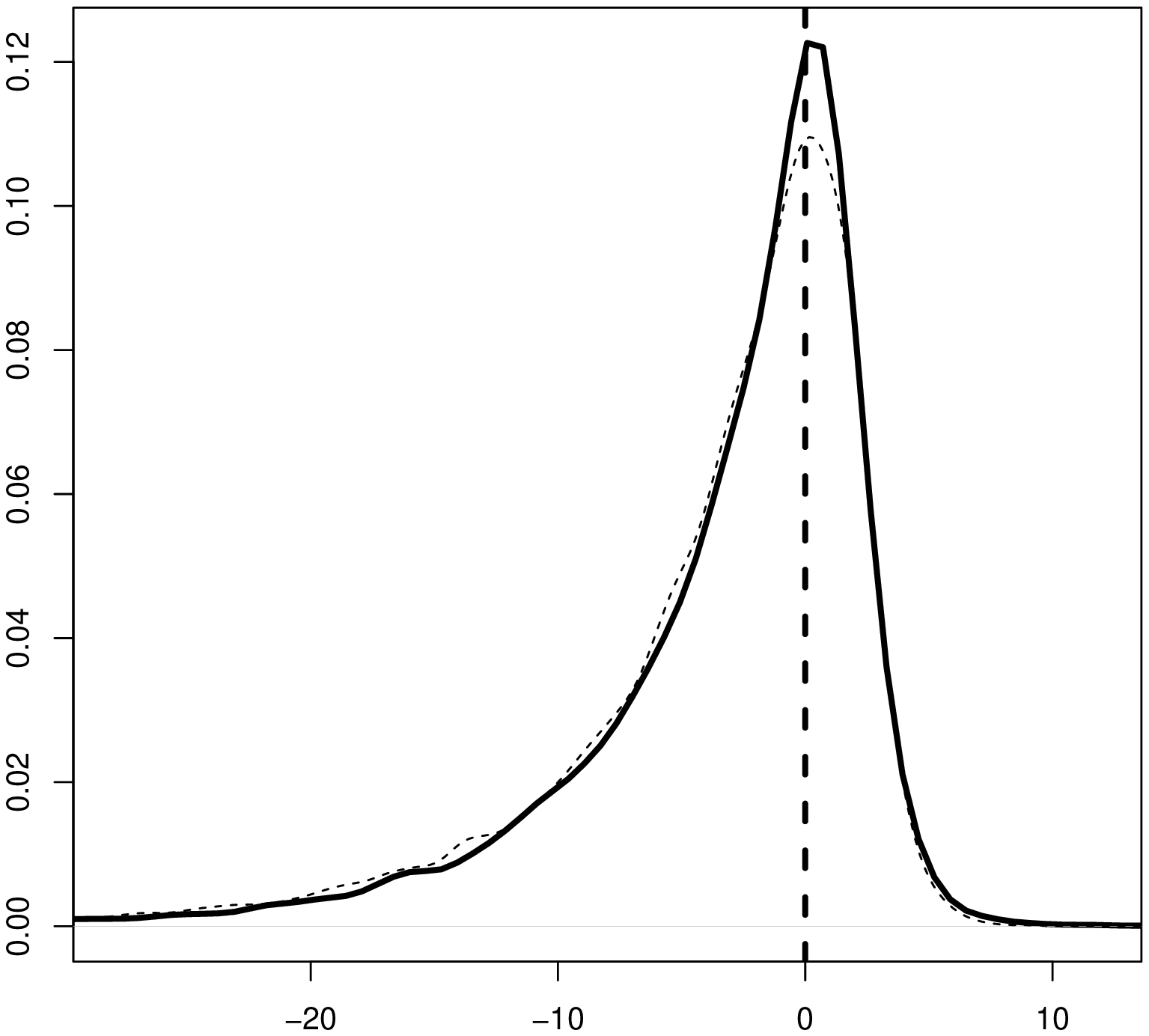}
\end{subfigure}\vspace{1.5cm}
\caption{{\small Distribution of $T(\hat{\beta}_1(\hat{\tau}_T)-\beta_1)$ when $\tau_0=0.5$, $\beta_2=0.8$, $\beta_1=1-1/T$. Left:
$\{\varepsilon_t\}_{t=1}^T\sim t(3)$; Right: $\{\varepsilon_t\}_{t=1}^T\sim t(2)$.}}
\end{figure}

\begin{figure}[H]
\label{fig19}
\begin{subfigure}[b]{0.5\textwidth}
\vspace{0.7cm}\centering
\includegraphics[trim=5cm 5cm 2cm 4cm, width=0.55\textwidth]{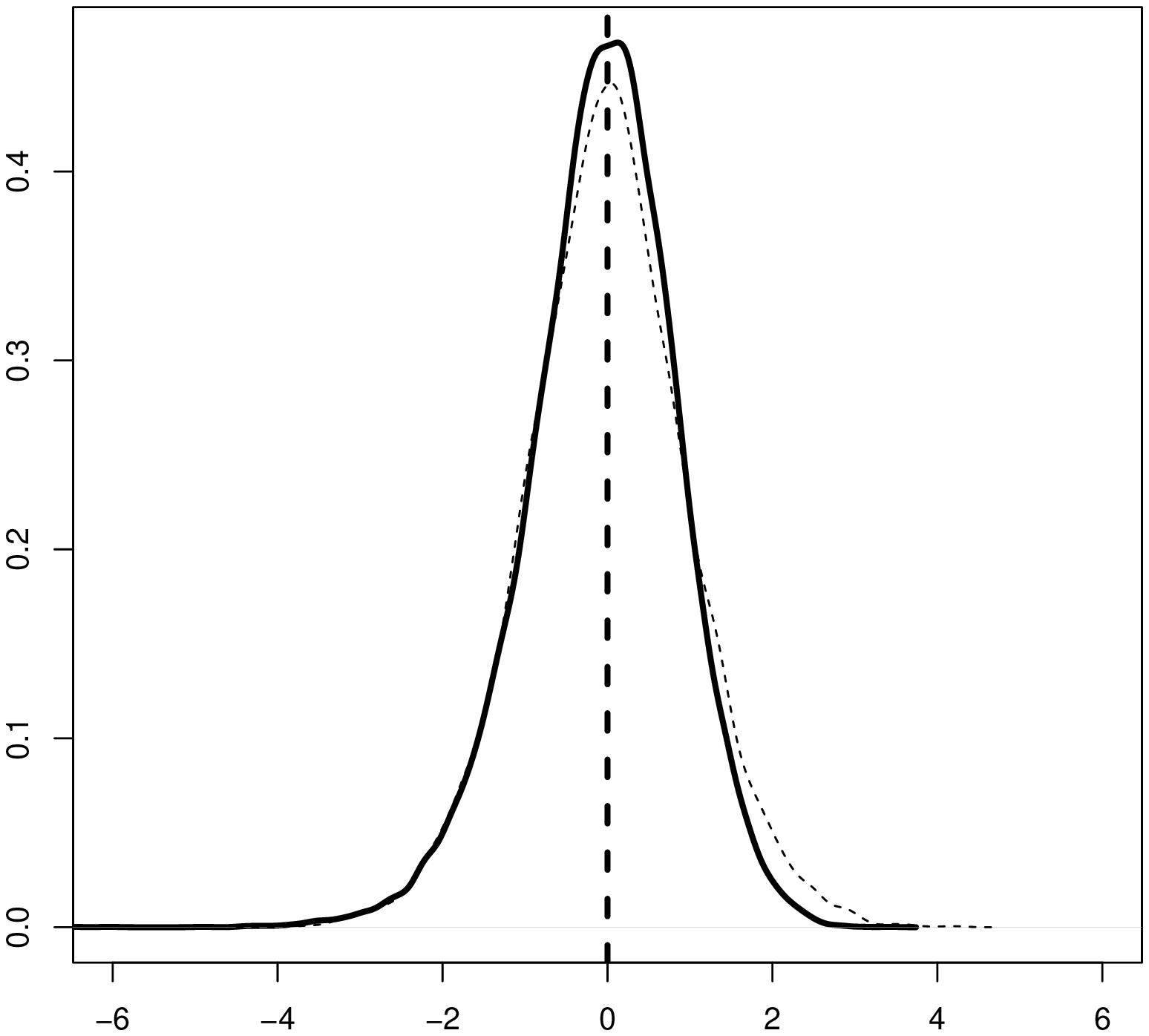}
\end{subfigure}\hspace{0.7cm}
\begin{subfigure}[b]{0.5\textwidth}
\vspace{0.7cm}\centering
\includegraphics[trim=5cm 5cm 2cm 4cm, width=0.55\textwidth]{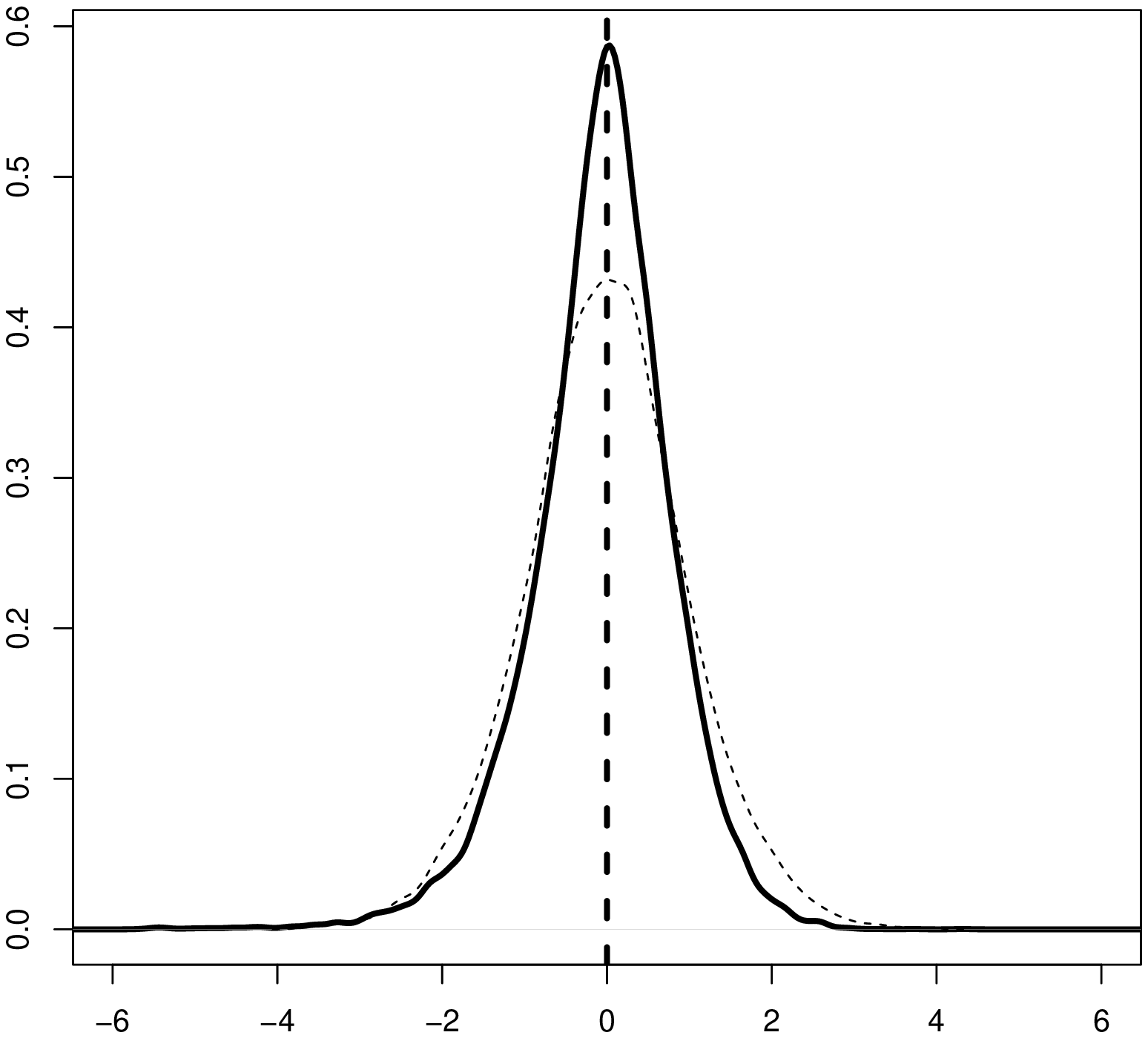}
\end{subfigure}\vspace{1.5cm}
\caption{{\small Distribution of $\sqrt{T}(\hat{\protect\beta}_2(\hat{\tau}_T)-\beta_2)$ when $\tau_0=0.3$, $\beta_2=0.5$, $\beta_1=1-1/T$. Left: $\{\varepsilon_t\}_{t=1}^T\sim t(3)$; Right: $\{\varepsilon_t\}_{t=1}^T\sim t(2)$.}}
\end{figure}

\begin{figure}[H]
\label{fig19}
\begin{subfigure}[b]{0.5\textwidth}
\vspace{0.7cm}\centering
\includegraphics[trim=5cm 5cm 2cm 4cm, width=0.55\textwidth]{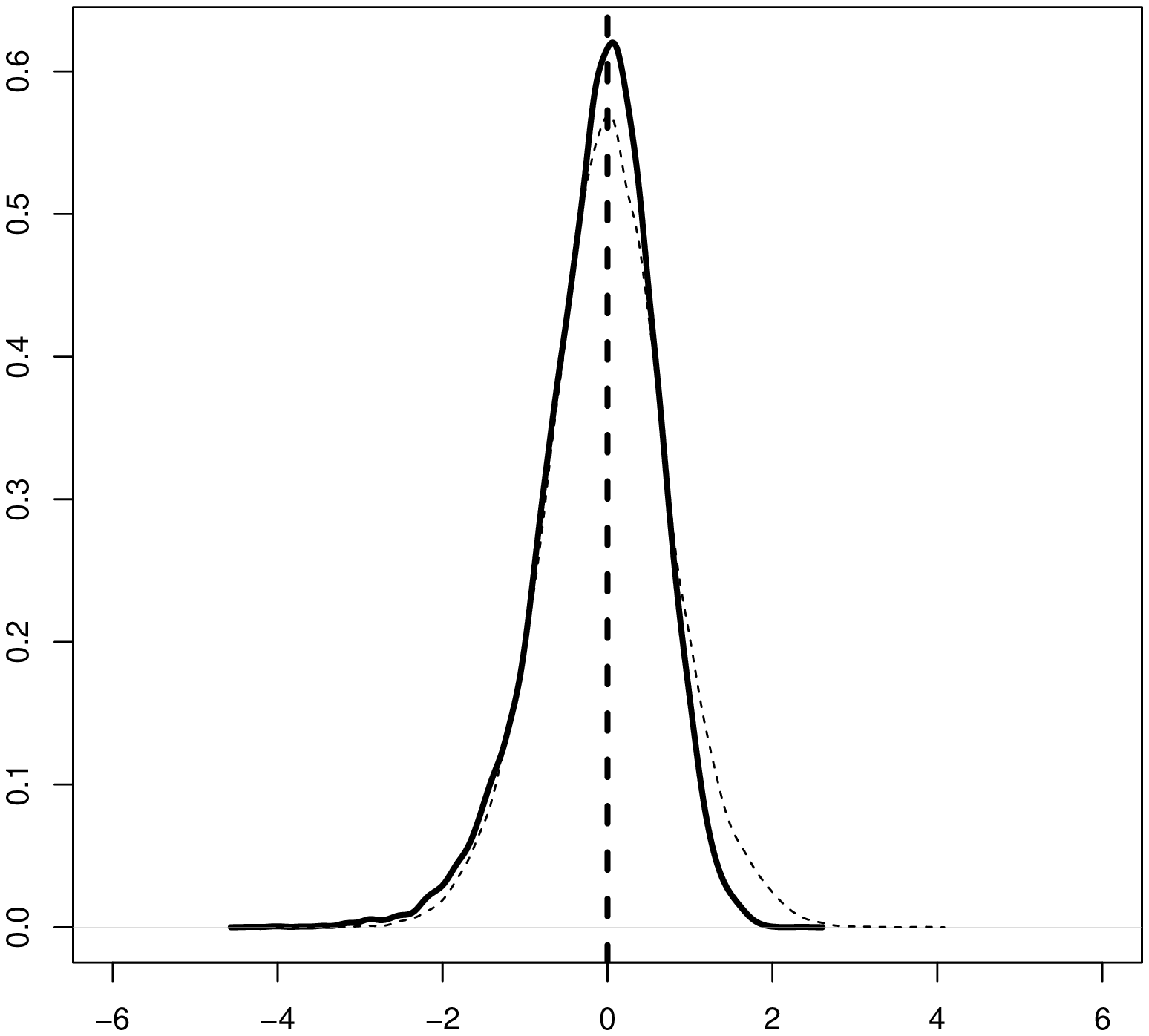}
\end{subfigure}\hspace{0.7cm}
\begin{subfigure}[b]{0.5\textwidth}
\vspace{0.7cm}\centering
\includegraphics[trim=5cm 5cm 2cm 4cm, width=0.55\textwidth]{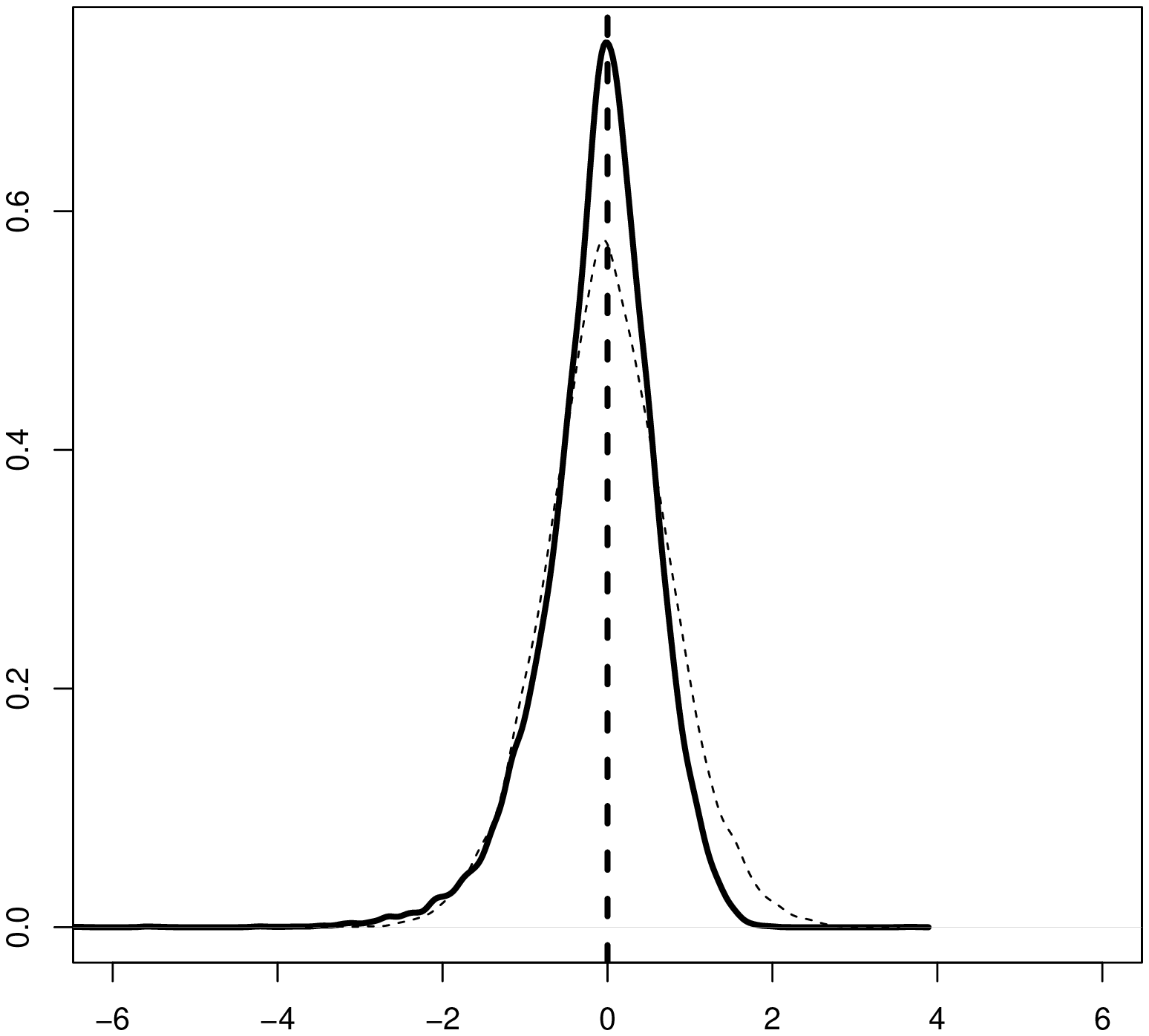}
\end{subfigure}\vspace{1.5cm}
\caption{{\small Distribution of $\sqrt{T}(\hat{\protect\beta}_2(\hat{\tau}_T)-\beta_2)$ when $\tau_0=0.3$, $\beta_2=0.75$, $\beta_1=1-1/T$. Left: $\{\varepsilon_t\}_{t=1}^T\sim t(3)$; Right: $\{\varepsilon_t\}_{t=1}^T\sim t(2)$.}}
\end{figure}

\begin{figure}[H]
\label{fig19}
\begin{subfigure}[b]{0.5\textwidth}
\vspace{0.7cm}\centering
\includegraphics[trim=5cm 5cm 2cm 4cm, width=0.55\textwidth]{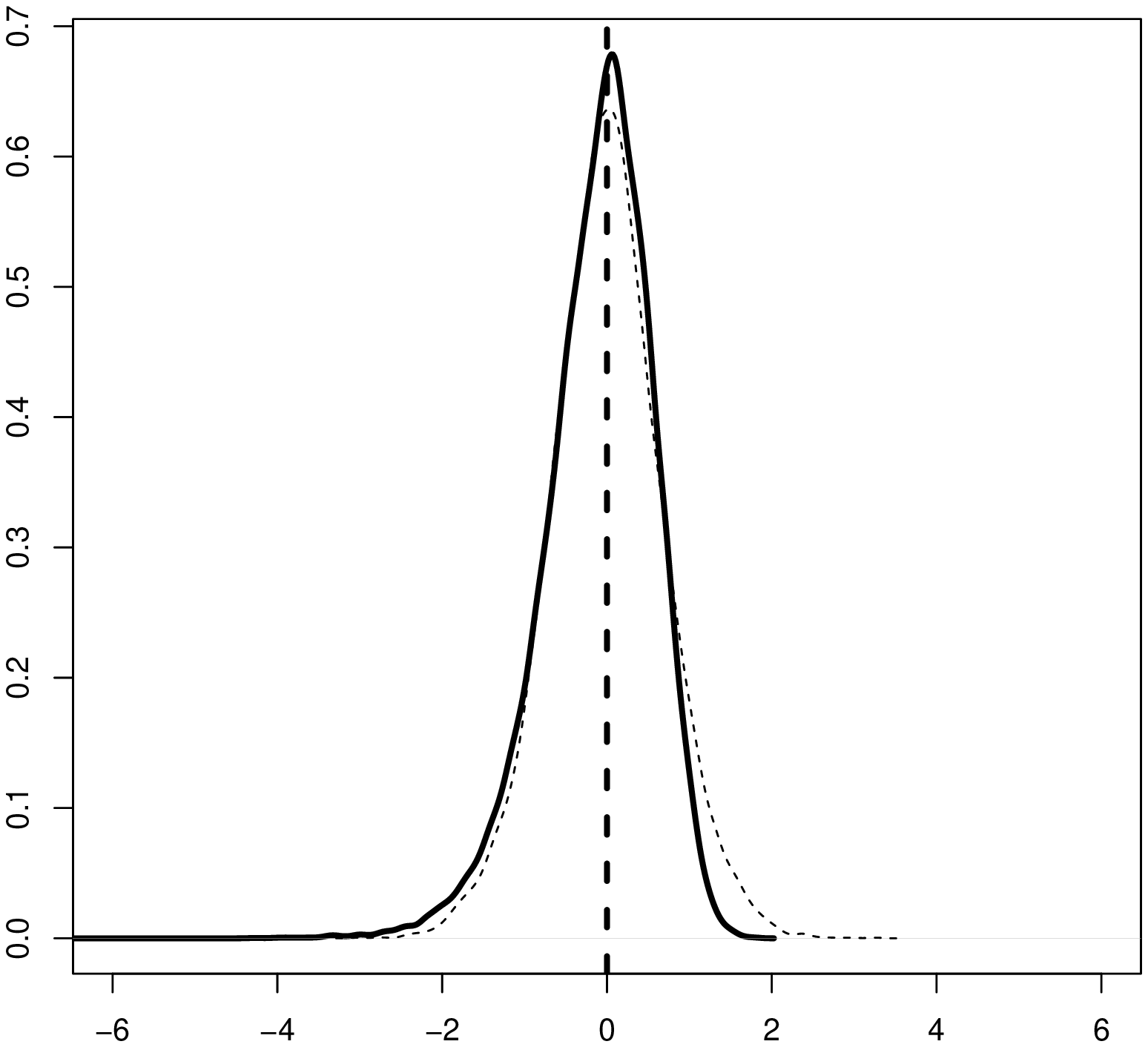}
\end{subfigure}\hspace{0.7cm}
\begin{subfigure}[b]{0.5\textwidth}
\vspace{0.7cm}\centering
\includegraphics[trim=5cm 5cm 2cm 4cm, width=0.55\textwidth]{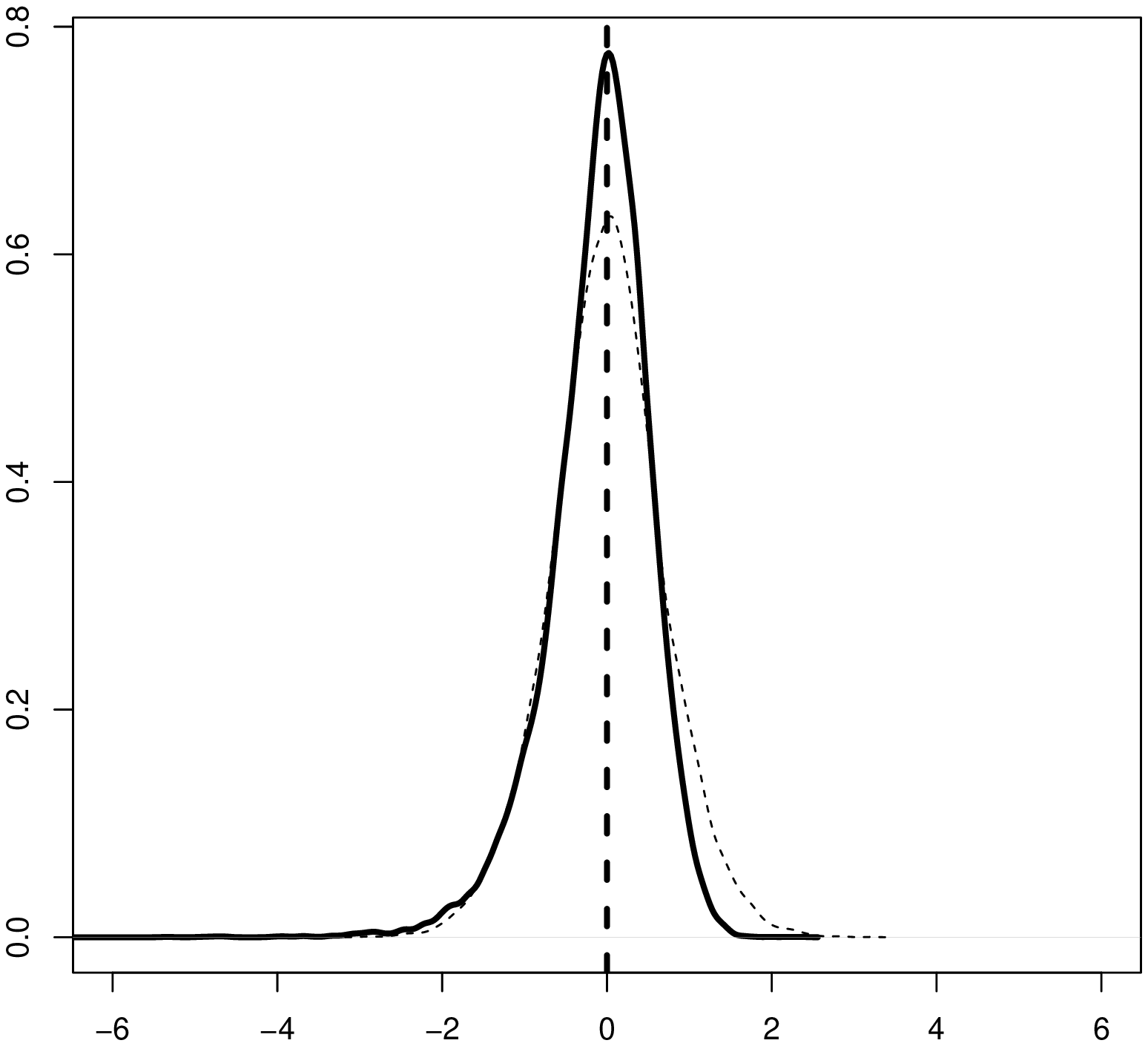}
\end{subfigure}\vspace{1.5cm}
\caption{{\small Distribution of $\sqrt{T}(\hat{\protect\beta}_2(\hat{\tau}_T)-\beta_2)$ when $\tau_0=0.3$, $\beta_2=0.8$, $\beta_1=1-1/T$. Left: $\{\varepsilon_t\}_{t=1}^T\sim t(3)$; Right: $\{\varepsilon_t\}_{t=1}^T\sim t(2)$.}}
\end{figure}

\begin{figure}[H]
\label{fig19}
\begin{subfigure}[b]{0.5\textwidth}
\vspace{0.7cm}\centering
\includegraphics[trim=5cm 5cm 2cm 4cm, width=0.55\textwidth]{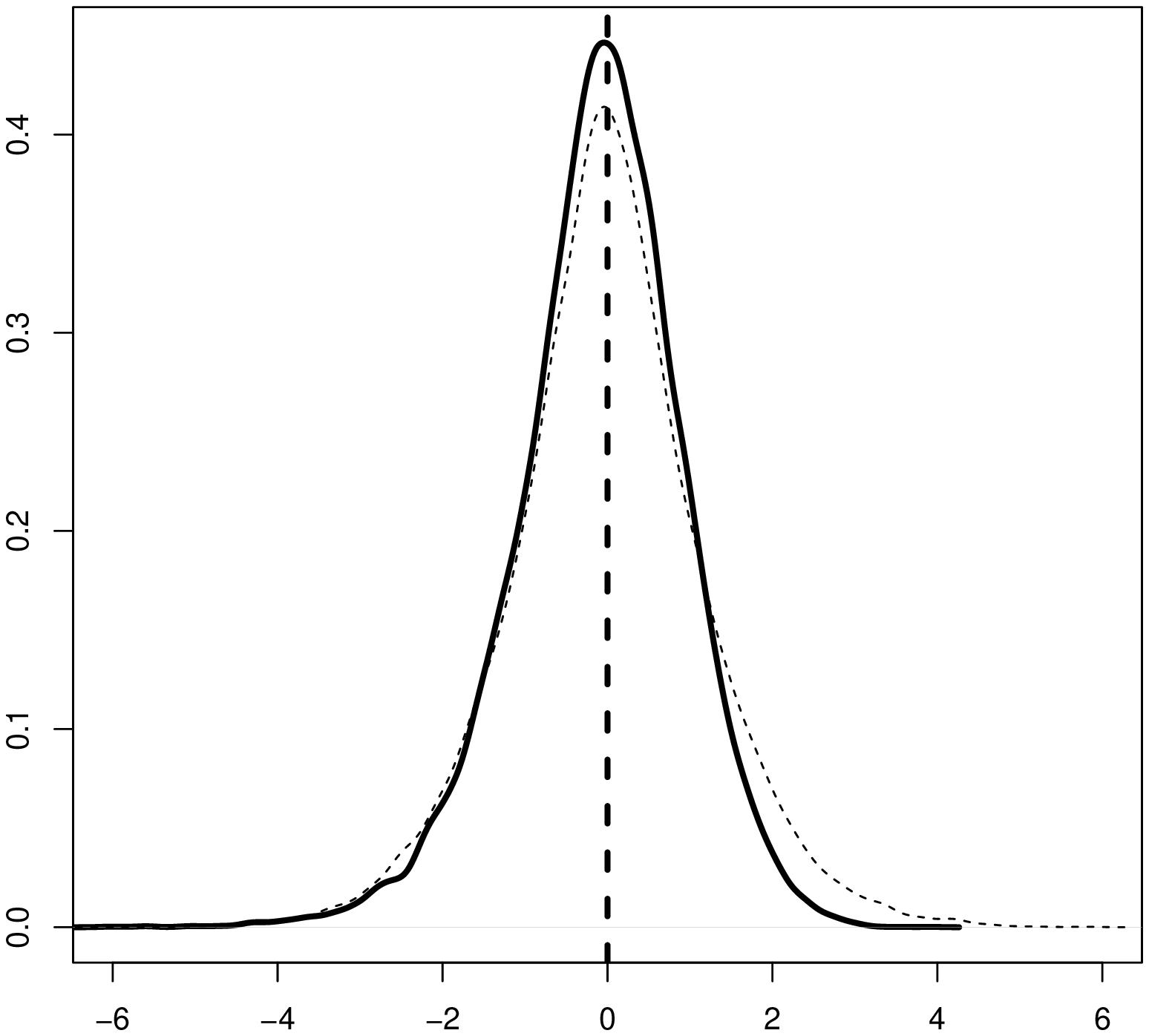}
\end{subfigure}\hspace{0.7cm}
\begin{subfigure}[b]{0.5\textwidth}
\vspace{0.7cm}\centering
\includegraphics[trim=5cm 5cm 2cm 4cm, width=0.55\textwidth]{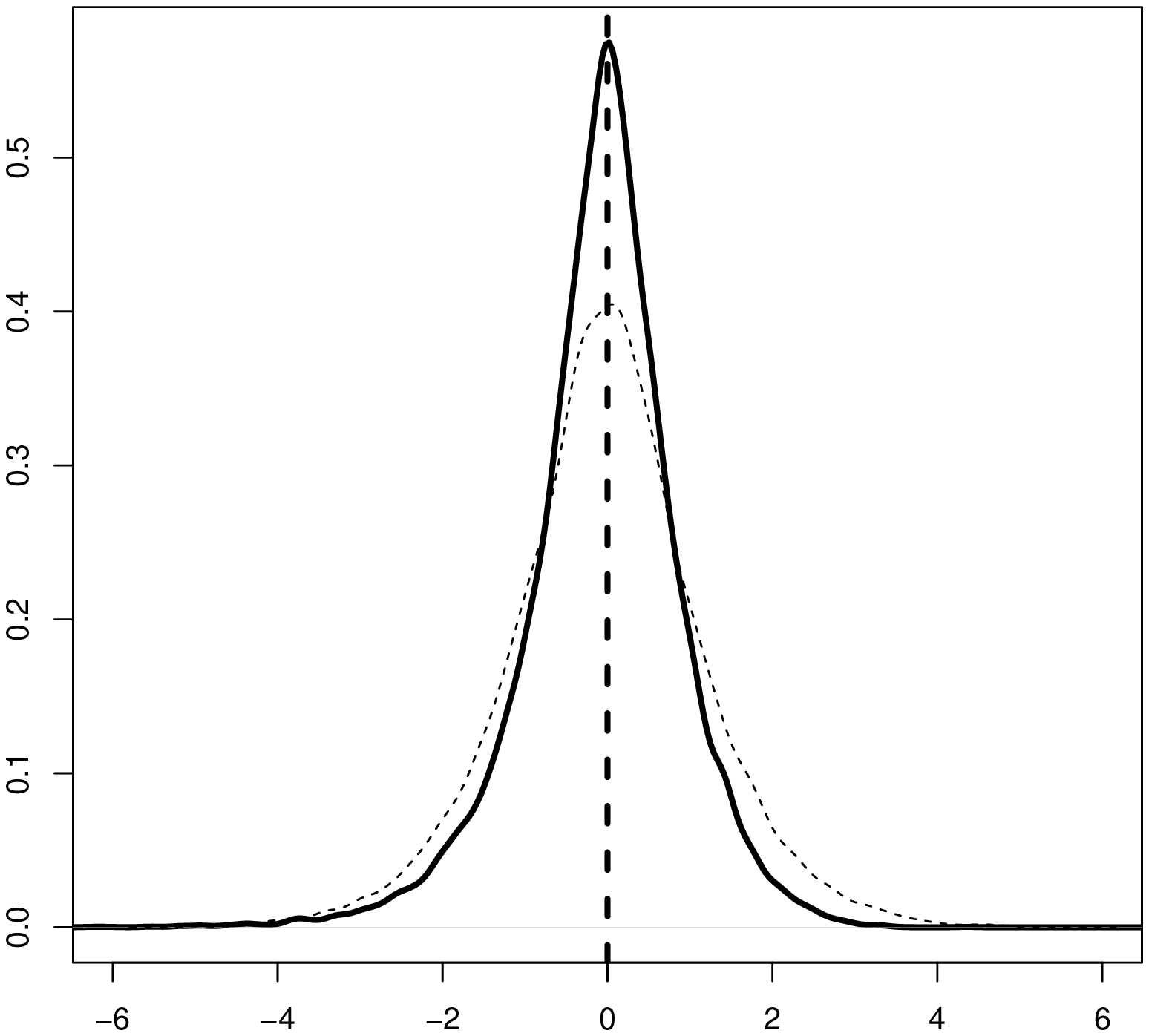}
\end{subfigure}\vspace{1.5cm}
\caption{{\small Distribution of $\sqrt{T}(\hat{\protect\beta}_2(\hat{\tau}_T)-\beta_2)$ when $\tau_0=0.5$, $\beta_2=0.5$, $\beta_1=1-1/T$. Left: $\{\varepsilon_t\}_{t=1}^T\sim t(3)$; Right: $\{\varepsilon_t\}_{t=1}^T\sim t(2)$.}}
\end{figure}

\begin{figure}[H]
\label{fig19}
\begin{subfigure}[b]{0.5\textwidth}
\vspace{0.7cm}\centering
\includegraphics[trim=5cm 5cm 2cm 4cm, width=0.55\textwidth]{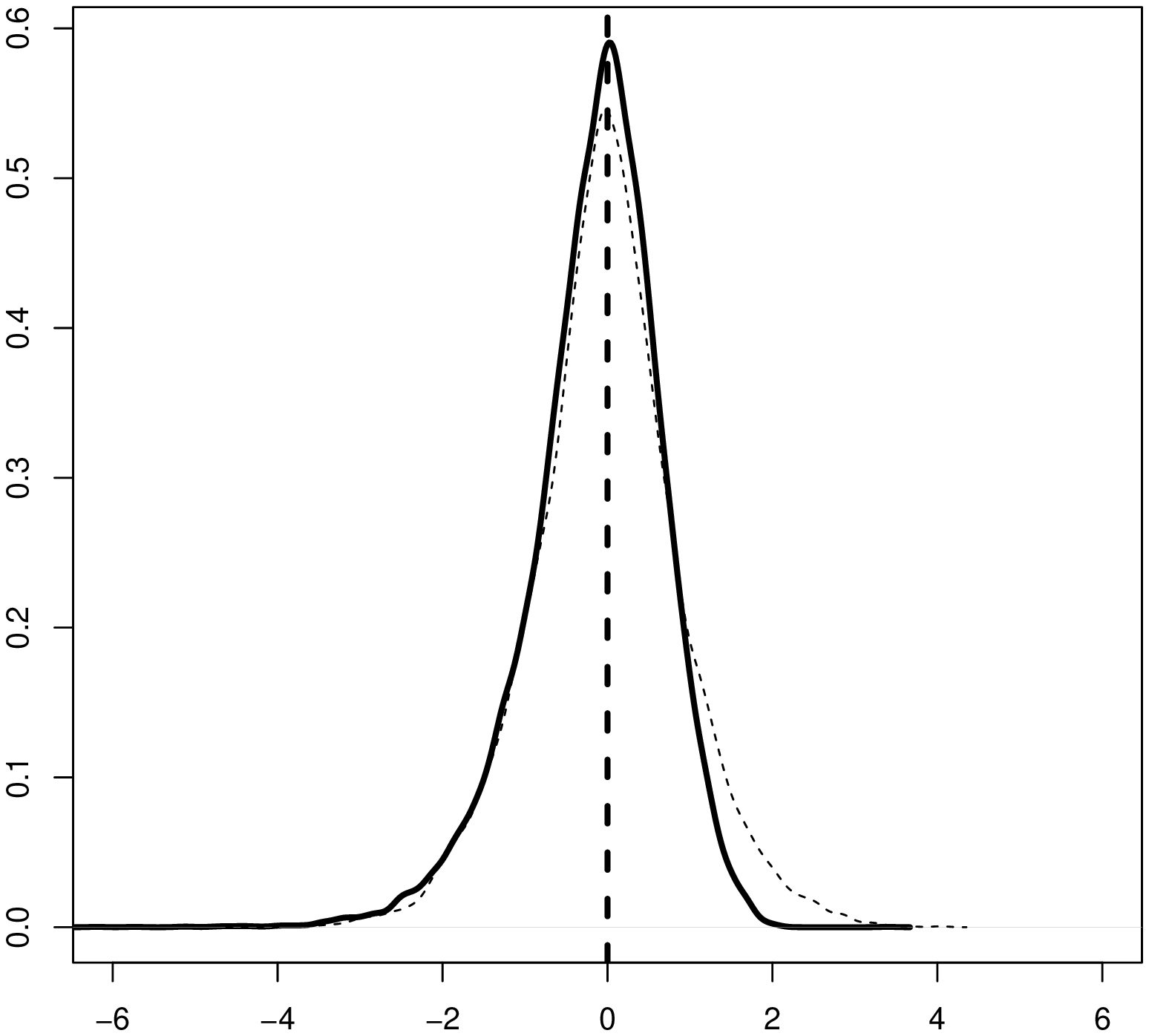}
\end{subfigure}\hspace{0.7cm}
\begin{subfigure}[b]{0.5\textwidth}
\vspace{0.7cm}\centering
\includegraphics[trim=5cm 5cm 2cm 4cm, width=0.55\textwidth]{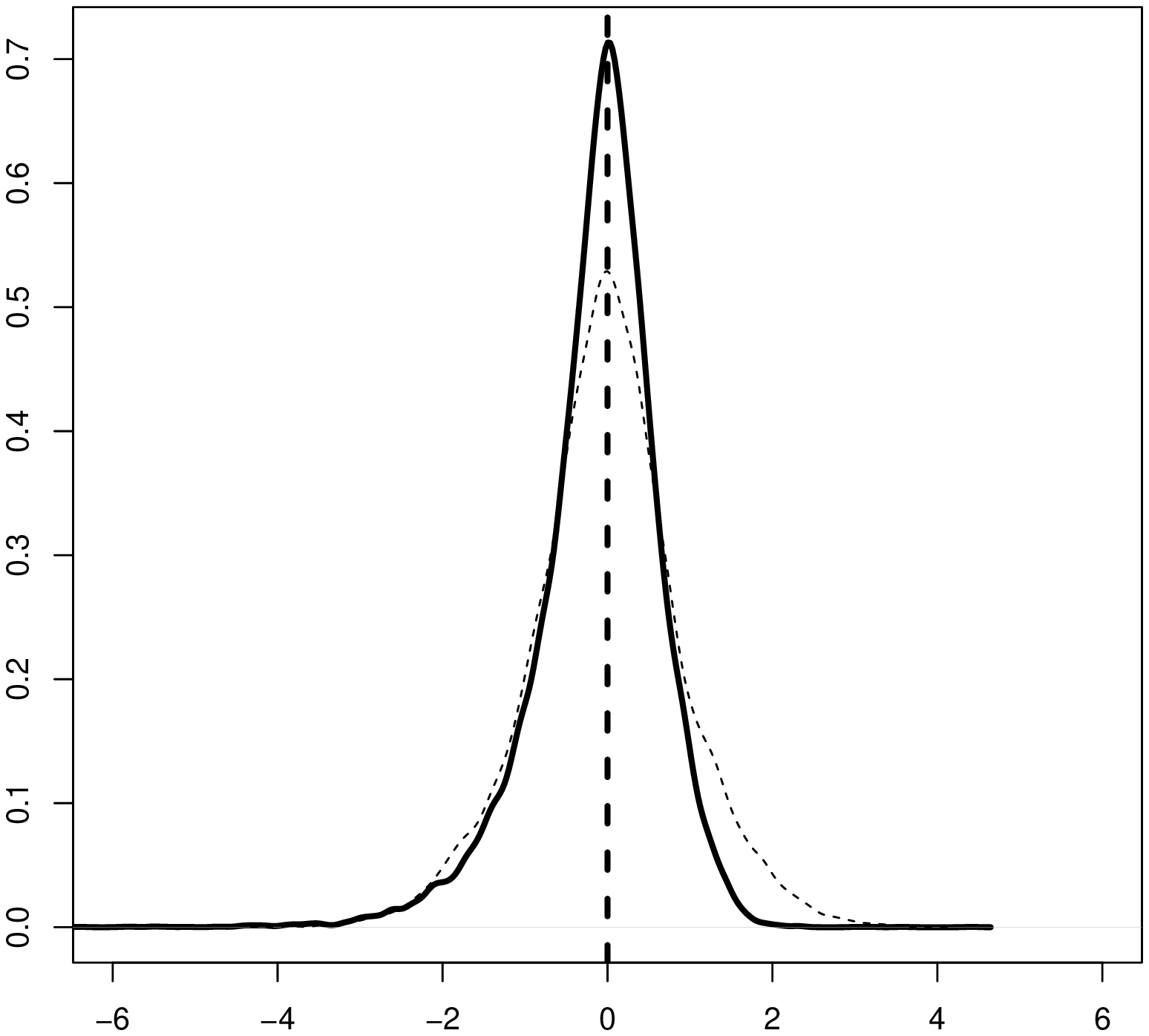}
\end{subfigure}\vspace{1.5cm}
\caption{{\small Distribution of $\sqrt{T}(\hat{\protect\beta}_2(\hat{\tau}_T)-\beta_2)$ when $\tau_0=0.5$, $\beta_2=0.75$, $\beta_1=1-1/T$. Left: $\{\varepsilon_t\}_{t=1}^T\sim t(3)$; Right: $\{\varepsilon_t\}_{t=1}^T\sim t(2)$.}}
\end{figure}

\begin{figure}[H]
\label{fig19}
\begin{subfigure}[b]{0.5\textwidth}
\vspace{0.7cm}\centering
\includegraphics[trim=5cm 5cm 2cm 4cm, width=0.55\textwidth]{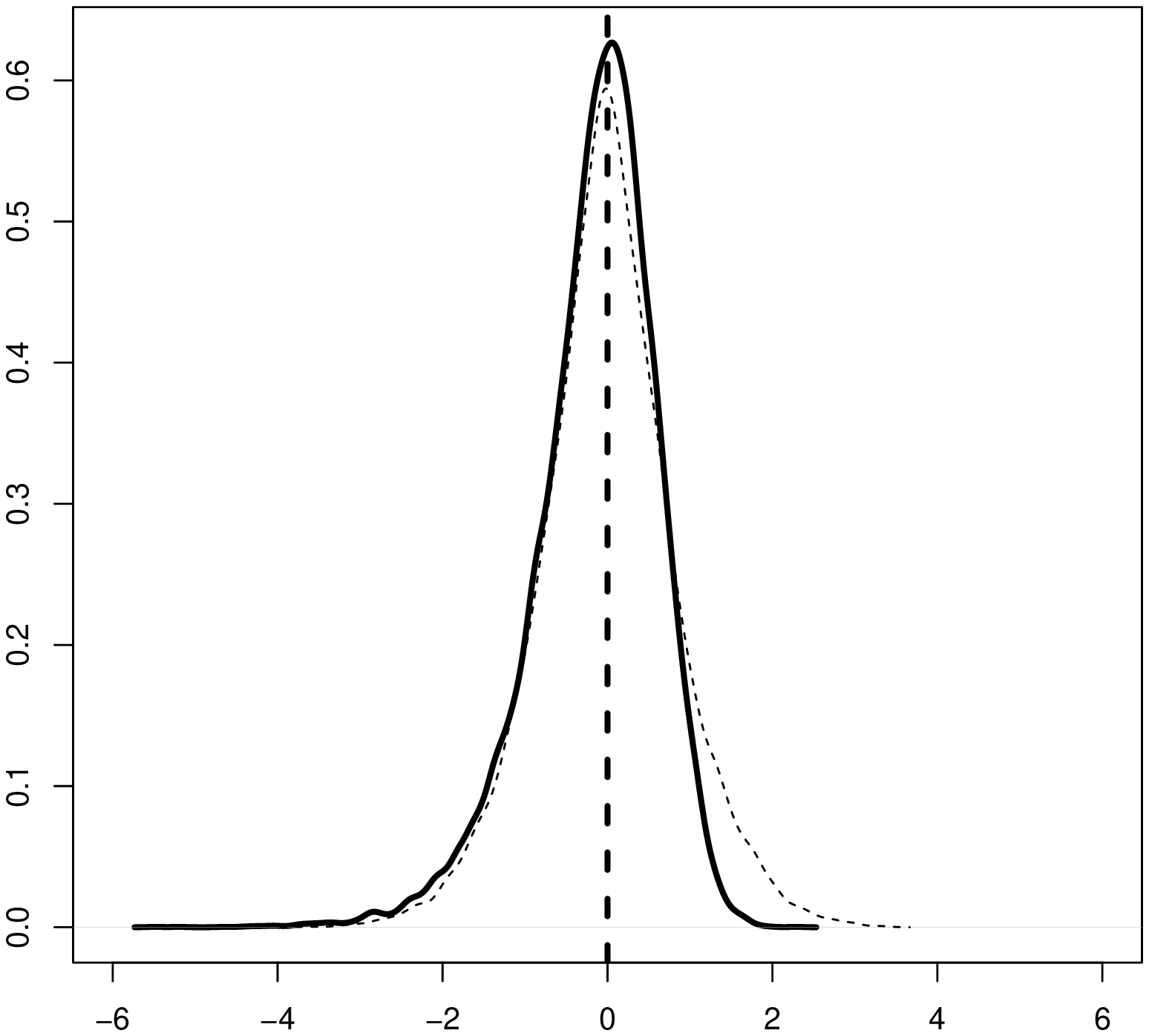}
\end{subfigure}\hspace{0.7cm}
\begin{subfigure}[b]{0.5\textwidth}
\vspace{0.7cm}\centering
\includegraphics[trim=5cm 5cm 2cm 4cm, width=0.55\textwidth]{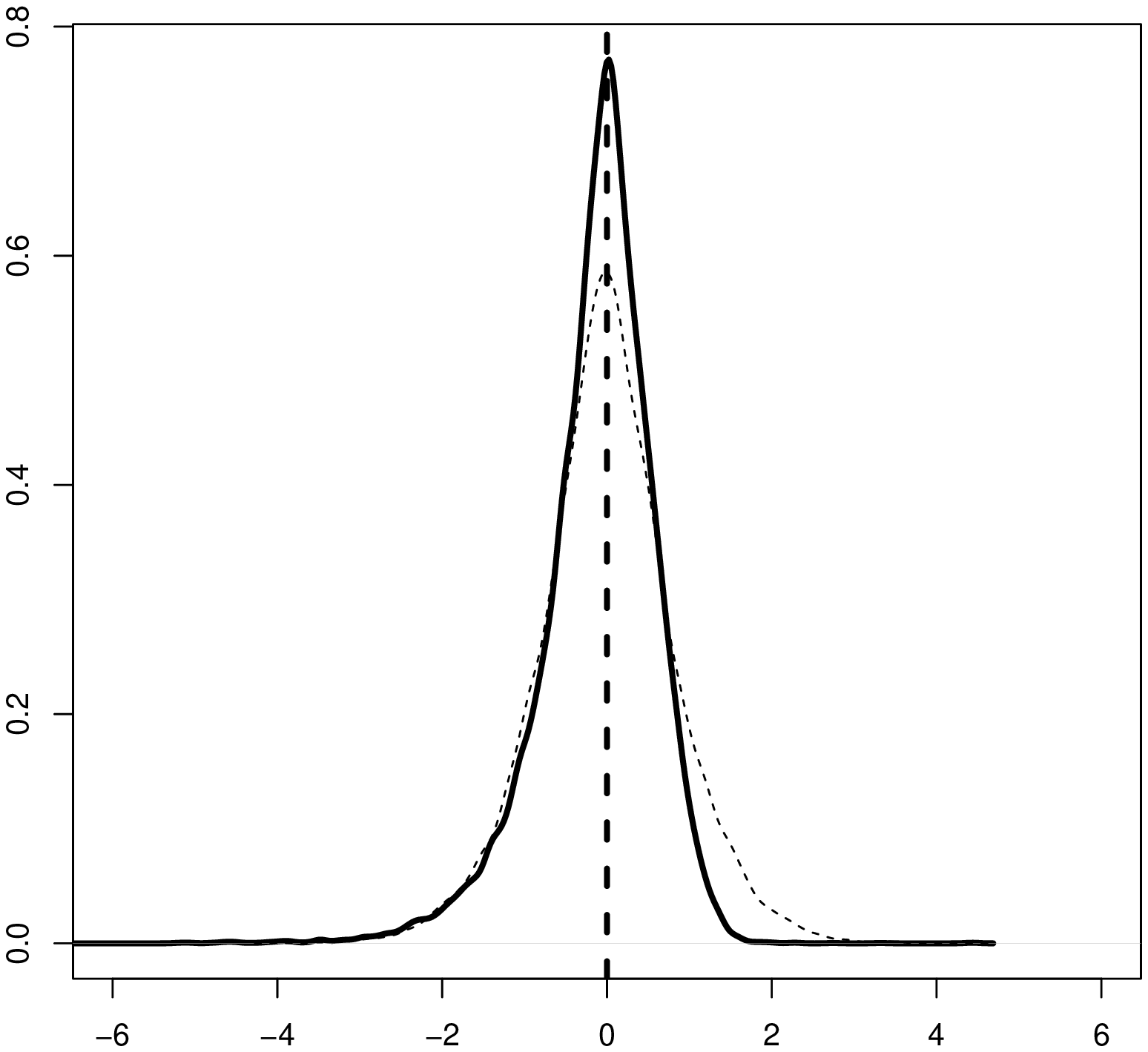}
\end{subfigure}\vspace{1.5cm}
\caption{{\small Distribution of $\sqrt{T}(\hat{\protect\beta}_2(\hat{\tau}_T)-\beta_2)$ when $\tau_0=0.5$, $\beta_2=0.8$, $\beta_1=1-1/T$. Left: $\{\varepsilon_t\}_{t=1}^T\sim t(3)$; Right: $\{\varepsilon_t\}_{t=1}^T\sim t(2)$.}}
\end{figure}

\section{Proof of Theorem \protect\ref{thm}}\label{sec3}

\setcounter{equation}{0} To prove Theorem \ref{thm} and Theorem \ref{thm2}
when the $\varepsilon _{t}$'s are heavy-tailed, we employ the truncation
technique in this paper. We let
\begin{equation*}
l(t)=E\varepsilon _{1}^{2}I\{|\varepsilon _{1}|\leq t\},~~~b=\inf \{t\geq
1:l(t)>0\},
\end{equation*}%
and
\begin{equation*}
\eta _{j}=\inf \{s:s\geq b+1,\frac{l(s)}{s^{2}}\leq \frac{1}{j}\},~~~\mathrm{%
for}~j=1,2,3,\cdots .
\end{equation*}%
Note that $Tl(\eta _{T})\leq \eta _{T}^{2}$ for all $T\geq 1$ and $\eta
_{T}^{2}\approx Tl(\eta _{T})$ for large $T$. In addition, for each $T$ we
let
\begin{equation*}
\left\{
\begin{array}{ll}
\varepsilon _{t}^{(1)}=\varepsilon _{t}I\{|\varepsilon _{t}|\leq \eta
_{T}\}-E\varepsilon _{t}I\{|\varepsilon _{t}|\leq \eta _{T}\}, &  \\
\varepsilon _{t}^{(2)}=\varepsilon _{t}I\{|\varepsilon _{t}|>\eta
_{T}\}-E\varepsilon _{t}I\{|\varepsilon _{t}|>\eta _{T}\}. &
\end{array}%
\right.
\end{equation*}%
for $t=1,\cdots ,T$.\newline

The following two lemmas are taken from Cs\"{o}rg\H{o} et al. (2003) and
Pang and Zhang (2011) respectively:

\begin{lemma}\label{lem1}
 Let $X$ be a random variable, and denote $l(x)=EX^{2}I\{|X|\leq
x\}$. The following statements are equivalent:\newline
\textrm{(1a)}~ $X$ is in the domain of attraction of the normal law,\newline
\textrm{(1b)}~ $x^{2}\mathbb{P}(|X|>x)=o(l(x))$,\newline
\textrm{(1c)}~ $xE(|X|I\{|X|>x\})=o(l(x))$,\newline
\textrm{(1d)}~ $E(|X|^{n}I\{|X|\leq x\})=o(x^{n-2}l(x))$ for $n>2$.
\end{lemma}

\begin{lemma}
\label{lem7} Suppose assumptions C1-C3 are satisfied, then in Model (\ref{model1}) with $|\beta _{1}|<1$, the following results hold jointly:\newline
\textrm{(2a)}~$\frac{1}{\sqrt{T}l(\eta _{T})}\sum_{t=1}^{[\tau
_{0}T]}y_{t-1}\varepsilon _{t}\Rightarrow N(0,\tau _{0}/(1-\beta _{1}^{2}))$,%
\newline
\textrm{(2b)}~$\frac{1}{Tl(\eta _{T})}\sum_{t=1}^{[\tau _{0}T]}y_{t-1}^{2}%
\overset{p}{\rightarrow }\tau _{0}/(1-\beta _{1}^{2})$.
\end{lemma}

The following two lemmas are useful in proving Theorem \ref{thm}:

\begin{lemma}
\label{lem6} Suppose assumption C1 is satisfied and $\beta_{2}=\beta_{2T}=1-c/T$, where $c$ is a fixed constant, then for any $\tau_{0}\leq
\tau \leq 1$,
\begin{eqnarray*}
X_{T}(\tau_{0},\tau ) &:=&\frac{1}{\sqrt{Tl(\eta_{T})}}\sum_{t=[\tau
_{0}T]+1}^{[\tau T]}\beta _{2T}^{T-t}\varepsilon _{t} \\
&\Rightarrow &e^{-c(1-\tau )}(W(\tau )-W(\tau_{0}))-c\int_{\tau_{0}}^{\tau
}e^{-c(1-t)}(W(t)-W(\tau _{0}))dt=F(W,c,\tau_{0},\tau )
\end{eqnarray*}%
and for any $0\leq s\leq 1-\tau _{0}$
\begin{eqnarray*}
Z_{T}(\tau_{0},s) &:=&\frac{1}{\sqrt{Tl(\eta_{T})}}\sum_{t=[\tau_{0}T]+1}^{[\tau_{0}T]+[sT]}\beta_{2T}^{[\tau_{0}T]+[sT]-t}\varepsilon
_{t} \\
&\Rightarrow &W(\tau_{0}+s)-W(\tau _{0})-c\int_{\tau_{0}}^{\tau_{0}+s}e^{-c(\tau_{0}+s-t)}(W(t)-W(\tau_{0}))dt=I(W,c,\tau_{0},s).
\end{eqnarray*}
\end{lemma}

\noindent \textbf{Proof.} \ Denote
\begin{equation*}
S_{T}(\tau _{0},\tau )=\frac{1}{\sqrt{Tl(\eta _{T})}}\sum_{t=[\tau
_{0}T]+1}^{[\tau T]}\varepsilon _{t}.
\end{equation*}%
The limiting distribution of $S_{T}(\tau _{0},\tau )$ is given by $W(\tau
)-W(\tau _{0})$, (See Theorem 1 in Cs\"{o}rg\H{o} et al. (2003)). In
addition, we define $c_{T}^{\ast }$ by
\begin{equation*}
\beta _{2T}=1-c/T=\exp (-c_{T}^{\ast }/T).
\end{equation*}%
Obviously, $c_{T}^{\ast }/c\rightarrow 1$, and we have
\begin{eqnarray*}
&&\frac{1}{\sqrt{Tl(\eta _{T})}}\sum_{t=[\tau _{0}T]+1}^{[\tau T]}\beta
_{2T}^{T-t}\varepsilon _{t} \\
&=&\sum_{t=[\tau _{0}T]+1}^{[\tau T]}e^{-(T-t)c_{T}^{\ast }/T}\frac{%
\varepsilon _{t}}{\sqrt{Tl(\eta _{T})}} \\
&=&\sum_{t=[\tau _{0}T]+1}^{[\tau T]}e^{-(T-t)c_{T}^{\ast }/T}\int_{\frac{t-1%
}{T}}^{\frac{t}{T}}dS_{T}(s,\tau _{0}) \\
&=&\Big(\sum_{t=[\tau _{0}T]+1}^{[\tau T]}\int_{\frac{t-1}{T}}^{\frac{t}{T}%
}e^{-(1-s)c_{T}^{\ast }}dS_{T}(s,\tau _{0})\Big)\cdot (1+o_{p}(1)) \\
&=&\Big(\int_{\tau _{0}}^{\tau }e^{-(1-s)c_{T}^{\ast }}dS_{T}(s,\tau _{0})%
\Big)\cdot (1+o_{p}(1)) \\
&=&\Big(e^{-(1-\tau )c_{T}^{\ast }}S_{T}(\tau ,\tau _{0})-c_{T}^{\ast
}\int_{\tau _{0}}^{\tau }e^{-(1-s)c_{T}^{\ast }}S_{T}(s,\tau _{0})ds\Big)%
\cdot (1+o_{p}(1)) \\
&\Rightarrow &e^{-c(1-\tau )}(W(\tau )-W(\tau _{0}))-c\int_{\tau _{0}}^{\tau
}e^{-c(1-t)}(W(t)-W(\tau _{0}))dt,
\end{eqnarray*}%
which finishes the proof of the first result.

For the second result, note that
\begin{eqnarray*}
Z_{T}(\tau _{0},s) &:=&\frac{1}{\sqrt{Tl(\eta _{T})}}\sum_{t=[\tau
_{0}T]+1}^{[\tau _{0}T]+[sT]}\beta _{2T}^{[\tau _{0}T]+[sT]-t}\varepsilon
_{t} \\
&=&\frac{1+o_{p}(1)}{\sqrt{Tl(\eta _{T})}}\sum_{t=[\tau _{0}T]+1}^{[(\tau
_{0}+s)T]}\beta _{2T}^{T-t}\varepsilon _{t}\cdot \beta _{2T}^{[\tau
_{0}T]+[sT]-T},
\end{eqnarray*}%
applying the first result of this lemma, we have
\begin{eqnarray*}
Z_{T}(\tau _{0},s) &\Rightarrow &F(W,c,\tau _{0},\tau _{0}+s)e^{c(1-\tau
_{0}-s)} \\
&=&W(\tau _{0}+s)-W(\tau _{0})-c\int_{\tau _{0}}^{\tau _{0}+s}e^{-c(\tau
_{0}+s-t)}(W(t)-W(\tau _{0}))dt.
\end{eqnarray*}%
The proofs are complete. $\hfill \Box $\newline

\begin{lemma}
\label{lem2} Let $\{y_{t},t\geq 1\}$ be generated according to Model
(\ref{model1}) with $|\beta _{1}|<1$ and $\beta _{2}=\beta _{2T}=1-c/T$, where
$c$ is a fixed constant. Under the assumptions C1-C3, the following
results hold jointly:\newline
\textrm{(3a)}~$\frac{1}{Tl(\eta _{T})}\sum_{t=[\tau
_{0}T]+1}^{T}y_{t-1}\varepsilon _{t}\Rightarrow \frac{1}{2}F^{2}(W,c,\tau
_{0},1)+c\int_{\tau _{0}}^{1}e^{2c(1-t)}F^{2}(W,c,\tau _{0},t)dt-\frac{1}{2}%
(1-\tau _{0}),$\newline
\textrm{(3b)}~$\frac{1}{T^{2}l(\eta _{T})}\sum_{t=[\tau
_{0}T]+1}^{T}y_{t-1}^{2}\Rightarrow \int_{\tau
_{0}}^{1}e^{2c(1-t)}F^{2}(W,c,\tau _{0},t)dt.$
\end{lemma}

\noindent \textbf{Proof.} \ To prove (3a), we first square the equation $%
y_{t}=(1-c/T)y_{t-1}+\varepsilon _{t}$ and apply the summation from $[\tau
_{0}T]+1$ to $T$ to obtain
\begin{equation*}
2(1-\frac{c}{T})\sum_{t=[\tau _{0}T]+1}^{T}y_{t-1}\varepsilon
_{t}=y_{T}^{2}-y_{[\tau _{0}T]}^{2}+(\frac{2c}{T}-\frac{c^{2}}{T^{2}}%
)\sum_{t=[\tau _{0}T]+1}^{T}y_{t-1}^{2}-\sum_{t=[\tau
_{0}T]+1}^{T}\varepsilon _{t}^{2},
\end{equation*}%
which can be rewritten as
\begin{eqnarray}
\frac{\sum_{t=[\tau _{0}T]+1}^{T}y_{t-1}\varepsilon _{t}}{Tl(\eta _{T})} &=&%
\frac{T}{2(T-c)}\frac{y_{T}^{2}-y_{[\tau _{0}T]}^{2}}{Tl(\eta _{T})}+\frac{%
c(2T-c)}{2T(T-c)}\sum_{t=[\tau _{0}T]+1}^{T}\Big(\frac{y_{t-1}}{\sqrt{%
Tl(\eta _{T})}}\Big)^{2}  \notag  \label{2.3} \\
&&-\frac{T}{2(T-c)}\frac{1}{Tl(\eta _{T})}\sum_{t=[\tau
_{0}T]+1}^{T}\varepsilon _{t}^{2}.
\end{eqnarray}

Consider the first term in the right hand side of (\ref{2.3}). Write
\begin{eqnarray}
\frac{y_{T}^{2}-y_{[\tau _{0}T]}^{2}}{Tl(\eta _{T})} &=&\frac{(\sum_{t=[\tau
_{0}T]+1}^{T}\beta _{2T}^{T-t}\varepsilon _{t}+\beta _{2T}^{T-[\tau
_{0}T]}y_{[\tau _{0}T]})^{2}-y_{[\tau _{0}T]}^{2}}{Tl(\eta _{T})}  \notag
\label{2.1} \\
&=&\Big(X_{T}(\tau _{0},1)+\frac{\beta _{2T}^{T-[\tau _{0}T]}y_{[\tau _{0}T]}%
}{\sqrt{Tl(\eta _{T})}}\Big)^{2}-\frac{y_{[\tau _{0}T]}^{2}}{Tl(\eta _{T})}.
\end{eqnarray}%
Since $y_{[\tau _{0}T]}$ is generated from a stationary AR(1) model with
initial condition $y_{0}=o_{p}(\sqrt{T})$, it is obvious that
\begin{equation*}
\frac{|y_{[\tau _{0}T]}|}{\sqrt{Tl(\eta _{T})}}\overset{p}{\rightarrow } 0,
\end{equation*}%
which together with Lemma \ref{lem6} and (\ref{2.1}) lead to
\begin{equation*}
\frac{y_{T}^{2}-y_{[\tau _{0}T]}^{2}}{Tl(\eta _{T})}\Rightarrow
~F^{2}(W,c,\tau _{0},1)
\end{equation*}%
by noting that $|\beta _{2T}^{T-[\tau _{0}T]}|<C$. Or equivalently
\begin{equation}  \label{2.13}
\frac{T}{2(T-c)}\frac{y_{T}^{2}-y_{[\tau _{0}T]}^{2}}{Tl(\eta _{T})}%
\Rightarrow ~\frac{1}{2}F^{2}(W,c,\tau _{0},1).
\end{equation}

Consider the second term in the right hand side of (\ref{2.3}). Write
\begin{eqnarray}
\frac{1}{T}\sum_{t=[\tau _{0}T]+1}^{T}\Big(\frac{y_{t-1}}{\sqrt{Tl(\eta _{T})%
}}\Big)^{2} &=&\frac{1}{T}\sum_{t=[\tau _{0}T]+1}^{T}\Big(X_{T}(\tau _{0},%
\frac{t-1}{T})\beta _{2T}^{t-1-T}+\frac{\beta _{2T}^{t-[\tau
_{0}T]-1}y_{[\tau _{0}T]}}{\sqrt{Tl(\eta _{T})}}\Big)^{2}  \notag
\label{2.12} \\
 &:=&\frac{1}{T}\sum_{t=[\tau _{0}T]+1}^{T}(III+IV)^{2}.
\end{eqnarray}%
Note that
\begin{eqnarray}
\frac{1}{T}\sum_{t=[\tau _{0}T]+1}^{T}(IV)^{2} &=&\frac{y_{[\tau _{0}T]}^{2}%
}{T^{2}l(\eta _{T})}\sum_{t=[\tau _{0}T]+1}^{T}\beta _{2T}^{2(t-[\tau
_{0}T]-1)}  \notag  \label{2.15} \\
&=&\frac{y_{[\tau _{0}T]}^{2}}{Tl(\eta _{T})}\cdot \frac{(1-\beta
_{2T}^{2(T-[\tau _{0}T])})}{T(1-\beta _{2T}^{2})}\overset{p}{\rightarrow }0.
\end{eqnarray}

Consider the term $\frac{1}{T}\sum_{t=[\tau _{0}T]+1}^{T}(III)^{2}$. Write
\begin{eqnarray}
&&\frac{1}{T}\sum_{t=[\tau _{0}T]+1}^{T}(III)^{2}  \notag  \label{2.5} \\
&=&\frac{1}{T}\sum_{t=[\tau _{0}T]+1}^{T}X_{T}^{2}(\tau _{0},\frac{t-1}{T}%
)\beta _{2T}^{2(t-1-T)}-\int_{\tau _{0}}^{1}e^{2c(1-t)}X_{T}^{2}(\tau
_{0},t)dt+\int_{\tau _{0}}^{1}e^{2c(1-t)}X_{T}^{2}(\tau _{0},t)dt  \notag \\
 &:=&R_{T}(\tau _{0})+\int_{\tau _{0}}^{1}e^{2c(1-t)}X_{T}^{2}(\tau
_{0},t)dt.
\end{eqnarray}%
First, we shall employ the similar arguments used in Lemma 2.2 in Chan and
Wei (1987) to prove that
\begin{equation}
R_{T}(\tau _{0})\overset{p}{\rightarrow } 0.  \label{2.4}
\end{equation}%
To this end, write
\begin{eqnarray*}
&&|R_{T}(\tau _{0})| \\
&=&\Big|\sum_{t=[\tau _{0}T]}^{T-1}\int_{t/T}^{(t+1)/T}\beta
_{2T}^{2(t-T)}X_{T}^{2}(\tau _{0},\frac{t}{T})ds-\sum_{t=[\tau
_{0}T]}^{T-1}\int_{t/T}^{(t+1)/T}e^{2c(1-s)}X_{T}^{2}(\tau _{0},s)ds\Big|%
+o_{p}(1) \\
&=&\Big|\sum_{t=[\tau _{0}T]}^{T-1}\int_{t/T}^{(t+1)/T}\beta
_{2T}^{2(t-T)}X_{T}^{2}(\tau _{0},s)ds-\sum_{t=[\tau
_{0}T]}^{T-1}\int_{t/T}^{(t+1)/T}e^{2c(1-s)}X_{T}^{2}(\tau _{0},s)ds\Big|%
+o_{p}(1) \\
&\leq &\sum_{t=[\tau _{0}T]}^{T-1}\int_{t/T}^{(t+1)/T}\Big|\beta
_{2T}^{2(t-T)}-e^{2c(1-s)}\Big|X_{T}^{2}(\tau _{0},s)ds+o_{p}(1) \\
&\leq &\max_{[\tau _{0}T]\leq t\leq T-1}\sup_{t/T\leq s\leq (t+1)/T}|\beta
_{2}^{2(t-T)}-e^{2c(1-s)}|\cdot \sum_{t=[\tau
_{0}T]}^{T-1}\int_{t/T}^{(t+1)/T}X_{T}^{2}(\tau _{0},s)ds+o_{p}(1).
\end{eqnarray*}%
Noting that
\begin{equation*}
\max_{\lbrack \tau _{0}T]\leq t\leq T-1}\sup_{t/T\leq s\leq (t+1)/T}|\beta
_{2}^{2(t-T)}-e^{2c(1-s)}|\rightarrow 0
\end{equation*}%
by (2.7) in Chan and Wei (1987), and
\begin{equation*}
\sum_{t=[\tau _{0}T]}^{T-1}\int_{t/T}^{(t+1)/T}X_{T}^{2}(\tau
_{0},s)ds=\int_{\tau _{0}}^{1}X_{T}^{2}(\tau _{0},s)ds+o_{p}(1)\Rightarrow
\int_{\tau _{0}}^{1}F^{2}(W,c,\tau _{0},t)dt
\end{equation*}%
by Lemma \ref{lem6}, hence (\ref{2.4}) is proved. Combining (\ref{2.5}) and (%
\ref{2.4}), we have
\begin{equation*}
\frac{1}{T}\sum_{t=[\tau _{0}T]+1}^{T}(III)^{2}\Rightarrow ~\int_{\tau
_{0}}^{1}e^{2c(1-t)}F^{2}(W,c,\tau _{0},t)dt,
\end{equation*}%
which together with (\ref{2.12}), (\ref{2.15}) and Cauchy-Schwarz inequality
immediately yield
\begin{equation}
\frac{1}{T}\sum_{t=[\tau _{0}T]+1}^{T}\Big(\frac{y_{t-1}}{\sqrt{Tl(\eta _{T})%
}}\Big)^{2}\Rightarrow ~\int_{\tau _{0}}^{1}e^{2c(1-t)}F^{2}(W,c,\tau
_{0},t)dt.  \label{2.8}
\end{equation}%
This proves (3b) of Lemma 3.4. Note that the above can be written as
\begin{equation}
\frac{c(2T-c)}{2T(T-c)}\sum_{t=[\tau _{0}T]+1}^{T}\Big(\frac{y_{t-1}}{\sqrt{%
Tl(\eta _{T})}}\Big)^{2}\Rightarrow ~c\int_{\tau
_{0}}^{1}e^{2c(1-t)}F^{2}(W,c,\tau _{0},t)dt.  \label{2.20}
\end{equation}

For the third term in the right hand side of (\ref{2.3}), it is obvious that
\begin{equation}  \label{2.6}
\frac{T}{2(T-c)}\frac{1}{Tl(\eta _{T})}\sum_{t=[\tau
_{0}T]+1}^{T}\varepsilon _{t}^{2}\overset{p}{\rightarrow } \frac{1-\tau _{0}%
}{2}
\end{equation}
when $\varepsilon _{t}^{\prime }s$ are i.i.d. and from DAN.

Combining (\ref{2.13}), (\ref{2.20}) and (\ref{2.6}) together gives (3a) of
Lemma 3.4.$\hfill \Box $\newline

\noindent\textbf{Proof of Theorem \ref{thm}}. Consider the first part of (%
\ref{result1}). Along the lines of the proof of Theorem 3 in Chong (2001),
it is sufficient to prove that
\begin{eqnarray}  \label{2.30}
\left\{%
\begin{array}{ll}
\displaystyle A_1=\frac{\sum_{t=[\tau_0 T]+1}^T y_{t-1}\varepsilon_t}{%
\sum_{t=[\tau_0 T]+1}^T y^2_{t-1}}=o_p(1), &  \\
\displaystyle A_2=\sup_{m\in D_{1T}}\frac{\sum_{t=m+1}^{[\tau_0
T]}y_{t-1}\varepsilon_t}{\sum_{t=m+1}^{[\tau_0 T]}y^2_{t-1}}=o_p(1), &  \\
\displaystyle A_3=\sup_{m\in D_{1T}}\Big|\frac{\sum_{t=m+1}^{T}y^2_{t-1}}{%
\sum_{t=[\tau_0 T]+1}^T y^2_{t-1}\sum_{t=m+1}^{[\tau_0 T]}y^2_{t-1}}%
\Lambda_T(\frac{m}{T}) \Big|=o_p(1), &  \\
\displaystyle A_4=\frac{\sum_{t=1}^{[\tau_0 T]} y_{t-1}\varepsilon_t}{%
\sum_{t=1}^{[\tau_0 T]} y^2_{t-1}}=o_p(1), &  \\
\displaystyle A_5=\sup_{m\in D_{2T}}\frac{\sum_{t=[\tau_0
T]+1}^{m}y_{t-1}\varepsilon_t}{\sum_{t=[\tau_0 T]+1}^{m}y^2_{t-1}}=o_p(1), &
\\
\displaystyle A_6=\sup_{m\in D_{2T}}\Big|\frac{\sum_{t=1}^{m}y^2_{t-1}}{%
\sum_{t=[\tau_0 T]+1}^m y^2_{t-1}\sum_{t=1}^{[\tau_0 T]}y^2_{t-1}}\Lambda_T(%
\frac{m}{T}) \Big|=o_p(1) &
\end{array}%
\right.
\end{eqnarray}
in the case of $|\beta_1|<1$ and $\beta_2=\beta_{2T}=1-c/T$, where $c$ is a
fixed constant,
\begin{eqnarray}  \label{2.25}
\Lambda_T(\frac{m}{T})&=&\frac{(\sum_{t=1}^{[\tau_0 T]}
y_{t-1}\varepsilon_t)^2}{\sum_{t=1}^{[\tau_0 T]} y^2_{t-1}}\bigg(1-\Big(%
\frac{\sum_{t=1}^{m} y_{t-1}\varepsilon_t}{\sum_{t=1}^{[\tau_0 T]}
y_{t-1}\varepsilon_t}\Big)^2\frac{\sum_{t=1}^{[\tau_0 T]} y^2_{t-1}}{%
\sum_{t=1}^{m} y^2_{t-1}}\bigg)  \notag \\
&&+\frac{(\sum_{t=[\tau_0 T]+1}^T y_{t-1}\varepsilon_t)^2}{\sum_{t=[\tau_0
T]+1}^T y^2_{t-1}}\bigg(1-\Big(\frac{\sum_{t=m+1}^{T} y_{t-1}\varepsilon_t}{%
\sum_{t=[\tau_0 T]+1}^{T} y_{t-1}\varepsilon_t}\Big)^2\frac{\sum_{t=[\tau_0
T]+1}^{T} y^2_{t-1}}{\sum_{t=m+1}^{T} y^2_{t-1}}\bigg)
\end{eqnarray}
and
\begin{eqnarray*}
\left\{%
\begin{array}{ll}
D_{1T}=\{m: m\in Z_T, m< [\tau_0 T]-M_T\}, &  \\
D_{2T}=\{m: m\in Z_T, m> [\tau_0 T]+M_T\} &
\end{array}%
\right.
\end{eqnarray*}
with $M_T>0$ such that $M_T\rightarrow \infty$ and $M_T/T\rightarrow 0$,
where $Z_T$ denotes the set $\{0, 1, 2, \cdots, T\}$.

We apply Lemmas \ref{lem7} and \ref{lem2} to prove (\ref{2.30}). The results of $A_1$ and $A_4$ are obvious. For
$A_2$, Applying the uniform law of large numbers in Andrews (1987,
Theorem 1) yields
\begin{equation*}
|A_{2}|=\Big|\sup_{m\in D_{1T}}\frac{\sum_{t=m+1}^{[\tau
_{0}T]}y_{t-1}\varepsilon _{t}}{\sum_{t=m+1}^{[\tau _{0}T]}y_{t-1}^{2}}\Big|%
\leq O_{p}\Big(\frac{1}{\sqrt{M_{T}}}\big)=o_{p}(1).
\end{equation*}%
For $A_3, A_5$ and $A_6$, note that
\begin{eqnarray*}
A_{3} &=&\sup_{m\in D_{1T}}\Big|\frac{\sum_{t=m+1}^{T}y_{t-1}^{2}}{%
\sum_{t=[\tau _{0}T]+1}^{T}y_{t-1}^{2}\sum_{t=m+1}^{[\tau _{0}T]}y_{t-1}^{2}}%
\Lambda _{T}(\frac{m}{T})\Big| \\
&\leq &\bigg(\frac{1}{\sum_{t=[\tau _{0}T]+1}^{T}y_{t-1}^{2}}+\frac{1}{%
\sum_{t=[\tau _{0}T]-M_{T}}^{[\tau _{0}T]}y_{t-1}^{2}}\bigg)\bigg(\sup_{m\in
D_{1T}}\Big|\frac{(\sum_{t=1}^{[\tau _{0}T]}y_{t-1}\varepsilon _{t})^{2}}{%
\sum_{t=1}^{[\tau _{0}T]}y_{t-1}^{2}}-\frac{(\sum_{t=1}^{m}y_{t-1}%
\varepsilon _{t})^{2}}{\sum_{t=1}^{m}y_{t-1}^{2}}\Big| \\
&&~~~~+\sup_{m\in D_{1T}}\Big|\frac{(\sum_{t=[\tau
_{0}T]+1}^{T}y_{t-1}\varepsilon _{t})^{2}}{\sum_{t=[\tau
_{0}T]+1}^{T}y_{t-1}^{2}}-\frac{(\sum_{t=m+1}^{T}y_{t-1}\varepsilon _{t})^{2}%
}{\sum_{t=m+1}^{T}y_{t-1}^{2}}\Big|\bigg) \\
&=&\Big(O_{p}(\frac{1}{T^{2}l(\eta _{T})})+O_{p}(\frac{1}{M_{T}l(\eta _{T})})%
\Big)\cdot O_{p}(l(\eta _{T})) \\
&=&O_{p}(\frac{1}{M_{T}})=o_{p}(1).
\end{eqnarray*}%
For the fifth part of (\ref{2.30}), we have
\begin{eqnarray*}
A_{5} &=&\sup_{m\in D_{2T}}\frac{\sum_{t=[\tau
_{0}T]+1}^{m}y_{t-1}\varepsilon _{t}}{\sum_{t=[\tau _{0}T]+1}^{m}y_{t-1}^{2}}
\\
&=&\sup_{m\in D_{2T}}O_{p}\Big(\frac{1}{m-[\tau _{0}T]}\Big) \\
&=&O_{p}(\frac{1}{M_{T}})=o_{p}(1).
\end{eqnarray*}%
We now show the last part of (\ref{2.30}).
\begin{eqnarray*}
A_{6} &=&\sup_{m\in D_{2T}}\Big|\frac{\sum_{t=1}^{m}y_{t-1}^{2}}{%
\sum_{t=[\tau _{0}T]+1}^{m}y_{t-1}^{2}\sum_{t=1}^{[\tau _{0}T]}y_{t-1}^{2}}%
\Lambda _{T}(\frac{m}{T})\Big| \\
&\leq &\bigg(\frac{1}{\sum_{t=1}^{[\tau _{0}T]}y_{t-1}^{2}}+\frac{1}{%
\sum_{t=[\tau _{0}T]+1}^{[\tau _{0}T]+M_{T}}y_{t-1}^{2}}\bigg)\bigg(%
\sup_{m\in D_{2T}}\Big|\frac{(\sum_{t=1}^{[\tau _{0}T]}y_{t-1}\varepsilon
_{t})^{2}}{\sum_{t=1}^{[\tau _{0}T]}y_{t-1}^{2}}-\frac{%
(\sum_{t=1}^{m}y_{t-1}\varepsilon _{t})^{2}}{\sum_{t=1}^{m}y_{t-1}^{2}}\Big|
\\
&&~~~~+\sup_{m\in D_{2T}}\Big|\frac{(\sum_{t=[\tau
_{0}T]+1}^{T}y_{t-1}\varepsilon _{t})^{2}}{\sum_{t=[\tau
_{0}T]+1}^{T}y_{t-1}^{2}}-\frac{(\sum_{t=m+1}^{T}y_{t-1}\varepsilon _{t})^{2}%
}{\sum_{t=m+1}^{T}y_{t-1}^{2}}\Big|\bigg) \\
&=&\Big(O_{p}(\frac{1}{Tl(\eta _{T})})+O_{p}(\frac{1}{M_{T}^{2}l(\eta _{T})})%
\Big)\cdot O_{p}(l(\eta _{T})) \\
&=&o_{p}(\frac{1}{M_{T}})=o_{p}(1).
\end{eqnarray*}%
Hence, the first part of (\ref{result1}) is proved.

To find the limiting distribution of ${\hat{\beta}}_{1}({\hat{\tau}}_{T})$,
first note that $\hat{\tau}_{T}-\tau _{0}=O_{p}(1/T)$. Following Appendix G
in Chong (2001), we have
\begin{eqnarray*}
&&\sqrt{T}(\hat{\beta}_{1}(\hat{\tau}_{T})-\hat{\beta}_{1}(\tau _{0})) \\
&=&\sqrt{T}\bigg(\frac{\sum_{t=1}^{[\hat{\tau}_{T}T]}y_{t}y_{t-1}}{%
\sum_{t=1}^{[\hat{\tau}_{T}T]}y_{t}^{2}}-\frac{\sum_{t=1}^{[\tau
_{0}T]}y_{t}y_{t-1}}{\sum_{t=1}^{[\tau _{0}T]}y_{t}^{2}}\bigg) \\
&=&I\{\hat{\tau}_{T}\leq \tau _{0}\}\sqrt{T}\bigg(\frac{\sum_{t=[\hat{\tau}%
_{T}T]+1}^{[\tau _{0}T]}y_{t-1}^{2}}{\sum_{t=1}^{[\hat{\tau}%
_{T}T]}y_{t-1}^{2}}\frac{\sum_{t=1}^{[\tau _{0}T]}y_{t-1}\varepsilon _{t}}{%
\sum_{t=1}^{[\tau _{0}T]}y_{t-1}^{2}}-\frac{\sum_{t=[\hat{\tau}%
_{T}T]+1}^{[\tau _{0}T]}y_{t-1}\varepsilon _{t}}{\sum_{t=1}^{[\hat{\tau}%
_{T}T]}y_{t-1}^{2}}\bigg) \\
&&+I\{\hat{\tau}_{T}>\tau _{0}\}\sqrt{T}\bigg(-\frac{\sum_{t=[\tau
_{0}T]+1}^{[\hat{\tau}_{T}T]}y_{t-1}^{2}}{\sum_{t=1}^{[\hat{\tau}%
_{T}T]}y_{t-1}^{2}}\frac{\sum_{t=1}^{[\tau _{0}T]}y_{t-1}\varepsilon _{t}}{%
\sum_{t=1}^{[\tau _{0}T]}y_{t-1}^{2}}+\frac{\sum_{t=[\tau _{0}T]+1}^{[\hat{%
\tau}_{T}T]}y_{t-1}\varepsilon _{t}}{\sum_{t=1}^{[\hat{\tau}%
_{T}T]}y_{t-1}^{2}} \\
&&~~~~~~~~~~~~~~~~~~~~~~~~~~~+(\beta _{2T}-\beta _{1})\frac{\sum_{t=[\tau
_{0}T]+1}^{[\hat{\tau}_{T}T]}y_{t-1}^{2}}{\sum_{t=1}^{[\hat{\tau}%
_{T}T]}y_{t-1}^{2}}\bigg) \\
&=&I\{\hat{\tau}_{T}\leq \tau _{0}\}\sqrt{T}\Big(O_{p}(\frac{l(\eta _{T})}{%
Tl(\eta _{T})})O_{p}(\frac{1}{\sqrt{T}})-O_{p}(\frac{l(\eta _{T})}{Tl(\eta
_{T})})\Big) \\
&&+I\{\hat{\tau}_{T}>\tau _{0}\}\sqrt{T}\Big(-O_{p}(\frac{l(\eta _{T})}{%
Tl(\eta _{T})})O_{p}(\frac{1}{\sqrt{T}})+O_{p}(\frac{l(\eta _{T})}{Tl(\eta
_{T})})+O_{p}(\frac{l(\eta _{T})}{Tl(\eta _{T})})\Big) \\
&=&o_{p}(1),
\end{eqnarray*}%
one is referred to Chong (2001) for more details. Thus, $\hat{\beta}_{1}(%
\hat{\tau}_{T})$ and $\hat{\beta}_{1}(\tau _{0})$ have the same asymptotic
distribution. Applying Lemma \ref{lem7}, we have
\begin{equation*}
\sqrt{T}(\hat{\beta}_{1}(\hat{\tau}_{T})-\beta _{1})\overset{d}{=} \sqrt{T}(%
\hat{\beta}_{1}(\tau _{0})-\beta _{1})=\frac{\frac{1}{\sqrt{T}l(\eta _{T})}%
\sum_{t=1}^{[\tau _{0}T]}y_{t-1}\varepsilon _{t}}{\frac{1}{Tl(\eta _{T})}%
\sum_{t=1}^{[\tau _{0}T]}y_{t-1}^{2}}\Rightarrow ~N\big(0,\frac{1-\beta
_{1}^{2}}{\tau _{0}}\big).
\end{equation*}

Similarly, we have
\begin{eqnarray*}
&&T(\hat{\beta}_{2}(\hat{\tau}_{T})-\hat{\beta _{2}}(\tau _{0})) \\
&=&T\bigg(\frac{\sum_{t=[\hat{\tau}_{T}T]+1}^{T}y_{t}y_{t-1}}{\sum_{t=[\hat{%
\tau}_{T}T]+1}^{T}y_{t}^{2}}-\frac{\sum_{t=[\tau _{0}T]+1}^{T}y_{t}y_{t-1}}{%
\sum_{t=[\tau _{0}T]+1}^{T}y_{t}^{2}}\bigg) \\
&=&I\{\hat{\tau}_{T}\leq \tau _{0}\}T\bigg(-\frac{\sum_{t=[\hat{\tau}%
_{T}T]+1}^{[\tau _{0}T]}y_{t-1}^{2}}{\sum_{t=[\hat{\tau}%
_{T}T]+1}^{T}y_{t-1}^{2}}\frac{\sum_{t=[\tau _{0}T]+1}^{T}y_{t-1}\varepsilon
_{t}}{\sum_{t=[\tau _{0}T]+1}^{T}y_{t-1}^{2}}+\frac{\sum_{t=[\hat{\tau}%
_{T}T]+1}^{[\tau _{0}T]}y_{t-1}\varepsilon _{t}}{\sum_{t=[\hat{\tau}%
_{T}T]+1}^{T}y_{t-1}^{2}} \\
&&~~~~~~~~~~~~~~~~~~~~~~~~+(\beta _{1}-\beta _{2T})\frac{\sum_{t=[\hat{\tau}%
_{T}T]+1}^{[\tau _{0}T]}y_{t-1}^{2}}{\sum_{t=[\hat{\tau}%
_{T}T]+1}^{T}y_{t-1}^{2}}\bigg) \\
&&+I\{\hat{\tau}_{T}>\tau _{0}\}T\bigg(\frac{\sum_{t=[\tau _{0}T]+1}^{[\hat{%
\tau}_{T}T]}y_{t-1}^{2}}{\sum_{t=[\hat{\tau}_{T}T]+1}^{T}y_{t-1}^{2}}\frac{%
\sum_{t=[\tau _{0}T]+1}^{T}y_{t-1}\varepsilon _{t}}{\sum_{t=[\tau
_{0}T]+1}^{T}y_{t-1}^{2}}-\frac{\sum_{t=[\tau _{0}T]+1}^{[\hat{\tau}%
_{T}T]}y_{t-1}\varepsilon _{t}}{\sum_{t=[\hat{\tau}_{T}T]+1}^{T}y_{t-1}^{2}}%
\bigg) \\
&=&I\{\hat{\tau}_{T}\leq \tau _{0}\}T\Big(-O_{p}(\frac{l(\eta _{T})}{%
T^{2}l(\eta _{T})})O_{p}(\frac{1}{T})+O_{p}(\frac{l(\eta _{T})}{T^{2}l(\eta
_{T})})-O_{p}(\frac{l(\eta _{T})}{T^{2}l(\eta _{T})})\Big) \\
&&+I\{\hat{\tau}_{T}>\tau _{0}\}T\Big(O_{p}(\frac{l(\eta _{T})}{T^{2}l(\eta
_{T})})O_{p}(\frac{1}{T})-O_{p}(\frac{l(\eta _{T})}{T^{2}l(\eta _{T})})\Big)
\\
&=&o_{p}(1).
\end{eqnarray*}%
Thus, $\hat{\beta}_{2}(\hat{\tau}_{T})$ and $\hat{\beta}_{2}(\tau _{0})$
also have the same asymptotic distribution. Applying Lemma \ref{lem2}, we
have
\begin{eqnarray*}
T(\hat{\beta}_{2}(\hat{\tau}_{T})-\beta _{2}) &\overset{d}{=} &T(\hat{\beta}%
_{2}(\tau _{0})-\beta _{2}) \\
&=&\frac{\frac{1}{Tl(\eta _{T})}\sum_{t=[\tau
_{0}T]+1}^{T}y_{t-1}\varepsilon _{t}}{\frac{1}{T^{2}l(\eta _{T})}%
\sum_{t=[\tau _{0}T]+1}^{T}y_{t-1}^{2}} \\
&\Rightarrow &\frac{\frac{1}{2}F^{2}(W,c,\tau _{0},1)+c\int_{\tau
_{0}}^{1}e^{2c(1-t)}F^{2}(W,c,\tau _{0},t)dt-\frac{1}{2}(1-\tau _{0})}{%
\int_{\tau _{0}}^{1}e^{2c(1-t)}F^{2}(W,c,\tau _{0},t)dt}.
\end{eqnarray*}

To derive the limiting distribution of $\hat{\tau}_{T}$ for shrinking shift,
we let $\beta _{2}=\beta _{2T}=1-c/T$ and $\beta _{1}=\beta _{1T}=\beta
_{2T}-1/g(T)$ in the remaining proof of Theorem \ref{thm}, where $g(T)>0$
with $g(T)\rightarrow \infty $ and $g(T)/T\rightarrow 0$. Note that $\beta
_{1T}=1-1/g(T)+o(1/g(T))$. Hence, the sequence $\{y_t, 1\le t\le [\tau_0
T]\} $ is generated from a mildly integrated AR(1) model, as a result, the
results or ideas from Phillips and Magdalinos (2007) and Huang et al. (2012)
could be applied directly. Following Chong (2001), first, for $\tau =\tau
_{0}+\nu g(T)/T$ and $\nu \leq 0$, by recalling (\ref{2.25}), we have
\begin{eqnarray}
|\Lambda _{T}(\tau )| &=&O_{p}(l(\eta _{T}))\Big(%
1-(1-o_{p}(1))^{2}(1+o_{p}(1)\Big)+O_{p}(l(\eta _{T}))\Big(%
1-(1+o_{p}(1))^{2}(1-o_{p}(1))\Big)  \notag  \label{2.7} \\
&=&o_{p}(l(\eta _{T})).
\end{eqnarray}%
Second, for any $t=0,\cdots ,[|\nu |g(T)]-1,$ we have%
\begin{eqnarray}
\frac{y_{[\tau _{0}T]-t-1}}{\sqrt{g(T)l(\eta _{T})}} &=&\frac{1}{\sqrt{%
g(T)l(\eta _{T})}}\sum_{i=1}^{[\tau _{0}T]-t-1}\beta _{1T}^{[\tau
_{0}T]-t-1-i}\varepsilon _{i}^{(1)}  \notag  \label{2.14} \\
&&+\frac{1}{\sqrt{g(T)l(\eta _{T})}}\sum_{i=1}^{[\tau _{0}T]-t-1}\beta
_{1T}^{[\tau _{0}T]-t-1-i}\varepsilon _{i}^{(2)}+\frac{\beta _{1T}^{[\tau
_{0}T]-t-1}y_{0}}{\sqrt{g(T)l(\eta _{T})}}.
\end{eqnarray}%
It is not difficult to show that
\begin{eqnarray}  \label{3.7}
\frac{\beta_{1T}^{[\tau_0 T]-t-1}y_0}{\sqrt{g(T)l(\eta_T)}}=\sqrt{\frac{T}{%
g(T)}}\beta_{1T}^{[\tau_0 T]-t-1}\cdot\frac{y_0}{\sqrt{T l(\eta_T)}}\overset{%
p}{\rightarrow} 0
\end{eqnarray}
by $y_0=o_p(\sqrt{T})$ and $g(T)=o(T)$; see the proof of Proposition A.1 in
Phillips and Magdalinos (2007) for more details. In addition, it follows
from Lemma \ref{lem1} that
\begin{eqnarray*}
&&\frac{1}{\sqrt{g(T)l(\eta _{T})}}\cdot E\Big|\sum_{i=1}^{[\tau
_{0}T]-t-1}\beta _{1T}^{[\tau _{0}T]-t-1-i}\varepsilon _{i}^{(2)}\Big| \\
&=&\frac{1}{\sqrt{g(T)l(\eta _{T})}}\cdot \sum_{i=1}^{[\tau
_{0}T]-t-1}|\beta _{1T}|^{[\tau _{0}T]-t-1-i}\cdot o(\frac{l(\eta _{T})}{%
\eta _{T}})=o(1)
\end{eqnarray*}%
by recalling that $\eta _{T}^{2}\approx Tl(\eta _{T})$ for large $T$. Then,
for any $t=0,\cdots ,[|\nu |g(T)]-1$
\begin{eqnarray}
\frac{y_{[\tau _{0}T]-t-1}}{\sqrt{g(T)l(\eta _{T})}} &=&\frac{1}{\sqrt{%
g(T)l(\eta _{T})}}\sum_{i=1}^{[\tau _{0}T]-t-1}\beta _{1T}^{[\tau
_{0}T]-t-1-i}\varepsilon _{i}^{(1)}+o_{p}(1)  \notag  \label{2.16} \\
&\Rightarrow &\int_{0}^{\infty }e^{-s}dW_{1}(s)\overset{d}{=} B_{a}(\frac{1}{%
2}),
\end{eqnarray}%
see page 138 in Chong (2001) for details. Note also that
\begin{equation*}
\frac{1}{\sqrt{g(T)l(\eta _{T})}}\sum_{t=0}^{[|\nu |g(T)]-1}\varepsilon
_{\lbrack \tau _{0}T]-t}\Rightarrow W_{1}(|\nu |)
\end{equation*}%
by functional central limit theorem for i.i.d. random variables from DAN.
As a result, we have
\begin{equation}  \label{2.9}
\frac{1}{g(T)l(\eta _{T})}\sum_{t=0}^{[|\nu |g(T)]-1}y_{[\tau
_{0}T]-t-1}\varepsilon _{\lbrack \tau _{0}T]-t}\Rightarrow B_{a}(\frac{1}{2}%
)W_{1}(|\nu |)
\end{equation}%
and
\begin{equation}  \label{2.11}
\frac{1}{g^{2}(T)l(\eta _{T})}\sum_{t=0}^{[|\nu |g(T)]-1}y_{[\tau
_{0}T]-t-1}^{2}\Rightarrow |\nu |B_{a}^{2}(\frac{1}{2}).
\end{equation}%
Moreover, note that
\begin{equation}  \label{2.36}
\frac{(\beta _{2T}-\beta _{1T})}{l(\eta _{T})}\frac{\sum_{t=[\tau
T]+1}^{[\tau _{0}T]}y_{t-1}^{2}\sum_{t=[\tau _{0}T]+1}^{T}y_{t-1}\varepsilon
_{t}}{\sum_{t=[\tau T]+1}^{T}y_{t-1}^{2}}=\frac{O_{p}(g^2(T)l(\eta
_{T}))O_{p}(Tl(\eta _{T}))}{g(T)l(\eta _{T})O_{p}(T^{2}l(\eta _{T}))}%
=o_{p}(1)
\end{equation}%
and
\begin{equation}  \label{2.37}
\frac{\sum_{t=[\tau _{0}T]+1}^{T}y_{t-1}^{2}}{\sum_{t=[\tau
T]+1}^{T}y_{t-1}^{2}}\overset{p}{\rightarrow } 1.
\end{equation}%
Then, by recalling equation B.2 on page 117 in Chong (2001) when $\underline{%
\tau }\leq \tau \leq \tau _{0}$
\begin{eqnarray*}
&&RSS_{T}(\tau )-RSS_{T}(\tau _{0}) \\
&=&2(\beta _{2T}-\beta _{1T})\bigg(\frac{\sum_{t=[\tau T]+1}^{[\tau
_{0}T]}y_{t-1}^{2}\sum_{t=[\tau _{0}T]+1}^{T}y_{t-1}\varepsilon _{t}}{%
\sum_{t=[\tau T]+1}^{T}y_{t-1}^{2}}-\frac{\sum_{t=[\tau
_{0}T]+1}^{T}y_{t-1}^{2}\sum_{t=[\tau T]+1}^{[\tau _{0}T]}y_{t-1}\varepsilon
_{t}}{\sum_{t=[\tau T]+1}^{T}y_{t-1}^{2}}\bigg) \\
&&+(\beta _{2T}-\beta _{1T})^{2}\frac{\sum_{t=[\tau
_{0}T]+1}^{T}y_{t-1}^{2}\sum_{t=[\tau T]+1}^{[\tau _{0}T]}y_{t-1}^{2}}{%
\sum_{t=[\tau T]+1}^{T}y_{t-1}^{2}}+\Lambda _{T}(\tau ).
\end{eqnarray*}%
We have
\begin{eqnarray*}
&&\frac{RSS_{T}(\tau )-RSS_{T}(\tau _{0})}{l(\eta _{T})} \\
&=&-\frac{2}{g(T)l(\eta _{T})}\sum_{t=[\tau T]+1}^{[\tau
_{0}T]}y_{t-1}\varepsilon _{t}\cdot (1+o_{p}(1))+\frac{1}{g^{2}(T)l(\eta
_{T})}\sum_{t=[\tau T]+1}^{[\tau _{0}T]}y_{t-1}^{2}\cdot
(1+o_{p}(1))+o_{p}(1) \\
&\Rightarrow &-2B_{a}(\frac{1}{2})W_{1}(|\nu |)+|\nu |B_{a}^{2}(\frac{1}{2})
\end{eqnarray*}%
by (\ref{2.9})-(\ref{2.37}).

Similarly, for $\tau =\tau _{0}+\nu g(T)/T$ with $\nu >0$, we also have $%
\Lambda _{T}(\tau )=o_{p}(l(\eta _{T}))$. Moreover, it can be shown that
\begin{eqnarray}
&&\frac{(\beta _{2T}-\beta _{1T})}{l(\eta _{T})}\frac{\sum_{t=[\tau
_{0}T]+1}^{[\tau T]}y_{t-1}^{2}\sum_{t=1}^{[\tau _{0}T]}y_{t-1}\varepsilon
_{t}}{\sum_{t=1}^{[\tau T]}y_{t-1}^{2}}  \notag  \label{2.42} \\
&=&\frac{O_{p}(g^{2}(T)l(\eta _{T}))O_{p}(\sqrt{T g(T)} l(\eta _{T}))}{%
g(T)l(\eta _{T})\big(O_{p}(T g(T)l(\eta _{T}))+O_{p}(g^{2}(T)l(\eta _{T}))%
\big)}  \notag \\
&=&O_{p}\Big(\sqrt{\frac{g(T)}{T}}\Big)=o_{p}(1);
\end{eqnarray}

\begin{equation*}
\frac{\sum_{t=1}^{[\tau T]}y_{t-1}^{2}}{\sum_{t=1}^{[\tau _{0}T]}y_{t-1}^{2}}%
=1+\frac{\sum_{t=[\tau _{0}T]+1}^{[\tau T]}y_{t-1}^{2}}{\sum_{t=1}^{[\tau
_{0}T]}y_{t-1}^{2}}=1+\frac{O_{p}(g^{2}(T)l(\eta _{T}))}{O_{p}(Tg(T)l(\eta
_{T}))}=1+o_{p}(1);
\end{equation*}%
\begin{eqnarray*}
\frac{y_{[\tau _{0}T]}}{\sqrt{g(T)l(\eta _{T})}} &=&\Big(1-\frac{c}{T}-\frac{%
1}{g(T)}\Big)^{[\tau _{0}T]}\frac{y_{0}}{\sqrt{g(T)l(\eta _{T})}}%
+\sum_{t=0}^{[\tau _{0}T]-1}\Big(1-\frac{c}{T}-\frac{1}{g(T)}\Big)^{t}\frac{%
\varepsilon _{\lbrack \tau _{0}T]-t}}{\sqrt{g(T)l(\eta _{T})}} \\
&\Rightarrow &\int_{0}^{\infty }\exp {(-s)}dW_{1}(s)\overset{d}{=} B_{a}(%
\frac{1}{2})
\end{eqnarray*}%
whose proof is similar to those of (\ref{3.7}) and (\ref{2.16}) ;
\begin{eqnarray*}
&&\frac{1}{g^{2}(T)l(\eta _{T})}\sum_{t=0}^{[\nu g(T)]-1}y_{[\tau
_{0}T]+t}^{2} \\
&=&\frac{1}{g(T)}\sum_{t=0}^{[\nu g(T)]-1}\Big(\frac{1}{\sqrt{g(T)l(\eta
_{T})}}\sum_{i=0}^{t-1}\beta _{2T}^{i}\varepsilon _{\lbrack \tau _{0}T]+t-i}+%
\frac{\beta _{2T}^{t}y_{[\tau _{0}T]}}{\sqrt{g(T)l(\eta _{T})}}\Big)^{2} \\
&\Rightarrow &\int_{0}^{\nu }\Big(I(W_{2},c,\tau _{0},t)+B_{a}(\frac{1}{2})%
\Big)^{2}dt
\end{eqnarray*}%
and
\begin{equation*}
\frac{1}{\sqrt{g(T)}}\sum_{t=0}^{[\nu g(T)]-1}\frac{y_{[\tau _{0}T]+t}}{%
\sqrt{g(T)l(\eta _{T})}}\frac{\varepsilon _{\lbrack \tau _{0}T]+t+1}}{\sqrt{%
l(\eta _{T})}}\Rightarrow \int_{0}^{\nu }\Big(I(W_{2},c,\tau _{0},t)+B_{a}(%
\frac{1}{2})\Big)dI(W_{2},c,\tau _{0},t)
\end{equation*}%
by virtue of Lemma \ref{lem6}. Thus, by recalling the equation B.4 on page
120 in Chong (2001) when $\tau _{0}\leq \tau \leq \overline{\tau }$
\begin{eqnarray*}
&&RSS_{T}(\tau )-RSS_{T}(\tau _{0}) \\
&=&2(\beta _{2T}-\beta _{1T})\bigg(\frac{\sum_{t=1}^{[\tau
_{0}T]}y_{t-1}^{2}\sum_{t=[\tau _{0}T]+1}^{[\tau T]}y_{t-1}\varepsilon _{t}}{%
\sum_{t=1}^{[\tau T]}y_{t-1}^{2}}-\frac{\sum_{t=[\tau _{0}T]+1}^{[\tau
T]}y_{t-1}^{2}\sum_{t=1}^{[\tau _{0}T]}y_{t-1}\varepsilon _{t}}{%
\sum_{t=1}^{[\tau T]}y_{t-1}^{2}}\bigg) \\
&&+(\beta _{2T}-\beta _{1T})^{2}\frac{\sum_{t=[\tau _{0}T]+1}^{[\tau
T]}y_{t-1}^{2}\sum_{t=1}^{[\tau _{0}T]}y_{t-1}^{2}}{\sum_{t=1}^{[\tau
T]}y_{t-1}^{2}}+\Lambda _{T}(\tau ),
\end{eqnarray*}%
we have
\begin{eqnarray*}
&&\frac{RSS_{T}(\tau )-RSS_{T}(\tau _{0})}{l(\eta _{T})} \\
&=&\frac{2(\beta _{2T}-\beta _{1T})}{l(\eta _{T})}\cdot \sum_{t=[\tau
_{0}T]+1}^{[\tau T]}y_{t-1}\varepsilon _{t}\cdot (1+o_{p}(1))+\frac{(\beta
_{2T}-\beta _{1T})^{2}}{l(\eta _{T})}\sum_{t=[\tau _{0}T]+1}^{[\tau
T]}y_{t-1}^{2}\cdot (1+o_{p}(1))+o_{p}(1) \\
&=&\frac{2(1+o_{p}(1))}{g(T)l(\eta _{T})}\sum_{t=0}^{[\nu g(T)]-1}y_{[\tau
_{0}T]+t}\varepsilon _{\lbrack \tau _{0}T]+t+1}+\frac{(1+o_{p}(1))}{%
g^{2}(T)l(\eta _{T})}\sum_{t=0}^{[\nu g(T)]-1}y_{[\tau _{0}T]+t}^{2}+o_{p}(1)
\\
&\Rightarrow &2\int_{0}^{\nu }\Big(I(W_{2},c,\tau _{0},t)+B_{a}(\frac{1}{2})%
\Big)dI(W_{2},c,\tau _{0},t)+\int_{0}^{\nu }\Big(I(W_{2},c,\tau
_{0},t)+B_{a}(\frac{1}{2})\Big)^{2}dt \\
&=&-2B_{a}^{2}(\frac{1}{2})\cdot \Bigg\{-\frac{I(W_{2},c,\tau _{0},\nu )}{%
B_{a}(\frac{1}{2})}-\int_{0}^{\nu }\frac{I(W_{2},c,\tau _{0},t)}{B_{a}^{2}(%
\frac{1}{2})}dI(W_{2},c,\tau _{0},t)- \\
&&~~~~~~~~~~~~~~~~~~~~~\int_{0}^{\nu }\Big(\frac{I(W_{2},c,\tau _{0},t)}{%
2B_{a}(\frac{1}{2})}+1\Big)\frac{I(W_{2},c,\tau _{0},t)}{B_{a}(\frac{1}{2})}%
dt-\frac{\nu }{2}\Bigg\}.
\end{eqnarray*}%
Applying the continuous mapping theorem for argmax functionals (cf. Kim and
Pollard (1990)), we have
\begin{eqnarray*}
(\beta _{2T}-\beta _{1T})T(\hat{\tau}_{T}-\tau _{0}) &=&\hat{\nu}=%
\mathop{\arg\min}_{\nu \in R}\big(RSS_{T}(\tau )-RSS_{T}(\tau _{0})\big) \\
&=&\mathop{\arg\min}_{\nu \in R}\Big(\frac{RSS_{T}(\tau )-RSS_{T}(\tau _{0})%
}{l(\eta _{T})}\Big) \\
&\Rightarrow &\mathop{\arg\min}_{\nu \in R}\bigg\{-2B_{a}^{2}(\frac{1}{2})%
\bigg(\frac{C^{\ast }(\nu )}{B_{a}(\frac{1}{2})}-\frac{|\nu |}{2}\bigg)%
\bigg\} \\
&=&\mathop{\arg\max}_{\nu \in R}\bigg\{\frac{C^{\ast }(\nu )}{B_{a}(\frac{1}{%
2})}-\frac{|\nu |}{2}\bigg\},
\end{eqnarray*}%
where $C^{\ast }(\nu )$ is defined as in Theorem \ref{thm}. The proofs are
complete. $\hfill \Box $\newline

\bigskip

\section{Proof of Theorem \protect\ref{thm2}}

\setcounter{equation}{0}

The following lemmas will be used in the proof of Theorem \ref{thm2}:

\begin{lemma}
\label{lem8} Suppose assumption C1 is satisfied and $\beta _{1}=\beta
_{1T}=1-c/T$, where $c$ is a fixed constant, then for any $0\leq \tau \leq \tau _{0}$,
\begin{equation*}
Q_{T}(\tau ):=\frac{1}{\sqrt{Tl(\eta _{T})}}\sum_{t=1}^{[\tau T]}\beta
_{1T}^{T-t}\varepsilon _{t}\Rightarrow e^{-c(1-\tau )}W(\tau
)-c\int_{0}^{\tau }e^{-c(1-s)}W(s)ds=G(W,c,\tau ).
\end{equation*}
\end{lemma}

\noindent\textbf{Proof.} \; Denote
\begin{eqnarray*}
U_T(t)=\frac{1}{\sqrt{T l(\eta_T)}}\sum_{i=1}^{[t T]}\varepsilon_i,
\end{eqnarray*}
then it follows from the functional central limit theorem for i.i.d.
sequence from DAN that
\begin{eqnarray*}
U_T(t)\Rightarrow W(t).
\end{eqnarray*}
Now, by denoting $\beta_{1T}=e^{-c_T^{*}/T}$ with $c_T^{*}/c\rightarrow 1$,
we have
\begin{eqnarray*}
\frac{1}{\sqrt{T l(\eta_T)}}\sum_{t=1}^{[\tau
T]}\beta_{1T}^{T-t}\varepsilon_t^{(1)}&=&\sum_{t=1}^{[\tau T]}e^{-c_T^{*}%
\frac{T-t}{T}}\int_{\frac{t-1}{T}}^{\frac{t}{T}}d U_T(s)  \notag \\
&=&\Big(\sum_{t=1}^{[\tau T]}\int_{\frac{t-1}{T}}^{\frac{t}{T}%
}e^{-c_T^{*}(1-s)}d U_T(s)\Big)\cdot (1+o_p(1))  \notag \\
&=&\Big(\int_{0}^{\tau}e^{-c_T^{*}(1-s)}d U_T(s)\Big)\cdot (1+o_p(1))  \notag
\\
&=&\Big(e^{-c_T^{*}(1-\tau)}U_T(\tau)-c_T^{*}\int_{0}^{%
\tau}e^{-c_T^{*}(1-s)}U_T(s)dt\Big)\cdot (1+o_p(1))  \notag \\
&\Rightarrow&e^{-c(1-\tau)}W(\tau)-c\int_{0}^{\tau}e^{-c(1-s)}W(s)ds,
\end{eqnarray*}
as desired. $\hfill \Box$\newline

\begin{lemma}
\label{lem3} Let $\{y_{t}\}$ be generated according to Model (\ref{model1}),
where $\beta _{1}=\beta _{1T}=1-c/T$ for a constant $c$ and $|\beta _{2}|<1$%
. Under assumptions C1-C3, the following results hold jointly:\newline
\textrm{(4a)}~$\frac{1}{Tl(\eta _{T})}\sum_{t=1}^{[\tau
_{0}T]}y_{t-1}\varepsilon _{t}\Rightarrow \frac{1}{2}e^{2c(1-\tau
_{0})}G^{2}(W,c,\tau _{0})+c\int_{0}^{\tau _{0}}e^{2c(1-t)}G^{2}(W,c,t)dt-%
\frac{\tau _{0}}{2},$ \newline
\textrm{(4b)}~$\frac{1}{T^{2}l(\eta _{T})}\sum_{t=1}^{[\tau
_{0}T]}y_{t-1}^{2}\Rightarrow \int_{0}^{\tau _{0}}e^{2c(1-t)}G^{2}(W,c,t)dt.$
\end{lemma}

\noindent \textbf{Proof.} \ First, it can be shown that
\begin{equation*}
2(1-\frac{c}{T})\sum_{t=1}^{[\tau _{0}T]}y_{t-1}\varepsilon _{t}=y_{[\tau
_{0}T]}^{2}-y_{0}^{2}+(\frac{2c}{T}-\frac{c^{2}}{T^{2}})\sum_{t=1}^{[\tau
_{0}T]}y_{t-1}^{2}-\sum_{t=1}^{[\tau _{0}T]}\varepsilon _{t}^{2},
\end{equation*}%
which is equivalent to
\begin{eqnarray}
\frac{\sum_{t=1}^{[\tau _{0}T]}y_{t-1}\varepsilon _{t}}{Tl(\eta _{T})} &=&%
\frac{T}{2(T-c)}\frac{y_{[\tau _{0}T]}^{2}-y_{0}^{2}}{Tl(\eta _{T})}+\frac{%
c(2T-c)}{2T(T-c)}\sum_{t=1}^{[\tau _{0}T]}\Big(\frac{y_{t-1}}{\sqrt{Tl(\eta
_{T})}}\Big)^{2}  \notag  \label{3.10} \\
&&-\frac{T}{2(T-c)}\frac{1}{Tl(\eta _{T})}\sum_{t=1}^{[\tau
_{0}T]}\varepsilon _{t}^{2}.
\end{eqnarray}

Consider the first term in the right hand side of (\ref{3.10}). Clearly,
\begin{equation}  \label{3.12}
\frac{y_{0}}{\sqrt{Tl(\eta _{T})}}\overset{p}{\rightarrow } 0
\end{equation}%
and
\begin{equation*}
\frac{y_{[\tau _{0}T]}^{2}}{Tl(\eta _{T})}=\Big(\frac{\sum_{t=1}^{[\tau
_{0}T]}\beta _{1T}^{[\tau _{0}T]-t}\varepsilon _{t}+\beta _{1T}^{[\tau
_{0}T]}y_{0}}{\sqrt{Tl(\eta _{T})}}\Big)^{2}\Rightarrow e^{2c(1-\tau
_{0})}G^{2}(W,c,\tau _{0})
\end{equation*}%
by Lemma \ref{lem8} and the fact $|\beta _{1T}^{[\tau _{0}T]}|<C$. Thus,
\begin{equation}  \label{3.11}
\frac{T}{2(T-c)}\frac{y_{[\tau _{0}T]}^{2}-y_{0}^{2}}{Tl(\eta _{T})}%
\Rightarrow \frac{1}{2}e^{2c(1-\tau _{0})}G^{2}(W,c,\tau _{0}).
\end{equation}

Consider the second term in the right hand side of (\ref{3.10}). Write
\begin{eqnarray}
\frac{1}{T}\sum_{t=1}^{[\tau _{0}T]}\Big(\frac{y_{t-1}}{\sqrt{Tl(\eta _{T})}}%
\Big)^{2} &=&\frac{1}{T}\sum_{t=1}^{[\tau _{0}T]}\Big(\frac{%
\sum_{i=1}^{t-1}\beta _{1T}^{t-1-i}\varepsilon _{i}}{\sqrt{Tl(\eta _{T})}}+%
\frac{\beta _{1T}^{t-1}y_{0}}{\sqrt{Tl(\eta _{T})}}\Big)^{2}  \notag
\label{3.15} \\
 &:=&\frac{1}{T}\sum_{t=1}^{[\tau _{0}T]}(V+VI)^{2}.
\end{eqnarray}%
Invoke (\ref{3.12}), we have
\begin{equation*}
\frac{1}{T}\sum_{t=1}^{[\tau _{0}T]}(VI)^{2}\overset{p}{\rightarrow } 0
\end{equation*}%
since
\begin{equation*}
\frac{1}{T}\sum_{t=1}^{[\tau _{0}T]}\beta _{1T}^{2(t-1)}=\frac{1-\beta
_{1T}^{2[\tau _{0}T]}}{T(1-\beta _{1T}^{2})}\rightarrow \frac{1-e^{-2c\tau
_{0}}}{2c}<\infty .
\end{equation*}%
For the term $\frac{1}{T}\sum_{t=1}^{[\tau _{0}T]}(V)^{2}$, write
\begin{eqnarray}
\frac{1}{T}\sum_{t=1}^{[\tau _{0}T]}(V)^{2} &=&\frac{1}{T}\sum_{t=1}^{[\tau
_{0}T]}\Big(\frac{\sum_{i=1}^{t-1}\beta _{1T}^{t-1-i}\varepsilon _{i}}{\sqrt{%
Tl(\eta _{T})}}\Big)^{2}  \notag  \label{3.14} \\
&=&\frac{1}{T}\sum_{t=1}^{[\tau _{0}T]}\Big(\frac{\sum_{i=1}^{t-1}\beta
_{1T}^{T-i}\varepsilon _{i}}{\sqrt{Tl(\eta _{T})}}\Big)^{2}\beta
_{1T}^{2(t-1-T)}-\int_{0}^{\tau _{0}}e^{2c(1-t)}Q_{T}^{2}(t)dt  \notag \\
&&+\int_{0}^{\tau _{0}}e^{2c(1-t)}Q_{T}^{2}(t)dt  \notag \\
 &:=&V_{T}(\tau _{0})+\int_{0}^{\tau _{0}}e^{2c(1-t)}Q_{T}^{2}(t)dt.
\end{eqnarray}%
Next, we shall show\ that
\begin{equation}
V_{T}(\tau _{0})\overset{p}{\rightarrow } 0.  \label{3.13}
\end{equation}%
Note that
\begin{eqnarray*}
|V_{T}(\tau _{0})| &=&\Big|\frac{1}{T}\sum_{t=1}^{[\tau _{0}T]}\beta
_{1T}^{2(t-1-T)}Q_{T}^{2}(\frac{t-1}{T})-\int_{0}^{\tau
_{0}}e^{2c(1-t)}Q_{T}^{2}(t)dt\Big| \\
&=&\Big|\sum_{t=1}^{[\tau _{0}T]}\int_{\frac{t-1}{T}}^{\frac{t}{T}}\beta
_{1T}^{2(t-1-T)}Q_{T}^{2}(s)ds-\sum_{t=1}^{[\tau _{0}T]}\int_{\frac{t-1}{T}%
}^{\frac{t}{T}}e^{2c(1-s)}Q_{T}^{2}(s)ds\Big|+o_{p}(1) \\
&\leq &\sum_{t=1}^{[\tau _{0}T]}\int_{\frac{t-1}{T}}^{\frac{t}{T}}|\beta
_{1T}^{2(t-1-T)}-e^{2c(1-s)}|Q_{T}^{2}(s)ds+o_{p}(1) \\
&\leq &\max_{1\leq t\leq \lbrack \tau _{0}T]}\sup_{(t-1)/T\leq s\leq
t/T}|\beta _{1T}^{2(t-1-T)}-e^{2c(1-s)}|\cdot \sum_{t=1}^{[\tau _{0}T]}\int_{%
\frac{t-1}{T}}^{\frac{t}{T}}Q_{T}^{2}(s)ds+o_{p}(1).
\end{eqnarray*}%
Following the proof of (2.7) in Chan and Wei (1987) again, we have
\begin{equation*}
\max_{1\leq t\leq \lbrack \tau _{0}T]}\sup_{(t-1)/T\leq s\leq t/T}|\beta
_{1T}^{2(t-1-T)}-e^{2c(1-s)}|\rightarrow 0.
\end{equation*}%
Meanwhile, it follows from Lemma \ref{lem8} that
\begin{equation*}
\sum_{t=1}^{[\tau _{0}T]}\int_{\frac{t-1}{T}}^{\frac{t}{T}%
}Q_{T}^{2}(s)ds\Rightarrow \int_{0}^{\tau _{0}}G^{2}(W,c,s)ds.
\end{equation*}%
This proves (\ref{3.13}). Further, invoke (\ref{3.15})-(\ref{3.13}), Lemma %
\ref{lem8} and Cauchy-Schwarz inequality, we have
\begin{equation}
\frac{1}{T}\sum_{t=1}^{[\tau _{0}T]}\Big(\frac{y_{t-1}}{\sqrt{Tl(\eta _{T})}}%
\Big)^{2}\Rightarrow \int_{0}^{\tau _{0}}e^{2c(1-t)}G^{2}(W,c,t)dt.
\label{3.19}
\end{equation}%
and
\begin{equation}
\frac{c(2T-c)}{2T(T-c)}\sum_{t=1}^{[\tau _{0}T]}\Big(\frac{y_{t-1}}{\sqrt{%
Tl(\eta _{T})}}\Big)^{2}\Rightarrow c\int_{0}^{\tau
_{0}}e^{2c(1-t)}G^{2}(W,c,t)dt.  \label{3.17}
\end{equation}

For the third term in (\ref{3.10}), it is obvious that
\begin{equation}
\frac{T}{2(T-c)}\frac{1}{Tl(\eta _{T})}\sum_{t=1}^{[\tau _{0}T]}\varepsilon
_{t}^{2}\overset{p}{\rightarrow } \frac{\tau _{0}}{2}  \label{3.18}
\end{equation}%
when $\varepsilon _{t}^{\prime }s$ are i.i.d. random variables and from DAN.
Combining (\ref{3.11}), (\ref{3.17}) and (\ref{3.18}) together lead to (4a).
(4b) is obvious from (\ref{3.19}). It is easy to see that (4a) and (4b) hold
jointly. The proofs are complete. $\hfill \Box $\newline

\begin{lemma}
\label{lem4} Let $\{y_{t}\}$ be generated according to Model (\ref{model1}),
where $\beta _{1}=\beta _{1T}=1-c/T$ for a fixed constant $c$ and $|\beta
_{2}|<1$, Under assumptions C1-C3, the following results hold jointly:%
\newline
\textrm{(5a)}~$\frac{1}{\sqrt{T}l(\eta _{T})}\sum_{t=[\tau
_{0}T]+1}^{T}y_{t-1}\varepsilon _{t}\Rightarrow \overline{W}(B(c,\tau _{0}))/%
\sqrt{1-\beta _{2}^{2}},$\newline
\textrm{(5b)}~$\frac{1}{Tl(\eta _{T})}\sum_{t=[\tau
_{0}T]+1}^{T}y_{t-1}^{2}\Rightarrow \frac{1-\tau _{0}+e^{2c(1-\tau
_{0})}G^{2}(W,c,\tau _{0})}{1-\beta _{2}^{2}},$\newline
where $B(c,\tau _{0})=(1-e^{-2c\tau _{0}})/(2c)+1-\tau _{0}.$
\end{lemma}

\noindent \textbf{Proof.} \ To prove (5a), we make use of the following
decomposition,
\begin{eqnarray*}
&&\frac{1}{\sqrt{T}l(\eta _{T})}\sum_{t=[\tau
_{0}T]+1}^{T}y_{t-1}\varepsilon _{t} \\
&=&\frac{1}{\sqrt{T}l(\eta _{T})}\sum_{t=[\tau _{0}T]+1}^{T}\Big(\beta
_{2}^{t-[\tau _{0}T]-1}y_{[\tau _{0}T]}+\sum_{i=[\tau _{0}T]+1}^{t-1}\beta
_{2}^{t-i-1}\varepsilon _{i}\Big)\varepsilon _{t} \\
&=&\frac{1}{\sqrt{T}l(\eta _{T})}\Big(\Big(\beta _{1T}^{[\tau
_{0}T]}y_{0}+\sum_{j=1}^{[\tau _{0}T]}\beta _{1T}^{[\tau
_{0}T]-j}\varepsilon _{j}\Big)\sum_{t=[\tau _{0}T]+1}^{T}\beta _{2}^{t-[\tau
_{0}T]-1}\varepsilon _{t}+\sum_{t=[\tau _{0}T]+1}^{T}\varepsilon
_{t}\sum_{i=[\tau _{0}T]+1}^{t-1}\beta _{2}^{t-i-1}\varepsilon _{i}\Big) \\
&:=& VII+VIII,
\end{eqnarray*}%
where
\begin{equation*}
VII=\frac{1}{\sqrt{T}l(\eta _{T})}\Big(\beta _{1T}^{[\tau
_{0}T]}y_{0}+\sum_{j=1}^{[\tau _{0}T]}\beta _{1T}^{[\tau
_{0}T]-j}\varepsilon _{j}\Big)\sum_{t=[\tau _{0}T]+1}^{T}\beta _{2}^{t-[\tau
_{0}T]-1}\varepsilon _{t}
\end{equation*}%
and
\begin{equation*}
VIII=\frac{1}{\sqrt{T}l(\eta _{T})}\sum_{t=[\tau _{0}T]+1}^{T}\varepsilon
_{t}\sum_{i=[\tau _{0}T]+1}^{t-1}\beta _{2}^{t-i-1}\varepsilon _{i}.
\end{equation*}

Consider the term $VII$ first. Since $|\beta _{1T}^{[\tau _{0}T]}|<\infty $,
$y_{0}=o_p(\sqrt{T})$ and $E|\sum_{i=[\tau _{0}T]+1}^{t-1}\beta
_{2}^{t-i-1}\varepsilon _{i}|<\infty $, we have
\begin{equation*}
\frac{1}{\sqrt{T}l(\eta _{T})}\beta _{1T}^{[\tau _{0}T]}y_{0}\sum_{t=[\tau
_{0}T]+1}^{T}\beta _{2}^{t-[\tau _{0}T]-1}\varepsilon _{t}=o_{p}(1).
\end{equation*}%
Moreover, note that
\begin{eqnarray}
&&\frac{1}{\sqrt{T}l(\eta _{T})}\sum_{j=1}^{[\tau _{0}T]}\beta _{1T}^{[\tau
_{0}T]-j}\varepsilon _{j}\sum_{t=[\tau _{0}T]+1}^{T}\beta _{2}^{t-[\tau
_{0}T]-1}\varepsilon _{t}  \notag  \label{3.4} \\
&=&\frac{1}{\sqrt{T}l(\eta _{T})}\sum_{j=1}^{[\tau _{0}T]}\beta _{1T}^{[\tau
_{0}T]-j}(\varepsilon _{j}^{(1)}+\varepsilon _{j}^{(2)})\sum_{t=[\tau
_{0}T]+1}^{T}\beta _{2}^{t-[\tau _{0}T]-1}(\varepsilon
_{t}^{(1)}+\varepsilon _{t}^{(2)}).
\end{eqnarray}%
%
%
%
%
%
%
%
Applying Lemma \ref{lem1}, the leading term of (\ref{3.4}) will be
\begin{equation*}
\frac{1}{\sqrt{T}}\sum_{j=1}^{[\tau _{0}T]}\beta _{1T}^{[\tau _{0}T]-j}\frac{%
\varepsilon _{j}^{(1)}}{\sqrt{l(\eta _{T})}}\sum_{t=[\tau _{0}T]+1}^{T}\beta
_{2}^{t-[\tau _{0}T]-1}\frac{\varepsilon _{t}^{(1)}}{\sqrt{l(\eta _{T})}}
\end{equation*}%
and the other three terms are negligible. Note that $\varepsilon _{j}^{(1)}/%
\sqrt{l(\eta _{T})}$, for $j=1,\cdots ,T$, has zero mean and finite
variance, it follows from the central limit theorem for martingale
differences that
\begin{equation*}
\frac{1}{\sqrt{T}}\sum_{j=1}^{[\tau _{0}T]}\beta _{1T}^{[\tau _{0}T]-j}\frac{%
\varepsilon _{j}^{(1)}}{\sqrt{l(\eta _{T})}}\sum_{t=[\tau _{0}T]+1}^{T}\beta
_{2}^{t-[\tau _{0}T]-1}\frac{\varepsilon _{t}^{(1)}}{\sqrt{l(\eta _{T})}}%
\Rightarrow N\Big(0,\frac{1-e^{-2c\tau _{0}}}{2c(1-\beta _{2}^{2})}\Big),
\end{equation*}%
which yields
\begin{equation*}
VII\Rightarrow N\Big(0,\frac{1-e^{-2c\tau _{0}}}{2c(1-\beta _{2}^{2})}\Big).
\end{equation*}%
Similarly, applying the central limit theorem for martingale differences
again, we have
\begin{equation*}
VIII\Rightarrow N\Big(0,\frac{1-\tau _{0}}{1-\beta _{2}^{2}}\Big).
\end{equation*}%
By the independence of the two martingale difference sequences given
previously, we have
\begin{equation*}
\frac{1}{\sqrt{T}l(\eta _{T})}\sum_{t=[\tau _{0}T]+1}^{T}y_{t-1}\varepsilon
_{t}\Rightarrow \overline{W}\Big(\frac{1-e^{-2c\tau _{0}}}{2c(1-\beta
_{2}^{2})}+\frac{1-\tau _{0}}{1-\beta _{2}^{2}}\Big),
\end{equation*}%
which implies (5a).

To prove (5b), note that
\begin{eqnarray}
&&\frac{1}{Tl(\eta _{T})}\sum_{t=[\tau _{0}T]+1}^{T}y_{t-1}^{2}  \notag
\label{3.20} \\
&=&\frac{1}{Tl(\eta _{T})}\sum_{t=[\tau _{0}T]+1}^{T}\Big(\beta
_{2}^{t-[\tau _{0}T]-1}y_{[\tau _{0}T]}+\sum_{i=[\tau _{0}T]+1}^{t-1}\beta
_{2}^{t-i-1}\varepsilon _{i}\Big)^{2}  \notag \\
&=&\Big(\frac{y_{[\tau _{0}T]}}{\sqrt{Tl(\eta _{T})}}\Big)^{2}\sum_{t=[\tau
_{0}T]+1}^{T}\beta _{2}^{2(t-[\tau _{0}T]-1)}+\frac{2y_{[\tau _{0}T]}}{%
Tl(\eta _{T})}\sum_{t=[\tau _{0}T]+1}^{T}\big(\beta _{2}^{t-[\tau
_{0}T]-1}\sum_{i=[\tau _{0}T]+1}^{t-1}\beta _{2}^{t-i-1}\varepsilon _{i}\big)
\notag \\
&&+\frac{1}{Tl(\eta _{T})}\sum_{t=[\tau _{0}T]+1}^{T}\Big(\sum_{i=[\tau
_{0}T]+1}^{t-1}\beta _{2}^{t-i-1}\varepsilon _{i}\Big)^{2}.
\end{eqnarray}

Consider the first term in the right hand side of (\ref{3.20}). Since
\begin{eqnarray*}
\frac{y_{[\tau_0 T] }}{\sqrt{T l(\eta_T)}}=\frac{\beta_{1T}^{[\tau_0
T]}y_0+\sum_{t=1}^{[\tau_0 T]}\beta_{1T}^{[\tau_0 T]-t}\varepsilon_t}{\sqrt{%
T l(\eta_T)}}=\frac{\beta_{1T}^{[\tau_0 T]}y_0+\beta_{1T}^{[\tau_0
T]-T}\sum_{t=1}^{[\tau_0 T]}\beta_{1T}^{T-t}\varepsilon_t}{\sqrt{T l(\eta_T)}%
}
\end{eqnarray*}
and $\frac{\beta_{1T}^{[\tau_0 T]}y_0}{\sqrt{T l(\eta_T)}}\overset{p}{%
\rightarrow} 0$, it follows from Lemma \ref{lem8} that
\begin{eqnarray}  \label{3.3}
\frac{y_{[\tau_0 T] }}{\sqrt{T l(\eta_T)}}\Rightarrow e^{c(1-\tau_0)}G(W, c,
\tau_0),
\end{eqnarray}
which further implies that
\begin{eqnarray}  \label{3.5}
\Big(\frac{y_{[\tau_0 T] }}{\sqrt{T l(\eta_T)}}\Big)^2\sum_{t=[\tau_0 T]
+1}^{T}\beta_2^{2(t-[\tau_0 T] -1)}\Rightarrow \frac{e^{2c(1-\tau_0)}}{%
1-\beta_2^2}G^2(W, c, \tau_0).
\end{eqnarray}

For the second term of (\ref{3.20}), note that (\ref{3.3}) and
\begin{eqnarray*}
&&E\Big|\sum_{t=[\tau _{0}T]+1}^{T}\big(\beta _{2}^{t-[\tau
_{0}T]-1}\sum_{i=[\tau _{0}T]+1}^{t-1}\beta _{2}^{t-i-1}\varepsilon _{i}\big)%
\Big| \\
&\leq &C\Big|\sum_{t=[\tau _{0}T]+1}^{T}|\beta _{2}|^{t-[\tau
_{0}T]-1}\sum_{i=[\tau _{0}T]+1}^{t-1}|\beta _{2}|^{t-i-1}\Big| \\
&\leq &\frac{C}{(1-|\beta _{2}|)^{2}}+\frac{C}{(1-|\beta _{2}|)(1-\beta
_{2}^{2})}=O(1).
\end{eqnarray*}%
We have
\begin{equation}
\frac{2y_{[\tau _{0}T]}}{Tl(\eta _{T})}\sum_{t=[\tau _{0}T]+1}^{T}\big(\beta
_{2}^{t-[\tau _{0}T]-1}\sum_{i=[\tau _{0}T]+1}^{t-1}\beta
_{2}^{t-i-1}\varepsilon _{i}\big)=o_{p}(1).  \label{3.2}
\end{equation}

Finally, one can show, by the truncation arguments used in (\ref{3.4}), that
\begin{equation}  \label{3.6}
\frac{1}{Tl(\eta _{T})}\sum_{t=[\tau _{0}T]+1}^{T}\Big(\sum_{i=[\tau
_{0}T]+1}^{t-1}\beta _{2}^{t-i-1}\varepsilon _{i}\Big)^{2}\overset{p}{%
\rightarrow } \frac{1-\tau _{0}}{1-\beta _{2}^{2}}.
\end{equation}%
Thus, (5b) is proved by combining (\ref{3.5}), (\ref{3.2}) and (\ref{3.6}).
Note that (5a) and (5b) hold jointly and that the limiting behavior of (5b)
is determined by the first term in (\ref{3.20}), it is not difficult to see
that $\overline{W}(\cdot )$ and $W(\cdot )$ appeared in (5a) and (5b),
respectively, are independent. $\hfill \Box $\newline

\noindent \textbf{Proof of Theorem \ref{thm2}}. It is not difficult to show
that $\hat{\tau}_{T}$ is $T$-consistent by similar arguments in the proof of Theorem \ref{thm}. To
prove a stronger result in (\ref{result2}), we follow Appendix K in Chong (2001) with some
modifications. For $m=0,1,\cdots $ and $m/T\rightarrow 0$,
\begin{equation*}
\frac{1}{Tl(\eta _{T})}RSS_{T}(\tau _{0}-\frac{m}{T})\Rightarrow 1+h_{1}(m)
\end{equation*}%
with
\begin{equation*}
h_{1}(m)=\frac{(1-\tau _{0}+e^{2c(1-\tau _{0})}G^{2}(W,c,\tau _{0}))(1-\beta
_{2}^{2})me^{2c(1-\tau _{0})}G^{2}(W,c,\tau _{0})}{(1-\beta
_{2}^{2})me^{2c(1-\tau _{0})}G^{2}(W,c,\tau _{0})+1-\tau _{0}+e^{2c(1-\tau
_{0})}G^{2}(W,c,\tau _{0})}
\end{equation*}%
and
\begin{equation*}
\frac{1}{Tl(\eta _{T})}RSS_{T}(\tau _{0}+\frac{m}{T})\Rightarrow 1+h_{2}(m)
\end{equation*}%
with
\begin{equation*}
h_{2}(m)=\frac{(1-\beta _{2})e^{2c(1-\tau _{0})}G^{2}(W,c,\tau _{0})(1-\beta
_{2}^{2m})}{1+\beta _{2}}.
\end{equation*}%
Since both $h_{1}(m)$ and $h_{2}(m)$ are increasing functions with respect
to $m$, the first result in (\ref{result2}) can be proved. For more details,
one is referred to Chong (2001).

To show the second and third parts of (\ref{result2}), it follows easily
from the first result of (\ref{result2}) that $\hat{\beta}_{1}(\hat{\tau}%
_{T})$ and $\hat{\beta}_{1}(\tau _{0})$ have the same asymptotic
distribution, and so do $\hat{\beta}_{2}(\hat{\tau}_{T})$ and $\hat{\beta}%
_{2}(\tau _{0})$. By Lemma \ref{lem3}, we have
\begin{eqnarray*}
T(\hat{\beta}_{1}(\hat{\tau}_{T})-\beta _{1}) &\overset{d}{=} &T(\hat{\beta}%
_{1}(\tau _{0})-\beta _{1}) \\
&=&\frac{\frac{1}{Tl(\eta _{T})}\sum_{t=1}^{[\tau _{0}T]}y_{t-1}\varepsilon
_{t}}{\frac{1}{T^{2}l(\eta _{T})}\sum_{t=1}^{[\tau _{0}T]}y_{t-1}^{2}} \\
&\Rightarrow &\frac{\frac{1}{2}e^{2c(1-\tau _{0})}G^{2}(W,c,\tau
_{0})+c\int_{0}^{\tau _{0}}e^{2c(1-t)}G^{2}(W,c,t)dt-\frac{\tau _{0}}{2}}{%
\int_{0}^{\tau _{0}}e^{2c(1-t)}G^{2}(W,c,t)dt}.
\end{eqnarray*}%
Similarly, it follows from Lemma \ref{lem4} that
\begin{eqnarray*}
\sqrt{T}(\hat{\beta}_{2}(\hat{\tau}_{T})-\beta _{2}) &\overset{d}{=} &\sqrt{T%
}(\hat{\beta}_{2}(\tau _{0})-\beta _{2}) \\
&=&\frac{\frac{1}{\sqrt{T}l(\eta _{T})}\sum_{t=[\tau
_{0}T]+1}^{T}y_{t-1}\varepsilon _{t}}{\frac{1}{Tl(\eta _{T})}\sum_{t=[\tau
_{0}T]+1}^{T}y_{t-1}^{2}} \\
&\Rightarrow &\frac{\sqrt{1-\beta _{2}^{2}}\cdot \overline{W}(B(c,\tau _{0}))%
}{1-\tau _{0}+e^{2c(1-\tau _{0})}G^{2}(W,c,\tau _{0})}.
\end{eqnarray*}%
%
%
%
%
%
%
%

To derive the limiting distribution of $\hat{\tau}_{T}$ for shrinking break,
we let $\beta _{2}=\beta _{2T}=\beta _{1T}-1/\sqrt{Tg(T)}$, where $g(T)>0$
with $g(T)\rightarrow \infty $ and $g(T)/\sqrt{T}\rightarrow 0$. Moreover,
let $\nu $ be a constant. For $\tau =\tau _{0}+\nu g(T)/T$ and $\nu \leq 0$,
following Appendix K in Chong (2001), we have
\begin{equation*}
\frac{RSS_{T}(\tau )-RSS_{T}(\tau _{0})}{l(\eta _{T})}\Rightarrow
-2e^{c(1-\tau _{0})}G(W_{1},c,\tau _{0})W_{1}(|\nu |)+|\nu |e^{2c(1-\tau
_{0})}G^{2}(W_{1},c,\tau _{0}).
\end{equation*}%
Similarly, for $\tau =\tau _{0}+\nu g(T)/T$ and $\nu >0$, we have
\begin{equation*}
\frac{RSS_{T}(\tau )-RSS_{T}(\tau _{0})}{l(\eta _{T})}\Rightarrow
-2e^{c(1-\tau _{0})}G(W_{1},c,\tau _{0})W_{2}(\nu )+\nu e^{2c(1-\tau
_{0})}G^{2}(W_{1},c,\tau _{0}).
\end{equation*}%
For more details, one is referred to Chong (2001). Thus, by applying the
continuous mapping theorem for argmax functionals, we have
\begin{eqnarray*}
&&(\beta _{1T}-\beta _{2T})^{2}T^{2}(\hat{\tau}_{T}-\tau _{0}) \\
&=&\frac{T}{g(T)}(\hat{\tau}_{T}-\tau _{0})=\hat{\nu} \\
&=&\mathop{\arg\min}_{\nu \in R}\big(RSS_{T}(\tau )-RSS_{T}(\tau _{0})\big)
\\
&=&\mathop{\arg\min}_{\nu \in R}\Big(\frac{RSS_{T}(\tau )-RSS_{T}(\tau _{0})%
}{l(\eta _{T})}\Big) \\
&\Rightarrow &\mathop{\arg\min}_{\nu \in R}\bigg\{-2e^{2c(1-\tau
_{0})}G^{2}(W_{1},c,\tau _{0})\Big(\frac{B^{\ast }(\nu )}{e^{c(1-\tau
_{0})}G(W_{1},c,\tau _{0})}-\frac{|\nu |}{2}\Big)\bigg\} \\
&=&\mathop{\arg\max}_{\nu \in R}\bigg\{\frac{B^{\ast }(\nu )}{e^{c(1-\tau
_{0})}G(W_{1},c,\tau _{0})}-\frac{|\nu |}{2}\bigg\},
\end{eqnarray*}%
where $B^{\ast }(\nu )$ is a two-sided Brownian motion on $R$ defined to be $%
B^{\ast }(\nu )=W_{1}(-\nu )$ for $\nu \leq 0$ and $B^{\ast }(\nu
)=W_{2}(\nu )$ for $\nu >0$. This completes our proof. $\hfill \Box $\newline

\end{document}